\documentclass[11pt]{amsart}

\usepackage{eufrak}
\usepackage{euscript}
\usepackage{amsmath}
\usepackage{amsfonts}
\usepackage{amssymb}

\usepackage{amsxtra}
\usepackage{amscd}
\usepackage{amsthm}

\usepackage{graphicx}
\usepackage{latexsym}

\newtheorem{Proposition}{Proposition}[subsection]
\newtheorem{Lemma}{Lemma}[subsection]
\newtheorem{Theorem}{Theorem}[subsection]
\newtheorem{Corollary}{Corollary}[subsection]
\newtheorem{Remark}{Remark}[subsection]
\newtheorem{Example}{Example}[subsection]

\newcommand{\nc}{\newcommand}
\newcommand{\on}{\operatorname}

\nc{\BA}{{\mathbb{A}}}
\nc{\BC}{{\mathbb{C}}}
\nc{\BE}{{\mathbb{E}}}
\nc{\BK}{{\mathbb{K}}}
\nc{\BM}{{\mathbb{M}}}
\nc{\BN}{{\mathbb{N}}}
\nc{\BP}{{\mathbb{P}}}
\nc{\BR}{{\mathbb{R}}}
\nc{\BZ}{{\mathbb{Z}}}
\nc{\BS}{{\mathbb{S}}}

\nc{\CA}{{\mathcal{A}}}
\nc{\CB}{{\mathcal{B}}}
\nc{\CC}{{\mathcal{C}}}
\nc{\CK}{{\mathcal{K}}}
\nc{\CE}{{\mathcal{E}}}
\nc{\CF}{{\mathcal{F}}}
\nc{\CG}{{\mathcal{G}}}
\nc{\CI}{{\mathcal{I}}}
\nc{\CJ}{{\mathcal{J}}}
\nc{\CL}{{\mathcal{L}}}
\nc{\CM}{{\mathcal{M}}}
\nc{\CN}{{\mathcal{N}}}
\nc{\CO}{{\mathcal{O}}}
\nc{\CP}{{\mathcal{P}}}
\nc{\CQ}{{\mathcal{Q}}}
\nc{\CR}{{\mathcal{R}}}
\nc{\CS}{{\mathcal{S}}}
\nc{\CT}{{\mathcal{T}}}
\nc{\CU}{{\mathcal{U}}}
\nc{\CV}{{\mathcal{V}}}
\nc{\CW}{{\mathcal{W}}}
\nc{\CZ}{{\mathcal{Z}}}

\nc{\fa}{{\mathfrak{a}}}
\nc{\fb}{{\mathfrak{b}}}
\nc{\fg}{{\mathfrak{g}}}
\nc{\fh}{{\mathfrak{h}}}
\nc{\fj}{{\mathfrak{j}}}
\nc{\fk}{{\mathfrak{k}}}
\nc{\fl}{{\mathfrak{l}}}
\nc{\fm}{{\mathfrak{m}}}
\nc{\fn}{{\mathfrak{n}}}
\nc{\fo}{{\mathfrak{o}}}
\nc{\fu}{{\mathfrak{u}}}
\nc{\fp}{{\mathfrak{p}}}
\nc{\fr}{{\mathfrak{r}}}
\nc{\fs}{{\mathfrak{s}}}
\nc{\ft}{{\mathfrak{t}}}

\nc{\fA}{{\mathfrak{A}}}
\nc{\fB}{{\mathfrak{B}}}
\nc{\fD}{{\mathfrak{D}}}
\nc{\fE}{{\mathfrak{E}}}
\nc{\fF}{{\mathfrak{F}}}
\nc{\fG}{{\mathfrak{G}}}
\nc{\fK}{{\mathfrak{K}}}
\nc{\fL}{{\mathfrak{L}}}
\nc{\fM}{{\mathfrak{M}}}
\nc{\fN}{{\mathfrak{N}}}
\nc{\fP}{{\mathfrak{P}}}
\nc{\fU}{{\mathfrak{U}}}
\nc{\fV}{{\mathfrak{V}}}
\nc{\fZ}{{\mathfrak{Z}}}

\nc{\bb}{{\mathbf{b}}}
\nc{\bc}{{\mathbf{c}}}
\nc{\bd}{{\mathbf{d}}}
\nc{\be}{{\mathbf{e}}}
\nc{\bj}{{\mathbf{j}}}
\nc{\bn}{{\mathbf{n}}}
\nc{\bp}{{\mathbf{p}}}
\nc{\bq}{{\mathbf{q}}}
\nc{\bu}{{\mathbf{u}}}
\nc{\bv}{{\mathbf{v}}}
\nc{\bx}{{\mathbf{x}}}
\nc{\bs}{{\mathbf{s}}}
\nc{\by}{{\mathbf{y}}}
\nc{\bw}{{\mathbf{w}}}
\nc{\bA}{{\mathbf{A}}}
\nc{\bK}{{\mathbf{K}}}
\nc{\bB}{{\mathbf{B}}}
\nc{\bC}{{\mathbf{C}}}
\nc{\bD}{{\mathbf{D}}}
\nc{\bH}{{\mathbf{H}}}
\nc{\bM}{{\mathbf{M}}}
\nc{\bN}{{\mathbf{N}}}
\nc{\bP}{{\mathbf{P}}}
\nc{\bV}{{\mathbf{V}}}
\nc{\bW}{{\mathbf{W}}}
\nc{\bX}{{\mathbf{X}}}
\nc{\bZ}{{\mathbf{Z}}}
\nc{\bS}{{\mathbf{S}}}

\nc{\sA}{{\mathsf{A}}}
\nc{\sB}{{\mathsf{B}}}
\nc{\sC}{{\mathsf{C}}}
\nc{\sD}{{\mathsf{D}}}
\nc{\sF}{{\mathsf{F}}}
\nc{\sK}{{\mathsf{K}}}
\nc{\sM}{{\mathsf{M}}}
\nc{\sO}{{\mathsf{O}}}
\nc{\sS}{{\mathsf{S}}}
\nc{\sZ}{{\mathsf{Z}}}

\nc{\tZ}{{\tilde{\mathbb Z}}}

\nc{\bigimlat}{{\Sigma(\cowt)}}
\nc{\imlat}{{\sigma(\cowt)}}

\nc{\Aut}{\on{Aut}}

\nc{\Inn}{\on{Inn}}

\nc{\red}{{\on{red}}}

\nc{\loc}{{\on{loc}}}

\nc{\pos}{{\on{pos}}}

\nc{\ch}{\check}

\nc{\R}{\Bbb R}
\nc{\C}{\Bbb C}
\renewcommand{\H}{\on{H}}

\newcommand{\Ob}{\on{Ob}}
\newcommand{\Spec}{\on{Spec}}
\newcommand{\Rep}{\on{Rep}}
\nc{\Hom}{\on{Hom}}
\nc{\Comod}{\on{Comod}}
\nc{\GL}{\on{GL}}
\nc{\Gr}{\on{Gr}}
\newcommand{\affgr}{\on{Gr}}
\newcommand{\raffgr}{\on{Gr}_{\Bbb R}}

\nc{\eH}{{\EuScript H}}
\nc{\eG}{{\EuScript G}}

\nc{\chG}{{\check G}}
\nc{\chT}{{\check T}}

\nc{\sP}{{\mathsf{P}}}
\nc{\sQ}{{\mathsf{Q}}}
\nc{\Vect}{{\on{\mathsf Vect}}}

\nc{\risom}{\stackrel{\sim}{\rightarrow}}
\nc{\lisom}{\stackrel{\sim}{\leftarrow}}

\nc{\mv}{S^\nu}
\nc{\rmv}{{S}^\nu_\R}

\nc{\nadgr}{\mathsf{adGr}^{(n)}}
\nc{\adgr}{\mathsf{adGr}}

\nc{\twt}{{\on F}^\nu}
\nc{\swt}{{\on F}^\nu_\R}

\nc{\Lie}{\on{Lie}}

\newcommand{\taut}{\on {Aut}^\otimes}
\newcommand{\tangp}{{\taut(\Bbb H)}}

\newcommand{\liesl}{\EuFrak {sl}}

\newcommand{\atwo}{\EuFrak {a}_2}

\newcommand{\encurve}{\hat{X}}

\newcommand{\M}{\mathcal M}
\newcommand{\shL}{\mathcal L}

\newcommand{\base}{X^{(\omega)}_\R}
\newcommand{\ls}{\EuFrak L}
\newcommand{\id}{\mbox{Id}}

\newcommand{\Z}{{\Bbb Z}}

\newcommand{\hc}{{\Bbb H}}

\newcommand{\strat}{{\mathcal S}}

\newcommand{\pR}{{{}^p\hspace{-0.2em}\on{R}}}
\newcommand{\pH}{{{}^p\hspace{-0.2em}\on{H}}}
\newcommand{\ic}{\on{IC}}

\newcommand{\vect}{\on{\sf{Vect}}}

\newcommand{\ff}{\on F}
\newcommand{\ef}{\on e}

\newcommand{\GR}{G_{{\Bbb R}}}

\newcommand{\LG}{\check G}

\newcommand{\jet}{{\mathcal J}}

\newcommand{\catp}{\mathbf P}
\newcommand{\catd}{{\mathbf D}}

\newcommand{\edcat}{{\mathbf D}_{G({\mathcal O})}({\Gr})}
\newcommand{\epcat}{{\mathbf P}_{G({\mathcal O})}({\Gr})}

\newcommand{\redcat}{{\mathbf D}_{G_\R({\mathcal O}_\R)}({\Gr}_\R)}
\newcommand{\repcat}{{\mathbf P}_{G_\R({\mathcal O}_\R)}({\Gr}^+_\R)}

\newcommand{\catq}{{\mathbf Q}({\Gr}_\R)}

\newcommand{\tors}{\mathcal F}
\newcommand{\sh}{\mathcal F}

\newcommand{\naffgr}{\Gr^{(n)}}

\newcommand{\wraffgr}{{\Gr}_\R^{(\omega)}}

\newcommand{\sraffgr}{{\Gr}^{(\sigma)}_\R}
\newcommand{\stsraffgr}{{\Gr}^{(\sigma\times\sigma)}_\R}

\newcommand{\oneaffgr}{{\affgr}^{(1)}}
\newcommand{\oneraffgr}{{\affgr}^{(1)}_\R}

\newcommand{\twoaffgr}{{\Gr}^{(2)}}
\newcommand{\tworaffgr}{{\Gr}^{(2)}_\R}

\newcommand{\alg}{A}

\newcommand{\morph}{\mbox{Hom}}
\newcommand{\spec}{\mbox{Spec}}

\newcommand{\aut}{\mbox{Aut}}

\newcommand{\psh}{\mathcal P}

\newcommand{\repg}{\on{Rep}(\LG)}

\newcommand{\tchar}{\on {Ch}}

\newcommand{\co}{\mathcal{O}}
\newcommand{\ck}{\mathcal{K}}

\newcommand{\cor}{\mathcal{O}_\R}
\newcommand{\ckr}{\mathcal{K}_\R}

\newcommand{\cowt}{{\Lambda_T}}
\newcommand{\wt}{{\ch\Lambda_T}}
\newcommand{\dcowt}{{\Lambda^+_T}}

\newcommand{\rcowt}{{\Lambda_S}}
\newcommand{\rdcowt}{{\Lambda^+_S}}

\newcommand{\gstratalta}{{\Gr}^{\lambda_1}}
\newcommand{\gstrataltb}{{\Gr}^{\lambda_2}}

\newcommand{\grsc}{G^{(2)}(\co)}
\newcommand{\rgrsc}{G^{(\sigma)}_\R(\cor)}

\newcommand{\schar}{\on{Ch}_\R}

\newcommand{\rep}{\on{Rep}}

\newcommand{\lieg}{\EuFrak g}
\newcommand{\lielg}{\check{\EuFrak g}}
\newcommand{\lieh}{{\EuFrak h}}
\newcommand{\liegr}{{\EuFrak g}_\R}

\newcommand{\taffgr}{\Gr^{(2)}}

\newcommand{\xaffgr}{{}_x\hspace{-0.2em}\affgr}
\newcommand{\xraffgr}{{}_x\hspace{-0.2em}\raffgr}

\begin{document}



\title{Perverse Sheaves on Real Loop Grassmannians}

\author{\sc David Nadler }
\address{Department of Mathematics\\ 
University of Chicago\\ Chicago, IL 60637}
\email{nadler@math.uchicago.edu}

\begin{abstract}
The aim of this paper is to identify
a certain tensor category of perverse sheaves on 
the loop Grassmannian $\raffgr$ of a real form $G_\R$ of a 
connected reductive
complex algebraic group $G$ with
the category of finite-dimensional representations
of a connected reductive complex algebraic subgroup $\ch H$ of the dual group $\ch G$.
The root system of $\ch H$ is closely related to the restricted root system
of $G_\R$.
The fact that $\ch H$ is reductive 
implies that an interesting family of real algebraic
maps 
satisfies the conclusion of the Decomposition Theorem
of Beilinson-Bernstein-Deligne. 
\end{abstract}

\maketitle





\begin{section}{Introduction}\label{sintro}

It is a general principle that the representation theory
of a connected reductive complex algebraic group $G$ is reflected in the geometry of
its dual group $\LG$. Although it is simple to define $\LG$
-- it is the 
reductive group
with based root datum dual to the based root datum of $G$ --
the duality is nevertheless mysterious.
One way to concretely obtain $\LG$ from $G$ 
is to study perverse sheaves on the loop Grassmannian $\affgr$ of $G$.
A certain tensor category of perverse sheaves on $\affgr$ 
is equivalent to the category $\repg$ of finite-dimensional representations
of $\LG$~\cite{Lu83,Ginz96,BD,MV00,MV04}. (See \cite[Theorem 7.3]{MV04} 
for a final account.)
The result is fundamental in the geometric Langlands
program~\cite{BD,FGV01,BG02}.
It also leads to a construction of canonical bases~\cite{MV00, BG01, MV04},
and a deeper understanding of the Satake isomorphism~\cite{Gait01}, among
other applications~\cite{BFGM02}.  
It also may be interpreted as providing 
a down-to-earth perspective on the duality itself.
For example, according to this approach, 
the dual group of $\GL(V)$ is canonically $\GL(\H^*({\Bbb P}(V)))$,
where $\Bbb P(V)$ denotes the space of lines in $V$, 
and $\H^*(\Bbb P(V))$ its cohomology.

In this paper,
we study perverse sheaves on the loop Grassmannian
$\raffgr$ of a real form $G_\R$ of $G$.
Although the theory of perverse sheaves as developed in \cite{BBD82}
is primarily for complex algebraic spaces, it is possible
to consider certain perverse sheaves on $\raffgr$.
The reason is that the natural finite-dimensional stratification of $\raffgr$ 
is real even-codimensional.
It turns out that some statements about sheaves on 
complex algebraic spaces -- the frequent
vanishing of odd-dimensional
intersection Betti numbers, the Decomposition Theorem~\cite[Th\'eor\`eme 6.2.5]{BBD82} for pushforwards
of intersection cohomology sheaves -- also hold for certain sheaves
on $\raffgr$
though it is only real algebraic.   
The topological results of this paper 
are all consequences of the main result and its proof.
It states that a certain category of 
perverse sheaves on $\raffgr$ is a tensor category equivalent to the category $\Rep(\ch H)$
of finite-dimensional
representations of a connected reductive complex algebraic
subgroup $\ch H\subset \LG$.

In the remainder of the introduction,
we discuss the associated subgroup $\ch H\subset\ch G$,
and then sketch some of the geometry
involved in its construction. 


\begin{subsection}{The associated subgroup}
In general, the subgroup $\ch H\subset \ch G$ associated
to a real form $G_\R$ of $G$
may be realized as the identity component of the fixed points of an involution 
 of a certain Levi subgroup $\ch L_1\subset\ch G$.
The description simplifies in the following circumstances.
A real form $G_\R$ is called {\em quasi-split} if the complexification
$P\subset G$ of a minimal parabolic subgroup $P_\R\subset G_\R$
is a Borel subgroup of $G$. There is a canonical bijection
from the set of conjugacy classes of quasi-split real forms to 
the set of involutions of 
the based root datum of $G$, and so also to
the set of involutions of the based root datum of $\ch G$. 
Suppose that the real form $G_\R$ is quasi-split.
Then the Levi subgroup $\ch L_1$ is the entire dual group $\ch G$.
Suppose in addition that the involution of the root datum of $\ch G$
associated to the conjugacy class
of $G_\R$ has the property
that it fixes a node in each component of the Dynkin diagram of $\ch G$ that
it preserves. Then the subgroup $\ch H \subset\ch G$ is the identity
component of the fixed points
of a lift to $\ch G$ of the involution
of the based root datum of $\ch G$ associated
to the conjugacy class of $G_\R$. 

Section~\ref{sfpofauto}
contains a concrete description of the subgroup $\ch H\subset \ch G$ in general.
The discussion there is self-contained,
and the interested reader may consult it independently.

We list here some basic properties of the subgroup $\ch H\subset\ch G$.
The duality between the groups $G$ and $\ch G$
descends to a duality between
their Lie algebras $\lieg$ and $\lielg$. 
An isogeny of groups $G^1\to G^2$ with Lie algebra $\lieg$
is dual to an isogeny $\LG^2\to \LG^1$ of groups with Lie algebra $\lielg$. 
Similarly, the association of the subgroup $\ch H\subset \LG$ 
to the real form $G_\R$ descends to an association of the Lie
subalgebra $\ch \lieh\subset\lielg$ to the real form $\liegr$.
An isogeny $G^1\to G^2$ which commutes with conjugation 
leads to an isogeny $\LG^2\to \LG^1$ which restricts to an isogeny
$\ch H^2\to \ch H^1$.
In addition, the subgroup $\ch H\subset\ch G$ 
associated to a product of real forms 
is the product of the associated subgroups.

It is possible to read off some invariants of the root system of $\ch \lieh$ 
directly from the restricted root system of $\liegr$.
For example, the rank of $\ch\lieh$ is equal to the real rank of $\liegr$,
and the Weyl group of $\ch\lieh$ is isomorphic to the small Weyl group
of $\liegr$.
We refer the reader to Table~1 for a list of the associated
Lie algebras $\ch\lieh$ for non-compact real forms $\liegr$ 
with simple complexification $\lieg$.
It is worth pointing out that when $\lieg$ is simple, $\ch\lieh$ is simple as well.

When the real form $G_\R$ is compact, 
the real loop Grassmannian $\raffgr$ is a single point,
and the associated subgroup
is the identity $\ch H=\langle 1 \rangle\subset \ch G$.
A real form $G_\R$ is called {\em split} if there is a
torus $T_\R\subset G_\R$ isomorphic to a product of copies of $\R^\times$
such that its complexification $T\subset G$ is a maximal torus of $G$.
When the real form $G_\R$ is split,
the associated subgroup is the entire dual group $\ch H=\ch G$.
In this case, 
the fact that $\ch H$ coincides with $\ch G$
implies interesting topological statements.   
The case of a quasi-split real form $G_\R$ is also noteworthy for topological reasons,
and we shall return to it shortly.


\begin{table}[h!]\centering
$
\begin{array}{|l|llll|l|} 
\hline
& \liegr & \lieg & \ch\fg & \ch\fh & \mbox{Remarks}\\ 
\hline

\mbox{AI} & 
\EuFrak{sl}_n(\R) & \EuFrak{sl}_n(\C) & 
\EuFrak{sl}_n(\C) & \EuFrak{sl}_n(\C) & \mbox{split}\\

\mbox{AII} & 
\EuFrak{su}^*(2n) & \EuFrak{sl}_{2n}(\C) & 
\EuFrak{sl}_{2n}(\C) & \EuFrak{sl}_n(\C) &  \\

\mbox{AIII/AIV} & 
\EuFrak{su}(p,q)  & \EuFrak{sl}_{n}(\C) & 
\EuFrak{sl}_{n}(\C) & \EuFrak{sp}_{p}(\C) & p\leq q \\

& & & & & p+q=n \\

& & & & & \mbox{quasi-split if $q=p$} \\
& & & & & \mbox{ or $q=p+1$} \\


\mbox{BI/BII} & 
\EuFrak{so}(p,q) 
& \EuFrak{so}_{2n+1}(\C) & 
\EuFrak{sp}_{n}(\C) & \EuFrak{sp}_{p}(\C) & p< q \\
 
& & & & & p+q=2n+1 \\

& & & & & \mbox{split if $q=p+1$} \\


\mbox{CI} & 
\EuFrak{sp}_n(\R)  & \EuFrak{sp}_{n}(\C) & 
\EuFrak{so}_{2n+1}(\C) & \EuFrak{so}_{2n+1}(\C) & \mbox{split}\\

\mbox{CII} & 
\EuFrak{sp}(p,q) & \EuFrak{sp}_{n}(\C) & 
\EuFrak{so}_{2n+1}(\C) & \EuFrak{sp}_{p}(\C) & p\leq q\\

& & & & & p+q=n \\


\mbox{DI/DII} & 
\EuFrak{so}(n,n) & \EuFrak{so}_{2n}(\C) & 
\EuFrak{so}_{2n}(\C) & \EuFrak{so}_{2n}(\C) & \mbox{split}\\

& \EuFrak{so}(p,q)  & \EuFrak{so}_{2n}(\C) & 
\EuFrak{so}_{2n}(\C) & \EuFrak{so}_{2p+1}(\C) & p<q \\

& & & & & p+q=2n \\

& & & & & \mbox{quasi-split if $q=p+2$} \\

\mbox{DIII} & 
\EuFrak{so}^*(2n)  & \EuFrak{so}_{2n}(\C) & 
\EuFrak{so}_{2n}(\C) & \EuFrak{sp}_{p}(\C) & p=[n/2] \\


\mbox{EI} & 
\EuFrak{e}_{6(6)}  & \EuFrak{e}_{6}(\C) & 
\EuFrak{e}_{6}(\C) & \EuFrak{e}_{6}(\C) & \mbox{split} \\

\mbox{EII} & 
\EuFrak{e}_{6(2)}  & \EuFrak{e}_{6}(\C) & 
\EuFrak{e}_{6}(\C) & \EuFrak{f}_{4}(\C) & \mbox{quasi-split}\\

\mbox{EIII} & 
\EuFrak{e}_{6(-14)}  & \EuFrak{e}_{6}(\C) & 
\EuFrak{e}_{6}(\C) & \EuFrak{so}_{5}(\C) & \\

\mbox{EIV} & 
\EuFrak{e}_{6(-26)}  & \EuFrak{e}_{6}(\C) & 
\EuFrak{e}_{6}(\C) & \EuFrak{sl}_{3}(\C) & \\

\mbox{EV} & 
\EuFrak{e}_{7(7)}  & \EuFrak{e}_{7}(\C) & 
\EuFrak{e}_{7}(\C) & \EuFrak{e}_{7}(\C) & \mbox{split} \\

\mbox{EVI} & 
\EuFrak{e}_{7(-5)}  & \EuFrak{e}_{7}(\C) & 
\EuFrak{e}_{7}(\C) & \EuFrak{f}_{4}(\C) & \\

\mbox{EVII} & 
\EuFrak{e}_{7(-25)}  & \EuFrak{e}_{7}(\C) & 
\EuFrak{e}_{7}(\C) & \EuFrak{sp}_{3}(\C) & \\

\mbox{EVIII} & 
\EuFrak{e}_{8(8)}  & \EuFrak{e}_{8}(\C) & 
\EuFrak{e}_{8}(\C) & \EuFrak{e}_{8}(\C) & \mbox{split} \\

\mbox{EIX} & 
\EuFrak{e}_{8(-24)}  & \EuFrak{e}_{8}(\C) & 
\EuFrak{e}_{8}(\C) & \EuFrak{f}_{4}(\C) & \\

\mbox{FI} & 
\EuFrak{f}_{4(4)}  & \EuFrak{f}_{4}(\C) & 
\EuFrak{f}_{4}(\C) & \EuFrak{f}_{4}(\C) & \mbox{split} \\

\mbox{FII} & 
\EuFrak{f}_{4(-20)}  & \EuFrak{f}_{4}(\C) & 
\EuFrak{f}_{4}(\C) & \EuFrak{sl}_{2}(\C) & \\

\mbox{G} & 
\EuFrak{g}_{2(2)}  & \EuFrak{g}_{2}(\C) & 
\EuFrak{g}_{2}(\C) & \EuFrak{g}_{2}(\C) & \mbox{split} \\

\hline
\end{array}
$
\caption{Associated Lie algebras $\ch\fh$ for 
non-compact real Lie algebras $\liegr$ 
with simple complexifications $\lieg$.
Notation following \' E. Cartan, 
and~\cite{Helg78}.}
\end{table}


\end{subsection}

\begin{subsection}{Sketch of geometry}


%

%

Let $\CK=\C((t))$ be the field of formal Laurent series,
and let $\CO=\C[[t]]$ be the ring of formal power series.
Let $G(\CK)$
be the group of $\CK$-valued points of $G$, and let $G(\CO)$ be the group of $\CO$-valued points.
The quotient set $G(\CK)/G(\CO)$ is the $\BC$-points
of a (not necessarily reduced) ind-finite type complex algebraic ind-scheme. 
In this paper, we shall only be interested in the space of $\BC$-points of this ind-scheme 
equipped with its classical topology.
We call this space the loop Grassmannian of $G$,
and denote it by $\affgr$.
The filtration by order of pole exhibits $\affgr$ as 
an increasing union of projective varieties.
As a topological space, $\affgr$ is homeomorphic
to the space of based loops from a circle $\Bbb S^1$ to a compact form $G_c$
of $G$ whose Fourier expansions are polynomial. It
is homotopy equivalent to the space of continuous
based loops from $\Bbb S^1$ to $G$.
It is called a Grassmannian since it may also be realized as
a certain collection of subspaces in the infinite-dimensional
$\C$-vector space $\lieg(\ck)$. 
See~\cite[Section 11]{Lu83} and~\cite[Chapter 8]{PS86} for more details. 
Although the latter work in the context of polynomial loops,
the polynomial version of the loop Grassmannian is isomorphic to $\affgr$.


The action of $G(\co)$ on $\affgr$ by left-multiplication
refines the filtration by order of pole.
The orbits of the action are parameterized by the dominant
weights $\lambda$ of the dual group $\LG$. 
Lusztig~\cite[Section 11, c)]{Lu83} discovered that the local invariants of the intersection cohomology
sheaf $\ic^\lambda$ of an orbit closure 
coincide with the weight multiplicities of the corresponding 
irreducible representation $V_\lambda$
of $\LG$.
In particular, he showed that 
the dimension $\dim\hc(\affgr,\ic^\lambda)$ of the hypercohomology 
of the intersection cohomology sheaf is equal to the dimension $\dim V_\lambda$
of the corresponding representation.
Furthermore, he showed that the decomposition of the convolution
of such sheaves agrees with the decomposition of the tensor product
of irreducible representations. 

Thanks to further work by Ginzburg~\cite{Ginz96}, Beilinson-Drinfeld~\cite{BD}, 
and Mirkovi\' c-Vilonen \cite{MV00, MV04}, we have the following.

\begin{Theorem}[\cite{MV04}, Theorem 7.3, Proposition 6.3]
The category $\epcat$ of $G(\co)$-equivariant perverse sheaves
on the loop Grassmannian $\affgr$ of a connected reductive complex algebraic group
$G$ is a tensor category
equivalent 
to the category $\rep(\chG)$ of finite-dimensional representations
of the dual group $\chG$. 
Under the equivalence, the hypercohomology
of a perverse sheaf corresponds to the underlying vector space
of a representation.
\end{Theorem}

\begin{Remark}
In \cite[Theorem 12.1]{MV04}, it is shown that for
any commutative, unital, Noetherian ring $R$ of finite global dimension, the category
of $G(\co)$-equivariant perverse sheaves
with $R$-coefficients  on $\affgr$ is equivalent to 
the category of finite-dimensional $R$-representations
of the canonical smooth, split, reductive group scheme $\chG_R$ over $R$
whose root datum is dual to that of $G$.
In this paper, we work only with sheaves with complex coefficients.
\end{Remark}

We mention here two main ingredients in the proof of the theorem.
First, using the Beilinson-Drinfeld Grassmannian of $G$ over a curve $X$,  
it is possible to reformulate the convolution of perverse sheaves
in a form which is clearly commutative~\cite[Section 5]{MV04}.
Second, using 
the perverse cells of Mirkovi\'c-Vilonen, it is possible to 
decompose the hypercohomology of a perverse
sheaf into weight spaces~\cite[Theorem 3.6]{MV04}.
The ideas involved in these two constructions are 
essential to the arguments of this paper.
In an appendix, we explain why in the theorem one obtains the dual group.

It is also worth mentioning that the paper~\cite{MV04} cites this one
for two technical points. Lest there appear to be a logical loop,
we take a moment to comment on this. First, in Section~\ref{sprel}, we discuss 
the formalism we have adopted in order to work with sheaves on
infinite-dimensional spaces such as the loop Grassmannian $\affgr$.
The discussion is quite general and depends on no other results.
The paper~\cite{MV04} uses these conventions.
Second, in working with certain sheaves on the Beilinson-Drinfeld
Grassmannian over a curve $X$, it is often useful to formalize the fact that
they are constant along $X$. This is explained in Section~\ref{sconstpsh}.
When working over $\BA^1$,
it may also be accomplished by choosing a global coordinate.
This is the approach taken in \cite[Sections 5 and 6]{MV04}. They mention that to extend
their results, one could adopt the formalism explained in this paper.
In any case, this paper contains proofs of all 
assertions which are needed
but which are not explicitly in~\cite{MV04}.




Now let $G_\R$ be a real form of $G$.
The conjugation of $G$ with respect to $G_\R$
induces a conjugation of $\affgr$.
We call
the subspace of fixed points of the induced conjugation
the real loop Grassmannian of $G_\R$, and denote it by $\raffgr$. 
Let $\cor=\R[ [t] ]$ be 
the ring of
real formal power series, and 
let $\ckr=\R((t))$ be the field of
real formal Laurent series.
We may identify $\raffgr$ with the quotient
$
G_\R(\ckr)/G_\R(\cor)
$
where $G_\R(\ckr)$ is the group of $\ckr$-valued points of $G_\R$,
and $G(\cor)$ is the group of $\cor$-valued points.
Since the conjugation of $\affgr$ preserves the filtration by
order of pole, $\raffgr$ is an increasing union of real projective varieties.
It was first recognized by Quillen, then proved in~\cite[Theorems 5.1 and 5.2]{Mitch88},
that as a topological space, $\raffgr$ is
homotopy equivalent to the based loop space of the symmetric variety $G/K$
where $K$ is the complexification of a maximal compact
subgroup of $G_\R$. 
It also may be realized as a certain collection of subspaces 
in the infinite-dimensional $\R$-vector space $\liegr(\ckr)$.

Each $G_\R(\cor)$-orbit in $\raffgr$ is finite-dimensional 
and real even-codimensional in the closure of any other. Thus
the components of $\raffgr$ are of two types: those containing
only even-dimensional $G_\R(\cor)$-orbits, and those containing
only odd-dimensional $G_\R(\cor)$-orbits.
In this paper, we restrict our attention to sheaves supported on
the union $\raffgr^+$ of the components of $\raffgr$ containing even-dimensional
orbits. The topological aspects of the theory of perverse sheaves
hold for such a stratified space. For example, there is a self-dual
perversity, and the resulting category of perverse
sheaves is abelian with simple objects the intersection cohomology
sheaves of strata with irreducible local system coefficients.

The category $\repcat$ of 
$G_\R(\cor)$-equivariant perverse sheaves 
on $\raffgr^+$ is more
complicated than the category
$\epcat$ of $G(\co)$-equivariant sheaves on $\affgr$.
For example, it follows from the results of Lusztig~\cite[Section 11]{Lu83}  
that $\epcat$
is semisimple, and the hypercohomology
is an exact and faithful functor.  
Neither of these statements is true in general for $\repcat$.
In this paper, we restrict our attention to a certain
strict full subcategory $\catq$. 
Roughly speaking, we shall only consider
perverse sheaves built out of the intersection cohomology sheaves
of certain relevant strata of $\raffgr^+$ with constant coefficients.
In fact, it will turn out that we are only dealing with
semisimple perverse sheaves, but proving this was 
one of our original motivations.

To define the category $\catq$, let $X$ be a smooth complex algebraic curve with
$X_\R$ a nonempty real form of $X$.
There is a real algebraic family over $X$ 
with fiber $\affgr$ over
points of $X\setminus X_\R$, and fiber $\affgr_\R$ over points of $X_\R$.
It is a real form of the Beilinson-Drinfeld Grassmannian 
of $G$.
We define the specialization functor
$$\on{R}:\epcat\to\redcat$$
to be the nearby cycles in this family. 
The specialization takes a sheaf on $\affgr$ to a sheaf
supported on the components $\raffgr^+$.
Therefore we may define the perverse specialization functor 
$$\pR:\epcat\to\repcat$$
to be the direct sum of the perverse homology sheaves 
of the specialization
$$
\pR=\sum_k \pH^k\circ\on{R}.
$$
We take $\catq$ 
to be the strict full subcategory of $\repcat$
consisting of objects isomorphic to subquotients
of objects of the form $\pR(\psh)$, as $\psh$
runs through all objects in $\epcat$.
(For the sake of simplicity, we have ignored here the natural $\Z$-grading
on the perverse specialization $\pR$, but it plays an important role
in the identification of the category $\catq$.)

Our main result is the following.

\begin{Theorem}
The category $\catq$ is a tensor category equivalent
to the category $\Rep(\ch H)$ of finite-dimensional representations
of a connected reductive complex algebraic subgroup $\ch H\subset \ch G$.
Under the equivalence, the hypercohomology of a perverse sheaf
corresponds to the underlying vector space of a representation,
and the perverse specialization 
corresponds to
the restriction functor on representations.
\end{Theorem}

As mentioned earlier,
Section~\ref{sfpofauto}
contains a concrete description of the subgroup $\ch H\subset \ch G$.



The following topological corollary
was one of our original motivations for establishing the equivalence
of categories of the theorem.
It asserts in particular
that there is an interesting family of real algebraic maps
that satisfy the conclusion of the Decomposition Theorem of 
Beilinson-Bernstein-Deligne~\cite[Th\'eor\`eme 6.2.5]{BBD82}.

\begin{Corollary}
The category $\catq$ is semisimple.
Each object is isomorphic to a direct sum of 
intersection cohomology sheaves 
with constant coefficients.
In particular, the convolution and perverse
specialization are semisimple.
\end{Corollary}

Finally, we point out an interesting aspect of the construction of the category $\catq$.
It is a fundamental result
in the theory of perverse sheaves
that the nearby cycles in a complex algebraic family over a curve
are perverse. (See~\cite[Section 4.4]{BBD82} or \cite[Section 6.5]{GM83b}.) 
But the nearby
cycles in a real algebraic family are not in general perverse.

\begin{Theorem}
The specialization 
$$\on{R}:\epcat\to\catd_{G_\R(\cor)}(\raffgr^+)$$
is perverse
if and only if the real form $G_\R$ is quasi-split.
\end{Theorem}


\end{subsection}


\begin{subsection}{Organization of paper}
We conclude the introduction with a brief summary of
the contents of the other sections of the paper.

In Section~\ref{sprel}, we collect notation
used throughout the paper,
and then establish the formalism we have adopted in order to work with 
infinite-dimensional spaces such as the loop Grassmannian $\affgr$.
The upshot is that for our needs 
we may discuss
such spaces as if they were finite-dimensional.   
In Section~\ref{srlg}, we give a brief account of
perverse sheaves on the loop Grassmannian $\affgr$,
and then collect those results which extend easily
to the real loop Grassmannian $\raffgr$.
In Section~\ref{srbdg}, we describe
the Beilinson-Drinfeld Grassmannian over a curve
and its real forms. We then prove some basic results
concerning their natural embeddings and stratifications.
In Section~\ref{sspec}, we define specialization
functors which take sheaves on the loop Grassmannian $\affgr$
to sheaves on its real form $\raffgr$.
In Section~\ref{smon}, we prove that these specialization
functors are monoidal.
In Section~\ref{sim1}, we introduce the category $\catq$, and its
graded versions.
In Section~\ref{schar}, we recall the construction of
weight functors due to Mirkovi\' c-Vilonen. We describe their basic properties,
then apply them to understand the specialization functors.
In Section~\ref{sim2}, we show that $\catq$ is a neutral Tannakian category.
Finally, in Section~\ref{ssub}, we apply Tannakian formalism
to the category $\catq$ and identify the corresponding group $\ch H$.

We have also included an appendix sketching how one
identifies the Tannakian group of the category $\epcat$
of $G(\CO)$-equivariant perverse sheaves on $\affgr$
with the dual group $\ch G$.

\end{subsection}

\end{section}


\medskip

\noindent{\bf Acknowledgements.}
I would like to thank my thesis advisor, Robert MacPherson,
for his guidance and support.
He introduced me to the areas of mathematics discussed in this paper, 
and offered inspiring ways to think about them.
I would like to thank Jared Anderson and Kari Vilonen for
many helpful conversations.
I would also like to thank Robert Kottwitz and a referee
for suggesting I present the subgroup associated
to a real form
as the fixed points of an automorphism.
This work was undertaken at Princeton University,
and supported in part by the NSF.



\section{Preliminaries}\label{sprel}


\subsection{Notation}\label{ssnot}
Throughout this paper, $G$ will be a connected reductive complex algebraic group,
$\theta$ will be a conjugation of $G$,
and $G_\R$ will be the real form of $G$ with respect to $\theta$.

Choose once and for all a 
maximal split torus $S_\R\subset G_\R$, and
a maximal torus $T_{\R}\subset G_\R$ such that $S_\R\subset T_\R$.
By definition, the torus $S_\R$ is isomorphic to $(\R^\times)^r$,
for some $r$, and maximal among subgroups of $\GR$ 
with this property.
Let $S\subset G$ be the complexification of $S_\R$,
and let $T\subset G$ be the complexification of $T_\R$.

Choose a mimimal parabolic subgroup $P_\R\subset \GR$
such that $T_\R\subset P_\R$, and a Levi factor $M_\R\subset P_\R$ such that 
$T_\R\subset M_\R$.
Let $P\subset G$ be the complexification of $P_\R$, and let $M\subset G$
be the complexification of $M_\R$.

Choose a Borel subgroup $B \subset G$ such that $T\subset B$, and $B\subset P$.

\begin{subsubsection}{Notation for $G$}

Let $\cowt=\morph(\C^\times,T)$ be the lattice of coweights of $T$,
and let $\dcowt\subset\cowt$ be the semigroup of coweights dominant with respect to $B$.

Let $R\subset \cowt$ be the set of coroots of $G$ with respect to $T$,
let $R^\pos\subset R$ be the set of positive coroots with respect to $B$,
and let $\Delta_{B,T}\subset R^\pos$ be the set of simple coroots.

Let $Q\subset\cowt$ be the lattice generated by $R$, and let
$Q^\pos\subset Q$ be the semigroup generated by $R^\pos$.

The coweight lattice $\cowt$ is naturally ordered: $\lambda,\mu\in\cowt$
satisfy
$\mu\leq\lambda$ if and only if $\lambda-\mu$ is a non-negative integral
linear combination of positive coroots of $G$.

Let $\CW_G=N_G(T)/T$ be the Weyl group of $G$, where $N_G(T)\subset G$ 
is the normalizer of $T$. 
It acts naturally on the coweight lattice $\cowt$,
and there is a unique dominant coweight in each orbit.

Let $2\ch\rho$ be the sum of the positive roots of $G$,
and let $\langle 2\ch\rho,\lambda\rangle\in\Z$
be the natural pairing for a coweight $\lambda\in\cowt$.

\end{subsubsection}


\begin{subsubsection}{Notation for $G_\R$}

We call 
$
\rcowt=\morph(\C^\times,S)
$ 
the lattice of {real coweights}, 
and we call $\rdcowt=\rcowt\cap\dcowt$ the semigroup of {dominant real coweights}.

The real coweight lattice $\rcowt$ is naturally ordered by the restriction of the
ordering from the coweight lattice $\cowt$:
$\lambda,\mu\in\rcowt$
satisfy
$\mu\leq\lambda$ if and only if $\lambda-\mu$ is a non-negative integral
linear combination of positive coroots of $G$.

Let $\CW_{G_\R}=N_{G_\R}(S_\R)/Z_{G_\R}(S_\R)$ be the small Weyl group of $G_\R$, 
where $N_{G_\R}(S_\R)\subset G_\R$ 
is the normalizer of $S_\R$, and $Z_{G_\R}(S_\R)\subset G_\R$ 
is the centralizer of $S_\R$. 
It acts naturally on the real coweight lattice $\rcowt$,
and there is a unique dominant real coweight in each orbit.

The conjugation $\theta$ induces an involution of $\cowt$
which we also denote by $\theta$.
For a coweight $\lambda:\C^\times\to T$, the coweight $\theta(\lambda)$ is
defined to be the composition
$$
\theta(\lambda):\C^\times\stackrel{c}{\to}\C^\times\stackrel{\lambda}{\to}
T\stackrel{\theta}{\to} T
$$
where $c$ is the standard conjugation of $\C^\times$ with respect to $\R^\times$.

The real coweight lattice $\rcowt$ is the fixed points
of the involution $\theta$ of the coweight lattice $\cowt$.
Therefore we have the projection
$$
\sigma:\cowt\to\rcowt
$$
$$
\sigma(\lambda)=\theta(\lambda)+\lambda,
$$
whose image we denote by
$\imlat\subset\rcowt.$

Let $2\ch\rho_M$ be the sum of the positive roots of the Levi factor $M\subset P$.
We have the projection
$$
\Sigma:\cowt\to\rcowt\times\Z
$$
$$
\Sigma(\lambda)=(\theta(\lambda)+\lambda,\langle 2\ch\rho_M,\lambda\rangle),
$$
whose image we denote by $\bigimlat\subset\rcowt\times\Z.$

The following is a useful characterization of the ordering on the real coweight lattice $\rcowt$.

\begin{Lemma}\label{lrcowtord}
The intersection $\rcowt\cap Q^\pos$ is generated by the elements 
$\alpha\in\rcowt\cap R^\pos$, and $\sigma(\alpha)\in\rcowt$,
for $\alpha\in R^\pos$.
\end{Lemma}

\begin{proof}
Since $\rcowt$ pairs trivially with $2\ch\rho_M$,
we may write $\beta\in\rcowt\cap Q^\pos$ 
uniquely as a sum $\beta=\sum_i \alpha_i$ of simple coroots $\alpha_i\in\Delta_0$
of the Levi subgroup $L_0\subset G$ 
that centralizes $2\ch\rho_M$.
Since $\theta$ preserves the set $\Delta_0$, and the sum $\beta=\sum_i \alpha_i$ 
is unique, $\theta$ also
preserves the set $\{\alpha_i\}$ of simple coroots counted with multiplicities
appearing in the sum.
The assertion follows by induction on the size of this set.
\end{proof}

\end{subsubsection}


\begin{subsubsection}{Graded categories}

Let $\vect$ be the category of finite-dimensional $\C$-vector spaces.

For a lattice $\Lambda$, we have the category 
$\vect_\Lambda$ of finite-dimensional
$\Lambda$-graded vector spaces. It is canonically 
equivalent to the category $\Rep(\ch S_\Lambda)$
of finite-dimensional representations of the torus
$
\ch S_\Lambda=\Spec(\C[\Lambda]).
$

For a lattice homomorphism $\tau:\Lambda_1\to\Lambda_2$, we have the functor
$$
\tau:\vect_{\Lambda_1}\to\vect_{\Lambda_2}
$$
$$
\tau(V)^\lambda=\sum_{\tau(\mu)=\lambda} V^\mu.
$$
It is canonically isomorphic to the restriction functor 
$\Rep(\ch S_{\Lambda_1})\to \Rep(\ch S_{\Lambda_2})$ coming from the group
homomorphism $\ch S_{\Lambda_2}\to \ch S_{\Lambda_1}$
induced by the ring homomorphism $\C[\Lambda_1]\to\C[\Lambda_2]$. 

More generally,
for a lattice $\Lambda$, and a $\C$-linear abelian category $\sC$, we have the 
$\C$-linear abelian category
$\sC_\Lambda=\sC\otimes\vect_\Lambda$ whose objects are $\Lambda$-graded
sums of objects of $\sC$, and whose morphisms
are $\Lambda$-graded
sums of morphisms of $\sC$.

For a lattice homomorphism $\tau:\Lambda_1\to\Lambda_2$, we have the functor
$$
\tau:\sC_{\Lambda_1}\to\sC_{\Lambda_2}
$$
$$
\tau(X)^\lambda=\sum_{\tau(\mu)=\lambda} X^\mu.
$$
When $\tau:\langle 0\rangle\to\Lambda$ is the inclusion of zero,
we obtain the fully faithful functor $\ef:\sC\to\sC_\Lambda$ that places
objects and morphisms in degree zero, and when $\tau:\Lambda\to\langle 0\rangle$
is the projection to zero, we obtain the forgetful functor $\ff:\sC_\Lambda\to\sC$
that forgets the grading on objects and morphisms.

\end{subsubsection}


\subsection{Ind-schemes}\label{ssind}

Many of the spaces we discuss in this paper are 
infinite-dimensional.
However, all of the geometry we study takes place in finite-dimensional 
approximations.
The complication is that no single finite-dimensional
approximation is sufficient. Instead, we consider
compatible families of approximations.
In this section, we describe the approach
we have adopted to formalize this.
All varieties and schemes are either real or complex,
and not assumed to be irreducible. 

Our basic object is a 
{\em compatible family} of varieties
$$
\begin{array}{cccccccc}
 \vdots & & \vdots & & \vdots & & \vdots &  \\
\downarrow & & \downarrow & & \downarrow & & \downarrow &  \\
Z^{\ell+1}_0 & \to & Z^{\ell+1}_1 &\to\cdots\to & Z^{\ell+1}_k & \to & Z^{\ell+1}_{k+1} & \to\cdots \\
\downarrow & & \downarrow & & \downarrow & & \downarrow &  \\
Z^\ell_0 & \to & Z^{\ell}_1 &\to\cdots\to & Z^\ell_k & \to & Z^\ell_{k+1} & \to\cdots \\
\downarrow & & \downarrow && \downarrow & & \downarrow &  \\
 \vdots & & \vdots & & \vdots & & \vdots &  \\
\downarrow & & \downarrow && \downarrow & & \downarrow &  \\
Z^1_0 & \to & Z_1^1& \to\cdots\to & Z^1_k & \to & Z^1_{k+1} & \to\cdots \\
\downarrow & & \downarrow && \downarrow & & \downarrow &  \\
Z^0_0 & \to & Z_1^0& \to\cdots\to & Z^0_k & \to & Z^0_{k+1} & \to\cdots \\
\end{array}
$$
satisfying the requirements:
\begin{enumerate}
\item The horizontal maps are closed embeddings.
\item The vertical maps are smooth fibrations such that
their fibers are affine spaces for $\ell$ large relative to $k$.
\item All squares in the diagram are Cartesian.
\end{enumerate}

Informally speaking, we are willing to throw away
any bounded part of the diagram, and 
to consider the resulting diagram as equivalent.
To formalize this, we consider the inverse limit 
$
Z_k=\lim_{\leftarrow} Z_k^\ell,
$
and the direct limit 
$
Z=\lim_{\rightarrow}Z_k.
$
We use the terms scheme and ind-scheme to refer to
spaces which arise in this way.

We say that $Z$ is stratified if
the varieties $Z_k^\ell$ are compatibly stratified in the sense that
for any map in the diagram, the strata of the domain
are the inverse images of the strata of the codomain.
We only consider conjugations $\theta$ and real forms $Z_\R$ 
of complex $Z$ 
which are the limits of compatible families of conjugations and
real forms. Moreover, the conjugations we consider always respect the
stratifications considered, and the strata of $Z_\R$ 
are always taken to be the real forms of the 
strata of $Z$.   

A special case of such a compatible family
is when the vertical families are constant,
and the diagram simplifies to
$$
Z_0\to Z_1\to Z_2\to\cdots\to Z_k\to \cdots.
$$
In this case, the direct limit $Z$
is an ind-scheme of ind-finite type.

Another special case is when the horizontal families
are constant, and in addition each of the varieties in the family
$$
L^0\gets L^1 \gets L^2\gets \cdots\gets L^\ell \gets\cdots 
$$
is a linear algebraic group.
In this case, the inverse limit $L$
is a group-scheme, although not usually of finite type.
We only consider group-schemes of this form such that for large $\ell_1<\ell_2$, 
the kernel of the projection $L^{\ell_1}\gets L^{\ell_2}$  is unipotent. 
We call such group-schemes {stable}.

We only consider maps between ind-schemes
which are the limit of compatible families of maps defined for $\ell$ large 
relative to $k$.
By an action of a group-scheme
$
L =\lim_{\leftarrow } L^\ell
$
on an ind-scheme
$
Z =\lim_{\rightarrow}\lim_{\leftarrow} Z^\ell_k,
$ 
we mean a compatible family of actions 
of the group-schemes $L^\ell$ on the varieties $Z_k^\ell$
defined for $\ell$ large relative to $k$.
When the vertical maps 
$Z_k^0\gets Z^{\ell}_k$
are the projections
of $L^\ell$-torsors,
we say the ind-scheme $Z$ is an $L$-torsor
over the ind-scheme $Z_0$.

We next discuss categories of sheaves on ind-schemes.
We work with sheaves of $\C$-vector spaces in the classical topology.
For a variety $Z$ acted upon by a linear algebraic group $L$, 
we write $\catd_{L}(Z)$ for the derived category of $L$-equivariant 
sheaves on $Z$ with bounded cohomology sheaves.
If $\strat$ is a locally-trivial stratification of $Z$, 
we write $\catd_{L,\strat}(Z)$
for the derived category of $L$-equivariant sheaves
on $Z$ with bounded $\strat$-constructible cohomology sheaves.
See~\cite{BL94} for the definitions of these categories.
Note that in the case when the orbits of $L$ in $Z$
coincide with the strata of $\strat$, then the forgetful
functor
$
\catd_{L,\strat}(Z)\stackrel{}{\to}\catd_{L}(Z)
$ 
is an equivalence.
We freely use the term sheaf to mean a complex of sheaves,
and $\strat$-constructible to mean with $\strat$-constructible
cohomology sheaves.

If an ind-scheme
$
Z =\lim_{\rightarrow} Z_k
$
of ind-finite type is acted upon 
by a stable group scheme $L=\lim_{\leftarrow}L^\ell$,
we may construct a directed system of categories
$$
\cdots\to  \catd_L(Z_k)  \to   \catd_L(Z_{k+1})  \to\cdots 
$$
in which each of the maps is the direct image functor.
We call the direct limit
$
\catd_L(Z)=\lim_{\rightarrow}\catd_L(Z_k)
$ 
the derived category of $L$-equivariant sheaves on $Z$.
If $Z$ is stratified, we similarly define 
the derived category $\catd_{L,\strat}(Z)$ of $L$-equivariant,
$\strat$-constructible sheaves on $Z$.

To work with sheaves on an ind-scheme
$
Z=\lim_{\rightarrow}\lim_{\leftarrow} Z^\ell_k
$
not necessarily of ind-finite type,
we assume that $Z$ is stratified.
If in addition $Z$ is acted upon by a stable group-scheme
$
L 
$
then we have a directed system of categories
$$
\cdots\to \catd_{L^\ell,\strat}(Z^\ell_k) \to \catd_{L^{\ell+1},\strat}(Z^{\ell+1}_k) \to \cdots
$$
in which each of the maps is the pull-back functor.
This sequence of categories stabilizes, because
the projections $Z_k^\ell\gets Z_k^{\ell+1}$ are compatibly stratified and
acyclic for $\ell$ large relative to $k$,
and 
the kernels of the projections $L^{\ell_1}\gets L^{\ell_2}$ are unipotent
for $\ell_1<\ell_2$ large. 
We have the limit derived category 
$
\catd_{L,\strat}(Z_k)=\lim_{\rightarrow}\catd_{L^\ell,\strat}(Z^\ell_k)
$ 
of
$L$-equivariant, $\strat$-constructible sheaves on $Z_k$.
Thanks to the stability,
we have a fully faithful directed system of categories
$$
\cdots\to  \catd_{L,\strat}(Z_k)  \to   \catd_{L,\strat}(Z_{k+1})  \to\cdots 
$$
in which each of the maps is the direct image functor.
We call the direct limit
$
\catd_{L,\strat}(Z)=\lim_{\rightarrow}\catd_{L,\strat}(Z_k)
$ 
the derived category of
$L$-equivariant, $\strat$-constructible sheaves on $Z$.

We only consider functors between such categories of sheaves which
are the limit of compatible families of functors
defined for $\ell$ large relative to $k$.
For a map 
$f: Y\to Z$ which is the limit of maps $f^\ell_k: Y^\ell_k\to Z^\ell_k$,
to obtain well-defined limit functors $f_*,f_!,f^*,f^!$,
one must require certain properties of the maps $f^\ell_k$
and the maps in the diagrams of $Y$ and $Z$.
For example, for the pushforwards $f_*,f_!$, one must have a Cartesian diagram
$$
\begin{array}{ccc}
Y^{\ell+1}_k& \stackrel{f^{\ell+1}_k}{\to}& Z^{\ell+1}_k \\
\downarrow && \downarrow\\
Y^{\ell}_k & \stackrel{f^{\ell}_k}{\to}& Z^\ell_k.\\ 
\end{array}
$$
Then by the standard identity for the composition of maps
and the smooth base change isomorphism, the limit functors $f_*,f_!$
are well-defined.
%
All of our maps satisfy the necessary requirements where we 
use such functors.
The standard identities among such functors
hold in this setting.

For an ind-scheme $Z$ of ind-finite type, there is in general no dualizing object,
but we do have the Verdier duality functor since the horizontal maps in the diagram
of $Z$ are closed embeddings.

We shall often make use of the following construction.
Let $X$ and $Y$ be stratified ind-schemes of ind-finite type,
let $L$ and $M$ be stable group-schemes, 
and let $p:Z\to X$ be a stratified $L$-torsor.
Assume that $M$ acts on $Z$ and $X$ such that $p:Z\to X$ is $M$-equivariant.
Then the pull-back functor 
$
p^*:\catd_{M,\strat}(X)\stackrel{}{\to} \catd_{M\times L,\strat}(Z)
$
is an equivalence.
If $L$ also acts on $Y$, then we 
may form 
the twisted product 
$$
X\tilde\times Y=\lim_{\rightarrow}\lim_{\leftarrow}Z_{k}^\ell\times_{L^\ell} Y_k
$$
where the variety $Z_{k}^\ell\times_{L^\ell} Y_k$ makes sense for $\ell$ large relative to $k$. 
The twisted product $X\tilde\times Y$ is stratified by the twisted products of 
the strata of $X$ and $Y$. 
In this situation, we may define a functor 
$$
\tilde\boxtimes:\catd_{M,\strat}(X)\times\catd_{L,\strat}(Y)\to\catd_{M,\strat}(X\tilde\times Y)
$$
by the requirement that
$$
q^*(\sh_1\tilde\boxtimes\sh_2)=p^*(\sh_1)\boxtimes\sh_2
$$
where $q:Z\times Y\to X\tilde\times Y$ is the natural projection.

Finally, suppose a real or complex ind-scheme 
$
Z =\lim_{\rightarrow} Z_k
$ 
of ind-finite type
is stratified such the strata 
are of even real dimension.
If in addition $Z$ is acted upon by a stable group-scheme
$
L, 
$
then we call the direct limit
$
\catp_{L,\strat}(Z)=\lim_{\rightarrow} \catp_{L,\strat}(Z_k)
$
the category 
of $L$-equivariant, $\strat$-constructible perverse sheaves on $Z$.
Here $\catp_{L,\strat}(Z_k)$ is the heart of $\catd_{L,\strat}(Z_k)$
with respect to the perverse $t$-structure 
defined by pulling back the perverse $t$-structure
from $\catd_\strat(Z_k)$ via the forgetful
functor $\catd_{L,\strat}(Z_k)\to\catd_\strat(Z_k)$.  
Since the direct image functor is $t$-exact for a closed embedding,
there is an induced $t$-structure on the limit category $\catd_{L,\strat}(Z)$,
and $\catp_{L,\strat}(Z)$ is the heart of the category $\catd_{L,\strat}(Z)$
with respect to the induced $t$-structure.
For any integer $k$, we have the functor
$
\pH^k:\catd_{L,\strat}(Z)\to\catp_{L,\strat}(Z)
$
which assigns to a sheaf its $k$-th perverse
homology sheaf. 
%



\begin{section}{Real loop Grassmannians}\label{srlg}
 
In this section, we recall some basic results about the loop Grassmannian $\affgr$
of the connected reductive complex algebraic group $G$,
then describe those results whish extend readily to the real loop Grassmannian $\raffgr$
of the real form $G_\R$.


\begin{subsection}{Loop Grassmannians}


Let $\CO=\C[ [ t ] ]$ 
be the ring of formal power series,
and for $\ell\geq 0$, let $\CO^\ell\subset\CO$ 
be the ideal generated by $t^\ell\in\CO$,
and let $\CJ^\ell=\CO/\CO^\ell$ be the finite-dimensional quotient.

The group $\Aut(\CO)$ of automorphisms of $\CO$ 
is the inverse limit
of the linear algebraic groups $\Aut(\CJ^\ell)$,
and for any $m\geq \ell>0$, the kernel of the
projection $\Aut(\CJ^{m})\to \Aut(\CJ^{\ell})$ 
is unipotent.
Let $\Aut(\CO^\ell)\subset \Aut(\CO)$ 
be the kernel 
of the projection $\Aut(\CO)\to \Aut(\CJ^\ell)$.

The group $G(\CO)$ of $\CO$-valued points of $G$
is the inverse limit
of the linear algebraic groups $G(\CJ^\ell)$,
and for any $m\geq \ell>0$, the kernel of the
projection $G(\CJ^{m})\to G(\CJ^{\ell})$ 
is unipotent.
Let $G(\CO^\ell)\subset G(\CO)$
be the kernel 
of the projection $G(\CO)\to G(\CJ^\ell)$.




Let $\CK=\C(( t ))$ 
be the field of formal Laurent series,
and 
for $k\geq 0$, let $\CK_k\subset\CK$ be the $\CO$-submodule generated by $t^{-k}\in\CK$.

The group $\GL_N(\CK)$ of $\CK$-valued points of $\GL_N$
is the direct limit of the schemes
$\GL_N(\CK)_k$ of matrices $g\in\GL_N(\CK)$ such that the
entries of $g$ and $g^{-1}$ lie in $\CK_k$.

Choose an embedding $G\subset\GL_N$.
The group $G(\CK)$ of $\CK$-valued points of $G$ is the direct limit 
of the schemes
$
G(\CK)_k
=G(\CK)\cap \GL_N(\CK)_k,
$
and the finite-dimensional varieties
$
G(\CK)^\ell_k=G(\CK)_k/G(\CO^\ell)
$
are compatibly stratified by the orbits of $G(\CO)\times G(\CO)$ acting
by multiplication on the left and right.
This realizes $G(\CK)$ as a stratified ind-scheme,
although not of ind-finite type.




The loop Grassmannian $\Gr=G(\CK)/G(\CO)$
is the
direct limit of 
the finite-dimensional projective varieties
$
\Gr_k=G(\CK)_k/G(\CO)
$
which are compatibly stratified by
the orbits of $G(\co)$ acting 
by left multiplication.
This realizes $\Gr$ as a stratified ind-scheme of ind-finite type.
%

\end{subsection}


\begin{subsection}{Spherical perverse sheaves}

To describe the $G(\CO)$-orbits in $\Gr$, observe
that each coweight $\lambda\in\cowt$ defines a point $\lambda\in\Gr$
via the embedding $\cowt\subset G(\CK)$.  
Let $\affgr^\lambda$ be the $G(\CO)$-orbit 
$G(\CO)\cdot\lambda\subset\Gr$ through $\lambda\in\Gr$.
The Cartan decomposition~\cite[Corollary 2.17]{IM65} 
states that each $G(\CO)$-orbit in $\Gr$
is of the form $\Gr^\lambda$ for some $\lambda\in\cowt$,
and two orbits $\Gr^{\lambda_1}$ and $\Gr^{\lambda_2}$ coincide if and only if 
$\lambda_1$ and $\lambda_2$ 
are in the same $\CW_G$-orbit in $\cowt$.
Note that the $G(\CO)$-orbit $\affgr^\lambda$ contains the $G$-orbit
$G\cdot\lambda$ which is isomorphic to the flag manifold $G/P^\lambda$
where the parabolic subgroup $P^\lambda\subset G$ is the stabilizer of $\lambda$.
As $z\in\BC^\times$ tends to $0$, 
the automorphism of $\affgr$ induced by the coordinate change $t\mapsto zt$
provides a retraction $\affgr^\lambda\to G/P^\lambda$.
The fibers of the retraction are the orbits of the pro-unipotent
congruence subgroup $G(\CO^1)\subset G(\CO)$.
For more details and arguments for the other assertions of the following proposition,
see for example \cite[Section 11, a) and b)]{Lu83} and \cite[Section 8.6]{PS86}.
Although the latter work in the context of polynomial loops, the polynomial
version of the loop Grassmannian is isomorphic to $\affgr$.

\begin{Proposition}\label{plgstrat}
The loop Grassmannian $\affgr$ is the disjoint union of
the $G(\co)$-orbits $\affgr^\lambda$ 
through the dominant coweights $\lambda\in\dcowt$.
The closure of the orbit $\affgr^\lambda$ is the union of the orbits 
through $\mu\in\dcowt$ with $\mu\leq\lambda$.
The orbit $\affgr^\lambda$ is a vector bundle
over the flag manifold $G/P^\lambda$, and its closure
is a projective variety of 
dimension $2\langle \ch\rho,\lambda \rangle$.
\end{Proposition}


For $k,m\geq 0$,
there exists $\ell\geq 0$ such that $G(\CO^\ell)\rtimes\Aut (\CO^\ell)$ 
acts trivially
on $G(\CK)^m_k$,
and the action of $G(\CO)\rtimes\Aut(\CO)$ 
passes to $G(\CJ^\ell)\rtimes\Aut(\CJ^\ell)$.
In particular, for $k\geq 0$,
there exists $\ell\geq 0$ such that $G(\CO^\ell)\rtimes\Aut (\CO^\ell)$ 
acts trivially
on $\Gr_k$,
and the action of $G(\CO)\rtimes\Aut(\CO)$ 
passes to $G(\CJ^\ell)\rtimes\Aut(\CJ^\ell)$.
Therefore $G(\CO)\rtimes\Aut(\CO)$ acts on 
$\Gr$ according to our conventions.

We have the category $\catp_{G(\CO)\rtimes\Aut(\CO)}(\Gr)$ 
of $G(\CO)\rtimes\Aut(\CO)$-equivariant perverse sheaves on $\Gr$,
and the category $\epcat$ of $G(\CO)$-equivariant perverse sheaves on $\Gr$.
We also have the category $\catp_{\strat}(\Gr)$ of perverse sheaves on $\Gr$
constructible with respect to the stratification by $G(\CO)$-orbits.
In the following, the first equivalence is in~\cite[Proposition 1]{Gait01}, 
and both are in~\cite[Appendix A]{MV04}. 
The proof given here shows
that for complex coefficients
the categories are in fact semisimple.

\begin{Proposition}\label{autoe}
The forgetful functors are equivalences
$$\catp_{G(\CO)\rtimes\Aut(\CO)}(\Gr)\risom\epcat
\risom\catp_{\strat}(\Gr).$$
\end{Proposition}

\begin{proof}
By Proposition~\ref{plgstrat}, each of the $G(\co)$-orbits in $\affgr$
is connected and simply-connected.
The stabilizer in $G(\co)$ of a coweight $\lambda\in\affgr$
is the parahoric subgroup $\EuScript P^\lambda$ which is connected.
The stabilizer in $G(\CO)\rtimes\Aut(\CO)$
is the semidirect product $\EuScript P^\lambda\rtimes\Aut(\CO)$ which is also connected.
Therefore each of the categories have the same simple objects,
and there are no self-extensions of simple objects.
By Proposition~\ref{plgstrat}, the $G(\co)$-orbits in
a given component of $\affgr$ 
are either all even-dimensional or all odd-dimensional.
By~\cite[Section 11, c)]{Lu83}, the stalks of  the simple objects in the categories
have the parity vanishing property: they are non-zero only in
the parity of the dimension of their support.
Therefore there are no other extensions, and we conclude that
the categories are semisimple.
\end{proof}



\begin{subsection}{Convolution}
First, recall that from the $G(\CO)$-torsor $p:G(\CK)\to\Gr$
and the $G(\CO)$-action on $\Gr$
we may construct the twisted product 
$\Gr\tilde \times\Gr$. 
To be concrete, it is the direct limit of the finite-dimensional projective varieties
$
\Gr_{k_1}\tilde\times\Gr_{k_2}=G(\CK)_{k_1}\times_{G(\CO)}\Gr_{k_2}
$
which are compatibly stratified by the twisted product of strata
$
\affgr^{\lambda_1}\tilde\times\affgr^{\lambda_2}=p^{-1}(\affgr^{\lambda_1})\times_{G(\co)}(\affgr^{\lambda_2}).
$
This realizes $\Gr\tilde\times\Gr$ as a stratified ind-scheme of ind-finite type
on which $G(\CO)$ acts by left multiplication.

Next, we have the convolution diagram
$$
\affgr\times\affgr \stackrel{p}{\gets}
G(\ck)\times\affgr \stackrel{q}{\to}
\affgr\tilde\times\affgr \stackrel{m}{\to}
\affgr.
$$
The map $p$
realizes $G(\ck)\times\affgr$ as a $G(\co)$-torsor
over $\affgr\times\affgr$.
The map $q$
realizes $G(\ck)\times\affgr$ as a $G(\co)$-torsor
over $\affgr\tilde \times\affgr$.
(Note this confirms that the
stratification of $\affgr\tilde\times\affgr$ is locally-trivial,
since the smooth map $q$ preserves strata, and the strata of 
$G(\ck)\times\affgr$ are locally-trivial group orbits.)
The map $m$ is the multiplication,
and is the direct limit of the maps
$
m:G(\CK)_{k_1}\times_{G(\CO)}\Gr_{k_2}\to\Gr_{k_1+k_2}.
$

Now to define the convolution product
$$
\odot:\epcat\times\epcat\to\epcat,
$$
recall that
for perverse sheaves $\psh_1,\psh_2$ in the category $\epcat$,
there is a unique perverse sheaf 
$\psh_1\tilde\boxtimes\psh_2$ in the category $\catp_{G(\co)}(\affgr\tilde \times\affgr)$
such that
$$
q^*(\psh_1\tilde\boxtimes\psh_2)\simeq
p^*(\psh_1\boxtimes\psh_2).
$$
The convolution is defined to be
$$
\psh_1\odot\psh_2=m_!(\psh_1\tilde\boxtimes\psh_2).
$$

Lusztig~\cite[Section 11]{Lu83} first proved that the convolution
takes perverse sheaves to perverse sheaves.
The assertion follows from his calculations
in the affine Hecke algebra.
The following result of Mirkovi{\'c}-Vilonen 
gives a direct geometric proof. See also \cite[Section 9]{NP01}.

\begin{Theorem} 
[\cite{MV04}, Lemma 4.3]\label{tsemismall}
The multiplication map $m:G(\CK)\times_{G(\CO)}\Gr {\to}\Gr$
is a stratified semismall map.
\end{Theorem}

The theorem refers to the stratification of
$\Gr$ by the $G(\CO)$-orbits $\Gr^\lambda$,
and the stratification of $\Gr\tilde\times\Gr$
by the twisted products $\gstratalta\tilde \times \gstrataltb$.
To be precise, it asserts that each of the maps
$
m:G(\ck)_{k_1}\times_{G(\co)}\affgr_{k_2}\to\affgr_{k_1+k_2}
$
is stratified semismall with respect to these stratifications.
See~\cite[Section 4]{MV04} for the notion of a stratified semsmall map.

\end{subsection}


\begin{subsection}{Tannakian formalism}

It is possible to place associativity and commutativity constraints
on the category $\epcat$ with respect to convolution,
and then to check that it is rigid. The hypercohomology
functor is an exact faithful tensor functor.
The most delicate part is the commutativity constraint, for which
see
\cite[Section 5]{MV04} for the most straightforward approach,
or~\cite[Remark in Section 1.1 and Theorem 1(b)]{Gait01} for an equivalent approach.
In the appendix, we sketch how one identifies the Tannakian
group of the category $\epcat$ with the dual group $\ch G$.

\begin{Theorem}[\cite{MV04}, Theorem 7.3, Proposition 6.3]\label{tfolk}
The category $\epcat$ is naturally a tensor category
with respect to convolution.
It is equivalent as a tensor category
to the category $\rep(\chG)$ of finite-dimensional representations
of the dual group $\chG$. 
Under the equivalence, the hypercohomology
of a perverse sheaf corresponds to the underlying vector space
of a representation.
\end{Theorem}


\begin{Example}\label{epsontorus}
As a topological space, the loop Grassmannian $\Gr_T$ 
of a torus $T$ is homeomorphic to 
the coweight lattice $\Lambda_T$.
The category $\catp_{T(\CO)}(\Gr_T)$ 
is equivalent as a tensor category to the category
$\Rep(\chT)$ 
of finite-dimensional representations of the dual torus $\chT$.
\end{Example}

\end{subsection}

\end{subsection}


\begin{subsection}{Real forms}\label{ssrf}
As mentioned in the preliminaries,
we only consider conjugations and real forms of ind-schemes
which are the limit of compatible families of conjugations
and real forms. In addition, the conjugations always respect
the stratifications considered, and
the strata of a real form 
are always taken to be the real forms of 
the strata of its complexification.
Thus once we have described
the ind-structure or stratification
of a complex ind-scheme, a real form of it inherits
an ind-structure and stratification by restriction.
In addition, we only consider the actions of real group-schemes
which are the limit of compatible families of actions of real forms.
Thus once we have described the action of a complex group-scheme
on a complex ind-scheme, a real form of the group-scheme inherits an action on 
the real form of the ind-scheme.

Let $\cor=\R[[t]]$ be the ring of real formal power series, and
let $\ckr=\R((t))$ be the field of real formal Laurent series.
Let $c$ be the conjugation of $\ck$ with respect to $\ckr$.
Let $\theta$ be the conjugation of the connected reductive complex algebraic group $G$
with respect to the real form $G_\R$.

The conjugation $\theta$ extends from $G$ to a conjugation of
$G(\CK)$ which we also denote by $\theta$.
For $g\in G(\ck)$ thought of as a map
$$
g:\Spec(\ck) \to G,
$$
the extended conjugation $\theta$ takes $g$ to the composite map
$$
\theta(g):\Spec(\ck) \stackrel{c}{\to} \Spec(\ck) \stackrel{g}{\to}G
   \stackrel{\theta}{\to}G.
$$
We may identify the resulting real form of $G(\ck)$
with the real stratified ind-scheme $G_\R(\ckr)$.

The conjugation of $G(\CK)$ restricts to a conjugation of $G(\CO)$,
and we may identify the resulting real form of $G(\co)$
with the real stable group-scheme $G_\R(\cor)$.

Since the conjugation of $G(\CK)$ preserves $G(\CO)$, it induces a conjugation of $\Gr$,
and we may identify the resulting real form $\raffgr$
with the real stratified ind-scheme of ind-finite type
$
G_\R(\ckr)/G_\R(\cor).
$
We call it the real loop Grassmannian of $G_\R$.
As a topological space, $\raffgr$
is known~\cite[Theorems 5.1 and 5.2]{Mitch88} to be homotopy equivalent to the based loop space
of the symmetric variety $G/K$ where $K$ is the complexification
of a maximal compact subgroup of $G_\R$.

\end{subsection}


\begin{subsection}{Real spherical sheaves}

The restriction of the stratification of $\affgr$ to
the real form $\raffgr$ coincides with the orbits
of the action of $G_\R(\cor)$ by left-multiplication.
The Cartan decomposition (which follows from the Bruhat decomposition of~\cite[Theorem 5.3]{Mitch88}) 
states
that the $G_\BR(\CO_\BR)$-orbit $\affgr^\lambda_\R\subset\raffgr$ is non-empty
if and only if $\lambda$ is a real coweight $\lambda\in\rcowt$,
and two orbits $\raffgr^{\lambda_1}$ and $\raffgr^{\lambda_2}$ coincide if and only if  
$\lambda_1$ and $\lambda_2$ 
are in the same $\CW_{G_\R}$-orbit in $\rcowt$.
Note that the $G_\R(\cor)$-orbit 
$\raffgr^\lambda$ contains the $G_\R$-orbit
$G_\R\cdot\lambda$ which is isomorphic to the real flag manifold $G_\R/P_\R^\lambda$
where the parabolic subgroup $P_\R^\lambda\subset G$ is the stabilizer of $\lambda$.

\begin{Proposition}\label{prstrat}
The real loop Grassmannian $\raffgr$ is the disjoint union of
the $G_\R(\cor)$-orbits $\raffgr^\lambda$ 
through the dominant real coweights $\lambda\in\rdcowt$.
The closure of the orbit $\raffgr^\lambda$ is the union of the orbits 
through $\mu\in\rdcowt$ with $\mu\leq\lambda$.
Each orbit $\raffgr^\lambda$ is a real vector bundle
over the real flag manifold $G_\R/P_\R^\lambda$, and its closure
is a real projective variety of 
dimension $2\langle \ch\rho,\lambda \rangle$.
\end{Proposition}

\begin{proof}
Thanks to Proposition~\ref{plgstrat} and the discussion preceding it,
it only remains to prove that if
$\lambda-\mu$
is a non-negative integral
linear combination of positive coroots of $G$, then
$\raffgr^\mu\subset\overline\raffgr^\lambda$.
By Lemma~\ref{lrcowtord}, it suffices to prove this when $\lambda-\mu$
is a positive coroot $\alpha \in R^\pos_{}$, or when $\lambda-\mu$ is of the form
$\theta(\alpha)+\alpha$, for a positive coroot $\alpha\in R^\pos_{}$,
but $\lambda-\mu$ is not a multiple of a positive coroot.
In the first case, we may find $\on{SL}_2(\R)\subset G_\R$ such that $\alpha$
is its positive coroot. Then the orbit through $\lambda$ of the one parameter subgroup
$U_{\alpha}(rt)\subset G_\R(\cor)$, for $r\in\R$, 
is isomorphic to $\R$, and its closure is isomorphic to $\R\BP^1$
with $\nu$ the point at infinity. 
In the second case, we may find $\on{SL}_2(\C)\subset G_\R$ such that $\theta(\alpha)+\alpha$
is its positive coroot. Then the orbit through $\lambda$ of the one parameter subgroup
$U_{\theta(\alpha)+\alpha}(ct)\subset G_\R(\cor)$, 
for $c\in\C$, is isomorphic to $\C$, and its closure is isomorphic to $\C\BP^1$
with $\nu$ the point at infinity. 
\end{proof}

The conjugation $c$ of $\CO$ induces a conjugation of $\Aut(\CO)$,
and we may identify the resulting real form 
with the real stable group-scheme $\Aut(\cor)$ of automorphisms of $\CO_\R$.
We define $\Aut^0(\cor)$ to be the connected component of the identity
consisting of orientation-preserving automorphisms.

We have the derived category 
$\catd_{G_\R(\cor)\rtimes\Aut^0(\cor)}(\raffgr)$ 
of $G_\R(\cor)\rtimes\Aut^0(\cor)$-equivariant sheaves on $\raffgr$, and
the derived category 
$\catd_{G_\R(\cor)}(\raffgr)$ of $G_\R(\cor)$-equivariant  
sheaves on $\raffgr$.
We also have the derived category $\catd_{G_\R,\strat}(\Gr)$ of $G_\R$-equivariant
sheaves on $\raffgr$
constructible with respect to the stratification by $G_\R(\cor)$-orbits.

The following lemma is useful in constructing
equivariant sheaves via functors which appear to only
provide sheaves constructible with respect to the 
orbit stratification.

\begin{Lemma}\label{lffeq}
The forgetful functors are equivalences
$$
\catd_{G_\R(\cor)\rtimes\Aut^0(\cor)}(\raffgr)\risom
\catd_{G_\R(\cor)}(\raffgr) \stackrel{\sim}{\to}  
\catd_{G_\R,{\mathcal S}}(\raffgr).
$$
\end{Lemma}

\begin{proof}
For the first equivalence, the group scheme $\Aut^0(\cor)$ is contractible,
and each $\Aut^0(\CO_\R)$-orbit in $\raffgr$ is contained in a $G_\R(\CO_\R)$-orbit. 
For the second,
the kernel of the projection $G_\R(\cor)\to G_\R$
is the inverse limit of unipotent groups.
\end{proof}

\end{subsection}


\begin{subsection}{Real convolution}

To define the convolution product
$$
\odot:\redcat\times\redcat\to\redcat,
$$
consider the real form of the convolution diagram
$$
\raffgr\times\raffgr \stackrel{p}{\gets}
G_\R(\ckr)\times\raffgr \stackrel{q}{\to}
\raffgr\tilde\times\raffgr \stackrel{m}{\to}
\raffgr.
$$
The product $G_\R(\ckr)\times\raffgr$ and twisted product 
$\raffgr\tilde\times\raffgr$ are real ind-schemes, stratified
by the real forms of the strata of their complexifications.

For sheaves $\sh_1,\sh_2$ in the category $\redcat$,
the convolution is defined to be
$$
\sh_1\odot\sh_2=m_!(\sh_1\tilde\boxtimes\sh_2)
$$
where $\sh_1\tilde\boxtimes\sh_2$
is the unique sheaf in the category $\catd_{G_\R(\cor)}(\raffgr\tilde \times\raffgr)$
such that 
$$
q^*(\sh_1\tilde\boxtimes\sh_2)\simeq
p^*(\sh_1\boxtimes\sh_2).
$$

\end{subsection}

\begin{subsection}{Real spherical perverse sheaves}

In this paper, 
we restrict our attention to sheaves supported
on certain connected components of $\raffgr$. 
The following
parameterization of $\pi_0(\raffgr)$ 
and the subsequent dimension assertion follow directly from
Proposition~\ref{prstrat}.

\begin{Lemma}\label{lparity}
The inclusion $\rcowt\subset\raffgr$ induces an isomorphism 
$
\rcowt/ (\rcowt\cap Q)\stackrel{\sim}{\to}\pi_0(\raffgr)
$
where $Q$ is the coroot lattice of $G$.
The strata
$\raffgr^\lambda$ in a connected component of $\raffgr$
are all even-dimensional or all odd-dimensional
depending on whether $2\langle\ch\rho,\lambda\rangle$ is congruent to $0$ or 
$1 \mod 2$.
\end{Lemma}

We define $\raffgr^+$ to be the union of the components of $\raffgr$
containing the even-dimensional strata.
Note by Proposition~\ref{prstrat} the action of
$G_\R(\cor)$ perserves each component of $\raffgr$.
Therefore we have the derived category
$\catd_{G_\R(\cor)}(\raffgr^+)$ of $G_\R(\cor)$-equivariant sheaves
on $\raffgr^+$.
Since the strata of $\raffgr^+$ are all even-dimensional,
there is a self-dual perverse $t$-structure on the 
category
$\catd_{G_\R(\cor)}(\raffgr^+)$. 
Its heart is the category
$\repcat$
of $G_\R(\cor)$-equivariant perverse sheaves on $\raffgr^+$.
All of the topological aspects of the theory of
perverse sheaves~\cite{BBD82,BL94} hold in this setting. For example,
the category $\repcat$ is abelian, and its simple objects
are the intersection cohomology sheaves of strata with coefficents
in irreducible $G_\R(\cor)$-equivariant local systems.

Note that the real convolution product takes sheaves supported on $\raffgr^+$
to sheaves supported on $\raffgr^+$ since the convolution of two even-dimensional
strata is also even-dimensional.
The real multiplication map $m:\raffgr\tilde\times\raffgr{\to}
\raffgr$ is a real form of the 
complex multiplication map, and the strata of its domain and codomain
are real forms of the strata of their complexifications, thus it
is a stratified semismall map by Theorem~\ref{tsemismall}. 
We conclude that the restriction of the real convolution product
to the category $\catd_{G_\R(\cor)}(\raffgr^+)$ is $t$-exact,
and preserves the category $\repcat$.

%

\end{subsection}


\begin{subsection}{Component refinement}\label{scompref}

In this paper, we in fact restrict our
attention to sheaves 
supported on only certain connected components of $\raffgr^+$.

Recall that $\affgr$ is homotopy equivalent to the based loop space of $G$,
and $\raffgr$ is homotopy equivalent to the based loop space of $G/K$
where $K$ is the complexification
of a maximal compact subgroup of $G_\R$. 
The fibration $K\to G\to G/K$ gives 
rise to an exact sequence of component groups
$$
\pi_0(\affgr)\stackrel{}{\to} \pi_0(\raffgr)\stackrel{}{\to} \pi_0(G_\R),
$$
which may be identified with the exact sequence
$$
\cowt/ Q\stackrel{\sigma}{\to} \rcowt/ (\rcowt\cap Q) 
\stackrel{\partial}{\to} \pi_0(G_\R)
$$
where $Q$ is the coroot lattice of $G$, $\sigma(\lambda)=\theta(\lambda)+\lambda$ is
the projection, 
and $\partial(\lambda)=[\lambda(-1)]$ is the boundary map.

We define $\raffgr^0$ be the union of the components
of $\raffgr$ in the image of $\sigma$, or equivalently in
the kernel of $\partial$. 

\begin{Lemma}\label{lcompinc}
The components $\raffgr^0$ are a subset of the components $\raffgr^+$.
\end{Lemma}

\begin{proof}
Recall that $M\subset G$ is the Levi factor of the complexification $P\subset G$
of the minimal parabolic subgroup $P_\R\subset G_\R$.
Using Lemma~\ref{lparity}, we calculate the parity of the dimension of strata
$$
\langle2\ch\rho,\sigma(\lambda)\rangle=\langle2\ch\rho-2\ch\rho_M,\theta(\lambda)+\lambda\rangle
=2\langle2\ch\rho-2\ch\rho_M,\lambda\rangle=0\mod 2.
$$
In the above equation,
the identity $\langle2\ch\rho_M,\sigma(\lambda)\rangle=0$ gives the first equality,
and the identity $\ch\theta(2\ch\rho-2\ch\rho_M)=2\ch\rho-2\ch\rho_M$ gives the second.
\end{proof}

\end{subsection}

\end{section}


\begin{section}{Real Beilinson-Drinfeld Grassmannians}\label{srbdg}

In this section, we recall some basic properties of 
Beilinson-Drinfeld Grassmannians, then
describe similar properties of their real forms.
We also collect some basic results about embeddings and stratifications.
Throughout this section and later sections, 
we fix a smooth projective complex algebraic curve $X$ with non-empty real form $X_\R$, 
and let $c$ denote the conjugation
of $X$ with respect to $X_\R$.

In what follows, we shall define certain ind-schemes by the functors 
from $\C$-algebras to sets which they represent.
The sets shall always be taken up to the natural notion of equivalence.
 
\begin{subsection}{Local loop Grassmannians}

Fix $x\in X$.
Let $\co_x$ be the completion of the local ring of $X$ at $x$,
and let $\ck_x$ be its fraction field.

The local loop Grassmannian $\xaffgr=G(\ck_x)/G(\co_x)$
is a stratified ind-scheme of ind-finite type isomorphic to the loop Grassmannian $\affgr$.
To see this, choose a formal coordinate at $x$.
The identification $\co_x\simeq\co$
induces a bijection
$
\xaffgr\stackrel{\sim}{\to}\affgr.
$
If we choose a different formal coordinate,
the resulting bijection differs by 
the action on $\affgr$
of an element of $\Aut(\co)$.
Since each of the varieties in the family $\affgr_k$ is $\Aut(\co)$-invariant,
we see that $\xaffgr$ is naturally an ind-scheme of ind-finite type,
and the bijection is an isomorphism for any choice of local coordinate.
Furthermore, 
the isomorphism
takes each $G(\co_x)$-orbit in $\xaffgr$ 
to a $G(\co)$-orbit in $\affgr$.
Since each $G(\co)$-orbit in $\affgr$ is $\aut(\co)$-invariant,
the orbit correspondence is  
independent of the choice of local coordinate.
We may therefore unambiguously index the $G(\co_x)$-orbits $\xaffgr^\lambda
\subset\xaffgr$ by dominant coweights $\lambda\in\dcowt$. 

Similarly, for $x\in X_\R$,
we have the real local loop Grassmannian $\xraffgr$,
an isomorphism $\xraffgr\risom\raffgr$ for every choice of real formal coordinate
at $x$,
and the strata $\xraffgr^\lambda\subset \xraffgr$, for
dominant real coweights $\lambda\in\rdcowt$.

The following proposition represents an important change
in perspective. It was first proved for $G=\on{SL}_n$ in \cite[Proposition 2.1]{BeauLa94}.
For a $\C$-algebra $\alg$ and scheme $Z$, let $Z_\alg$ denote 
the product $Z\times\spec(\alg)$.

\begin{Proposition}[\cite{LaSo97}, Proposition 3.10]\label{ploclg}
The local loop Grassmannian $\xaffgr$ represents the functor
$\alg \mapsto \{(\tors, \nu)\}$, 
where
$\tors$ is a $G$-torsor on  $X_\alg$,
and $\nu$ is a trivialization of $\tors$ over $(X\setminus x)_\alg.$
\end{Proposition}



\end{subsection}


\begin{subsection}{Beilinson-Drinfeld Grassmannians}
Some of the material discussed here may be found in~\cite[Sections 5.3.10, 5.3.11]{BD}
and \cite[Section 5]{MV04}.
The Beilinson-Drinfeld Grassmannian $\naffgr$ 
of the group $G$ 
is the ind-scheme of ind-finite type which represents
the functor 
$\alg\mapsto
\{\left(
(x_1,\ldots,x_n), 
\tors, 
\nu 
\right)\}$,
where
$(x_1,\ldots,x_n)\in X^n(\alg)$,
$\tors$ is a $G$-torsor on $X_\alg$,
and $\nu$ is a trivialization of $\tors$ over 
$X_\alg\setminus (x_1\cup\cdots\cup x_n)$.
Here we think of the points $x_i:\spec(\alg)\to X$
as subschemes of $X_\alg$ by taking their graphs.
See Section~\ref{semb} for confirmation that $\naffgr$ is indeed
an ind-scheme of ind-finite type.

One of the most important properties of $\naffgr$ is its
factorization with respect to the projection
$
\pi:\naffgr\to X^n.
$ 
To describe this, we introduce the incidence stratification of $X^n$. 
The strata are indexed by partitions of the set $\{1,\ldots,n\}$.
The stratum $T_\tau$ indexed by the partition $\tau$ 
consists of the points $(x_1,\ldots,x_n)\in X^n$
such that the coincidences specified by $\tau$ 
occur among the points
$x_1,\ldots,x_n\in X$.

\begin{Proposition}
\label{pfiber}
For each stratum $T_\tau\subset X^n$,
there is a canonical isomorphism
$$
\naffgr|T_\tau\stackrel{\sim}{\to} (\prod_{i=1}^k \oneaffgr)|T_0,
$$
where $k$ is the number of parts in the partition $\tau$,
and $T_0$ denotes the open stratum of distinct points $y_1,\ldots,y_k\in X$.
\end{Proposition}

\begin{proof}
The map is defined by
$
((x_1,\ldots,x_n),\tors,\nu)
\mapsto
\prod_{i=1}^k  (y_i,\tors_i,\nu_i),
$ 
where $y_1,\ldots,y_k\in X$ are the distinct points such that
there is an equality of sets
$
\{y_1, \ldots,y_k\}=\{x_1,\ldots,x_n\},
$
the torsor $\tors_i$ coincides with $\tors$ over
$X\setminus (y_1\cup\cdots\cup \hat y_i\cup\cdots \cup y_k),$
and the trivialization $\nu_i$
coincides with $\nu$ over 
$X\setminus (y_1\cup\cdots\cup y_k).$
We leave it to the reader to check that the map is an isomorphism.
\end{proof}

We next describe the standard stratification of $\naffgr$.
The strata of $\oneaffgr$
are indexed by the dominant coweights $\lambda\in\dcowt$.
The stratum of $\oneaffgr$ indexed by $\lambda$ 
consists of the union over all fibers $\oneaffgr|x$
of the points which 
map to the $G(\co_x)$-orbit $\xaffgr^\lambda$ in $\xaffgr$
via the canonical isomorphism
$\oneaffgr|x\stackrel{\sim}{\to} \xaffgr.$
In Section~\ref{sconstpsh}, we shall see that it is useful to realize
$\oneaffgr$ as the twisted product
of $\affgr$ and the curve $X$. From this perspective,
each stratum of $\oneaffgr$ is the twisted product
of a stratum of $\affgr$ with $X$.

In general, the strata of $\naffgr$ are indexed by labeled partitions
$(\tau,\{\lambda_1,\ldots,\lambda_k\})$,
where $\tau$ is a partition of $\{1,\ldots,n\}$ into $k$ parts, 
and $\lambda_1,\ldots,\lambda_k\in\dcowt$ are dominant coweights each assigned to a part
of the partition.
The stratum of $\naffgr$ indexed by $(\tau,\{\lambda_1,\ldots,\lambda_k\})$ consists of the points
which project to the stratum $T_\tau\subset X^n$,
and which map via the isomorphism
$
\naffgr|T_\tau\stackrel{\sim}{\to} (\prod_{i=1}^k \oneaffgr)|T_0
$
to the product stratum indexed by $\{\lambda_1,\ldots,\lambda_k\}$.

\end{subsection}


\begin{subsection}{Real forms}\label{ssrfbd}
Recall that $c$ denotes the conjugation of $X$ with respect to $X_\R$,
and $\theta$ denotes the conjugation of  $G$ 
with respect to $G_\R$.
The pair of conjugations induce a conjugation $\theta$ of $\naffgr$ as follows.
For data $((x_1,\ldots,x_n),\tors,\nu)$,
choose a cover $\{U_i\subset X\}_i$ so that we may realize $\tors$ 
as the $G$-torsor obtained from the disjoint union
$\sqcup_i (U_i\times G)$ via gluing maps
$\varphi_{ij}:U_i\cap U_j\to G$.
Then we define the conjugation to be
$$
\theta((x_1,\ldots,x_n),\tors,\nu)=((c(x_1),\ldots,c (x_n)),\tors_c^\theta,\nu_c^\theta).
$$
Here $\tors_c^\theta$ is the $G$-torsor obtained from 
the disjoint union
$\sqcup_i (c(U_i)\times G)$ via gluing maps
$$
c(U_i\cap U_j)\stackrel{c}{\to} U_i\cap U_j\stackrel{\varphi_{ij}}{\to} G\stackrel{\theta}{\to} G.
$$
Thinking of the trivialization $\nu$ as a section of $\tors$
over $X\setminus (x_1\cup\cdots\cup x_n)$, 
the trivialization $\nu_c^\theta$ is the composite section
$$
\nu_c^\theta:X\setminus c(x_1\cup\cdots\cup x_n)\stackrel{c}{\to}
X\setminus (x_1\cup\cdots\cup x_n)\stackrel{\nu}{\to} \tors
\stackrel{f^\theta_c}{\to} \tors_c^\theta,
$$
where $f^\theta_c:\tors\to\tors_c^\theta$ is the map defined by
$f^\theta_c((u_i,g))=(c(u_i), \theta(g)),$ for $u_i\in U_i$.

Now fix a permutation $\omega$ 
in the symmetric group $\Sigma_n$ such that $\omega^2=e$.
Then we obtain involutions of $X^n$ and of $\naffgr$ by permuting
the labels of the points 
$$
(x_1,\ldots,x_n)\mapsto (x_{\omega(1)},\ldots,x_{\omega(n)}).
$$ 
The composition of a conjugation and an involution is another
conjugation.
Let $c_{\omega}$ be the conjugation of $X^n$ defined by 
$$
c_{\omega}(x_1,\ldots,x_n)=(c(x_{\omega(1)}),\ldots,c(x_{\omega(n)})),
$$
and let $\theta_{\omega}$ be the conjugation of  
$\naffgr$ defined by
$$
\theta_\omega((x_1,\ldots,x_n),\tors,\nu) 
= (c(x_{\omega(1)}),\ldots,c(x_{\omega(n)}),\tors^\theta_c,\nu^\theta_c).
$$

Let $X^{(\omega)}_\R$ be the real form of $X^n$ with respect to the conjugation $c_\omega$,
and let $\wraffgr$ be the real form of $\naffgr$ with respect to 
the conjugation $\theta_\omega$.
When $\omega$ is the identity permutation, 
we write $X^n_\R$ in place of $X^{(\omega)}_\R$,
and $\naffgr_\R$ in place of $\wraffgr$.


As usual we stratify the real forms 
by taking the real forms
of the strata of their complexifications.
The strata of $\base$ are indexed
by those partitions $\tau$ of $\{1,\ldots,n\}$ such that the action of $\omega$
sends each part of $\tau$ to another part of $\tau$.
The strata of $\wraffgr$ are indexed by labelled partitions
$(\tau,\{\lambda_1,\ldots,\lambda_k\},\{\mu_1,\ldots,\mu_\ell\})$,
where $\tau$ is a partition of $\{1,\ldots,n\}$ such that the action of $\omega$
sends each part of $\tau$ to another part of $\tau$ with 
$k$ parts fixed and
$\ell$ pairs of parts exchanged,
$\lambda_1,\ldots,\lambda_k\in\rdcowt$ are dominant real coweights each assigned to a part
of $\tau$ fixed by $\omega$, and
$\mu_1,\ldots,\mu_\ell\in\dcowt$ are dominant coweights
each assigned to a pair of parts exchanged by $\omega$.

To describe the factorization of $\wraffgr$ with respect to the projection
$\pi:\wraffgr\to X^{(\omega)}$,
it is convenient to break the symmetry as assumed in the following.

\begin{Proposition}\label{prfiber}
Let $\tau$ be a partition such that $\omega$ sends each of its parts to another.
For each pair of parts of $\tau$ exchanged by $\omega$,
choose one of the parts.
Then there is a canonical strata-preserving isomorphism
$$
\wraffgr|T_\tau \simeq (\prod_{i=1}^k \oneraffgr\times\prod_{j=1}^\ell \oneaffgr)|T_0,
$$
where $k$ is the number of parts in $\tau$ fixed by $\omega$,
$\ell$ is the number of pairs of parts exchanged by $\omega$, 
and $T_0$ denotes the open stratum of distinct points   
$y_1,\ldots,y_k\in X_\R$, $z_1,\ldots,z_\ell\in X$.
\end{Proposition}

\begin{proof}
By Proposition~\ref{pfiber} and the definitions,
$\wraffgr|T_\tau$ is the real form of the product
$$
(\prod_{i=1}^k \oneaffgr\times\prod_{j=1}^\ell (\oneaffgr\times\oneaffgr))|T_0
$$
with respect to the standard conjugation on the factors of type $\oneaffgr$,
and the composition of the standard conjugation with the exchange involution
on the factors of type $\oneaffgr\times\oneaffgr$.
The choice of a part in each pair of parts exchanged by $\omega$ 
distinguishes one of the factors $\oneaffgr$ 
in each of the products $\oneaffgr\times\oneaffgr$.
Projection of the real form of each of the products $\oneaffgr\times\oneaffgr$
to its distinguished factor $\oneaffgr$ is an isomorphism.
\end{proof}

If we fix an identification
$$
X^{(\omega)}_\R\simeq X_\R^r\times X^s
$$ 
where $r$ is the number of fixed points of $\omega$,
$s$ is the number of pairs of points exchanged by $\omega$,
and $n=r+2s$,
then
by the proposition, for 
$(x_1,\ldots,x_r,z_1,\ldots,z_s)\in X^r_\R\times {X}^{s},$
we have a canonical isomorphism
$$
\wraffgr|(x_1,\ldots,x_r,z_1,\ldots,z_s)\stackrel{\sim}{\to}  
\prod_{i=1}^k {{}_{y_i}\hspace{-0.3em}\affgr_{\R}}\times
\prod_{j=1}^\ell {}_{w_j}\hspace{-0.3em}\affgr.
$$
where $y_1,\ldots,y_k\in X_\R, w_1,\ldots,w_\ell\in X\setminus X_\R$ are the distinct points with an
equality of sets
$$\{y_1,\ldots,y_k,w_1,\ldots, w_\ell\}=\{x_1,\ldots,x_r,z_1,\ldots,z_s\}.$$

\end{subsection}


\begin{subsection}{Embeddings}\label{semb}

It will be useful to have an embedding of the Beilinson-Drinfeld Grassmannian
$\naffgr$ in an ind-scheme that is the limit of smooth varieties.
We describe one such construction here. 

First, fix an embedding $G\subset\GL_N$ to obtain an embedding 
$\naffgr\subset\naffgr_N$ of the corresponding Beilinson-Drinfeld Grassmannians.

Next, 
let $\naffgr_{N,k}$ be the variety that represents the functor
$\alg\mapsto
\{(x_1,\ldots,x_n, \M)\},
$
where
$x_1,\ldots,x_n\in X(\alg),$ and 
$$
\M\subset \CO^{\oplus N}_{X_\alg}(k(x_1\cup\cdots\cup x_n))
$$ 
is a subsheaf of $\CO_{X_\alg}$-modules such that
$$
\co^{\oplus N}_{X_\alg}(-k(x_1\cup\cdots\cup x_n))\subset\M,
$$
and $\co^{\oplus N}_{X_\alg}(k(x_1\cup\cdots\cup x_n))/ \M$ is $\spec(\alg)$-flat.

The inclusions 
$\co^{\oplus N}_{X_\alg}(k(x_1\cup\cdots\cup x_n))\subset
\co^{\oplus N}_{X_\alg}((k+1)(x_1\cup\cdots\cup x_n))$
induce
canonical closed embeddings
$
\naffgr_{N,k}\to\naffgr_{N,k+1}.
$


\begin{Proposition}
The direct limit of the varieties $\naffgr_{N,k}$ is canonically isomorphic
to the Beilinson-Drinfeld Grassmannian $\naffgr_N$ of $\GL_N$.
\end{Proposition}

\begin{proof}
A $\GL_N$-torsor over $X$ with a trivialization over $X\setminus(x_1\cup\cdots\cup x_n)$
defines a vector bundle $\CV$ over $X$ and an isomorphism 
$\CO(\CV)\stackrel{\sim}{\to}\CO^{\oplus N}_X$ over $X\setminus(x_1\cup\cdots\cup x_n)$.
For $k$ sufficiently large, the isomorphism extends to a map 
$\CO(\CV)\to\co^{\oplus N}_{X}(k(x_1\cup\cdots\cup x_n))$ over $X$. 
And again for $k$ sufficiently large, the image of the extended map will 
contain $\co^{\oplus N}_{X}(-k(x_1\cup\cdots\cup x_n))$. 
We take the subsheaf $\CM\subset \co^{\oplus N}_{X}(k(x_1\cup\cdots\cup x_n))$
to be the image for $k$ large. We leave it to the reader to
check that this gives an isomorphism.
\end{proof}

Finally, 
let $\nadgr_{N,k}$ be the variety that represents the functor
$\alg\mapsto
\{(x_1,\ldots,x_n, \M)\},
$
where
$x_1,\ldots,x_n\in X(\alg),$ and 
$$
\M\subset \CO^{\oplus N}_{X_\alg}(k(x_1\cup\cdots\cup x_n))
$$ 
is a subsheaf of $\C$-modules such that
$$
\co^{\oplus N}_{X_\alg}(-k(x_1\cup\cdots\cup x_n))\subset\M,
$$
and $\co^{\oplus N}_{X_\alg}(k(x_1\cup\cdots\cup x_n))/ \M$ is $\spec(\alg)$-flat.

\begin{Lemma}
The variety $\nadgr_{N,k}$ is smooth.
\end{Lemma}

\begin{proof}
We may identify $\nadgr_{N,k}$ with 
the Grassmann bundle (of $\C$-subspaces of any dimension) of the vector bundle
$H^0(X^n,\M^N_k)\to X^n$,
for the quotient sheaf
$$
\CM^N_k=\co^{\oplus N}_{X}(k(x_1\cup\cdots\cup x_n))/
\co^{\oplus N}_{X}(-k(x_1\cup\cdots\cup x_n)).
$$
\end{proof}

Define the ind-scheme $\nadgr_N$ 
to be the direct limit of the varities $\nadgr_{N,k}.$
We conclude that for the choice of an embedding $G\subset\GL_N$,
we obtain embeddings of ind-schemes 
$$
\naffgr\subset\naffgr_N\subset\nadgr_N
$$
with $\nadgr_N$ the limit of a family of smooth varieties.

\begin{Remark}
The above discussion confirms
that the Beilinson-Drinfeld Grassmannian $\naffgr$
is an ind-scheme of ind-finite type.
\end{Remark}

Before continuing on,
we note that there is a factorization of $\adgr$
with respect to the projection $\pi:\adgr\to X^n$
which restricts to give the factorization of $\naffgr$ of 
Proposition~\ref{pfiber}.
For each stratum $T_\tau\subset X^n$, let 
$T^\tau \subset X^n$ be the
union of the strata of $X^n$
containing $T_\tau$ in their closures.

\begin{Proposition}\label{pexfiber}
For each stratum $T_\tau\subset X^n$,
there is a canonical isomorphism
$$
\nadgr_N|T^\tau\stackrel{\sim}{\to} 
(\prod_{i=1}^k \adgr^{(n_i)}_{N})|T^\tau 
$$
where $k$ is the number of parts of $\tau$,
and $n_i$ is the size of the part $\tau_i$.
\end{Proposition}

\begin{proof}
If the sheaf $\M\in\adgr^{(n)}_{N}$ 
is supported at points $x_1,\ldots,x_n\in X$ 
with $(x_1,\ldots,x_n)\in T^\tau$, 
then the coincidences among the points 
$x_1,\ldots,x_n$ are coarser than the coincidences specified by
$\tau$.
We may therefore define a map
$
\M\mapsto \prod_{i=1}^k \M_i,
$
requiring that 
$
\M\simeq \sum_{i=1}^k \M_i,
$
and that each sheaf $\M_i$ is supported at those points 
contained in the part $\tau_i$ of the partition $\tau$.
We leave it to the reader to check that this map is an isomorphism.
\end{proof}

\end{subsection}


\begin{subsection}{Stratifications}\label{ssstrat}
See \cite{Math70} or \cite[Section 1.2, and Section 2.A.1]{GM88} 
for the notions of a Whitney stratification of a manifold 
and of a Thom stratified map between manifolds.

We say that a stratification of a variety $V$ is 
a Whitney stratification if for some embedding $V\subset M$ into a smooth
manifold, the stratification of $M$ by the strata of $V$ 
and the complement $M\setminus V$ is a Whitney stratification.

We say that a stratification of an ind-variety $Z$ is a Whitney stratification
if the closure of each stratum of $Z$ is Whitney stratified by the strata of $Z$ in the
closure.

We say that a map $V\to N$, where $V$ is a stratified variety and $N$ is a smooth manifold,
is a Thom stratified map
if for some commutative diagram
$$
\begin{array}{ccc}
V & \subset & M \\
\downarrow & & \downarrow \\
N & = & N
\end{array}
$$
with $M$ smooth, the map $M\to N$ is a Thom stratified map
with respect to the stratification of $M$ by the strata of $V$ and the complement $M\setminus V$.

We say that a map $Z\to N$, where $Z$ is a stratified ind-variety and $N$ is a smooth
manifold, is a Thom stratified map
if the restriction of the map to the closure of each stratum of $Z$ is a Thom stratified
map with respect to the stratification of the closure by the strata of $Z$ in the closure.

\begin{Proposition}\label{pthom}
The stratifications of the Beilinson-Drinfeld Grassmannian $\naffgr$
and its real form $\wraffgr$
are Whitney stratifications. The projection
$
\pi:\naffgr\to X^n
$ 
and its restriction $\pi:\wraffgr\to\base$
are Thom stratified maps.
\end{Proposition}

\begin{proof} 
First, note that if we choose $G\subset GL_N$ to be equivariant
with respect to conjugation, then the embeddings
$$
\naffgr\subset\naffgr_N\subset\nadgr_N
$$
of the previous section
will also be equivariant with respect to conjugation.
Therefore, by
the following
easily-verified lemma, 
it suffices to prove the proposition in the complex case.

\begin{Lemma}\label{leqstrat}
Let $M$ and $N$ be Whitney stratified manifolds, and let $K$ be a compact
group acting smoothly on $M$ and $N$ 
such that the actions preserve the stratifications.
Then the fixed point manifolds $M^K$ and $N^K$ are 
Whitney stratified by the fixed  
points of the strata. If $M\to N$ is a Thom stratifed
$K$-equivariant map, then $M^K\to N^K$ is Thom stratified as well.
\end{Lemma}

Now, for the closure $V$ of a stratum of $\naffgr$,
we may find $k$ large so that we have an embedding
$$
 V\subset\naffgr_{N,k}\subset\nadgr_{N,k}.
$$

Recall that the factorization of $\nadgr_N$ restricts to give
that of $\naffgr$, and the strata of $\naffgr$ are by definition products
with respect to the factorization.
Therefore it suffices to 
verify Whitney's condition $B$ and Thom's condition $A_\pi$ are 
satisfied at a point $\M\in\adgr^{(n)}_{N,k}$ such
that
$\pi(\M)=(x,\ldots,x)$.
To accomplish this, consider the group $G^\loc_x$ of maps
$
X_x\to G,
$
where $X_x$ is the localization of $X$ at $x$.

\begin{Lemma}\label{llocsurj}
For $\ell\geq0$, the composite homomorphism
$
G^\loc_x\to G(\co_x)\to G(\jet_x^\ell) 
$
is surjective.~$\Box$
\end{Lemma}

\begin{proof}
In \cite[Lemma 4]{Gait01}, the analogous assertion is proven for Iwahori subgroups
which immediately implies this assertion.
\end{proof}

Using the lemma,
it is easy to confirm Thom's condition $A_\pi$ at the point $\CM$.

To verify Whitney's condition $B$ at the point $\CM$,
let $\M_i$ be a sequence coverging to $\M$ 
in the stratum $S$ containing $\M$, 
let $\shL_i$ be a sequence also
converging to $\M$ in some stratum $R$, and
suppose that the limits $\ell=\lim_i \overline{\M_i \shL_i}$ and
$\tau=\lim_i T_{\shL_i}R$ exist. 
Since $\pi(\M)=(x,\ldots,x)$, 
and $X$ is smooth, we may assume that the sequence $\M_i$
is in the same fiber as $\M$.
Then using the lemma, we may assume that the sequence 
$\M_i$ is in fact constant equal to $\M$.
Now the assertion that $\ell\subset\tau$ for such sequences follows 
from standard arguments such as the Curve Selection Lemma~\cite[Chapter 3]{Mi68}.
\end{proof}

\end{subsection}

\end{section}
%



\begin{section}{Specialization}\label{sspec}

In this section, we define a functor 
that takes perverse sheaves on
the loop Grassmannian $\affgr$ to perverse 
sheaves on the subspace $\raffgr^0$ of the real form $\raffgr$.

Note that the complement of the
real curve $X_\R$ in its complexification $X$ is the union of 
two connected components. Throughout what follows,
we distinguish one of the components and denote it by $X_+$.
This also distinguishes an orientation of $X_\R$
by the rule that
the complex structure of $X$ takes a positively-oriented tangent vector
to a tangent vector pointing into $X_+$.

\subsection{Global specialization}\label{ssglob}
To define the global specialization functor 
$$
\on{R}_X:\catd_{G,\strat}(\oneaffgr)\to\catd_{\GR,\strat}(\oneraffgr),
$$
we work with the real form
$\sraffgr$ of the Beilinson-Drinfeld Grassmannian
$\twoaffgr$, where $\sigma$ is the non-trivial
element of the symmetric group $\Sigma_2$.

We identify the projection 
$\sraffgr\to X^{(\sigma)}_\R$
with a map
$
\sraffgr\to X
$
via the isomorphism $X^{(\sigma)}_\R\stackrel{\sim}{\to} X$ defined by 
$(z,c(z))\mapsto z$.
By Proposition~\ref{prfiber},
we have canonical identifications
\begin{eqnarray*}
& \sraffgr|X_+\simeq \oneaffgr|X_+ & \\
& \sraffgr|X_\R \simeq \oneraffgr. &
\end{eqnarray*}

Consider the diagram
$$
\begin{array}{ccccccccc}
\oneaffgr|X_+  & \simeq & 
\sraffgr|X_+ & \stackrel{j}{\to} 
& \sraffgr|X & \stackrel{i}{\gets} & \sraffgr|X_\R 
& \simeq & \oneraffgr \\
 & \searrow & \downarrow & & \downarrow & & \downarrow & \swarrow &\\
  && X_+ & \to & X & \gets & X_\R && 
\end{array}
$$ 

Define the global specialization
by 
$$
\on{R}_X(\sh)= i^*j_*(\sh|X_+).
$$
Informally speaking, it is the nearby cycles in the family $\sraffgr\to X$.

By Proposition~\ref{pthom}, the stratification of $\sraffgr$
is a Whitney stratification and so locally-trivial.
Thus $\on{R}_X$ takes $\strat$-constructible sheaves to 
$\strat$-constructible
sheaves.
Since all of the maps in the 
above diagram are  $G_\R$-equivarant, $\on{R}_X$ takes $G$-equivariant
sheaves to $G_\R$-equivariant sheaves.


\begin{subsection}{Perverse sheaves constant along $X$}\label{sconstpsh}
Let $\encurve\to X$ be the $\Aut(\co)$-torsor
of smooth maps $\Spec(\co)\to X$.
It is the inverse limit of the $\Aut(\jet^\ell)$-torsors $\encurve^\ell\to X$
of smooth maps $\Spec(\CJ^\ell)\to X$.

Recall that for $x\in X$, we have an identification $\oneaffgr|x\simeq \xaffgr$, and 
for the choice of a formal coordinate at $x$, we  obtain an isomorphism
$\oneaffgr|x\simeq\affgr$. It follows that
$\oneaffgr$ is the twisted product obtained from the $\Aut(\CO)$-torsor
$\encurve\to X$ and the action of $\Aut(\CO)$ on $\affgr$.
In other words,  there is a canonical isomorphism
$$
\oneaffgr\simeq\encurve\times_{\Aut(\co)}\affgr.
$$

Consider the diagram
$$
X\times\affgr\stackrel{p}{\gets} 
\encurve \times\affgr\stackrel{q}{\to}
\encurve \times_{\Aut(\co)}\affgr\simeq\oneaffgr.
$$

By Lemma~\ref{autoe}, we may define a functor
$$
\rho:\catp_{G(\co)}(\affgr)\to\catd_{G, \mathcal S}(\oneaffgr)
$$
by the formula
$$
\rho(\psh)=\C_X\tilde\boxtimes\psh
$$
where 
$$
q^*(\C_X\tilde\boxtimes\psh)= p^*(\C_X\boxtimes\psh).
$$
(Although it may be more natural to shift here
so that $\rho$ is perverse,
we have found it convenient to avoid all such shifts.)
%

Define the category $\catp_{G, \hat{\mathcal S}}(\oneaffgr)$ to be the 
strict full subcategory of $\catd_{G, \mathcal S}(\oneaffgr)$ whose objects are
isomorphic to objects of the form $\rho(\psh)$, where $\psh$ runs through all objects
of $\epcat$. Note that objects in the
subcategory $\catp_{G, \hat{\mathcal S}}(\oneaffgr)$ 
are equivariant with respect to
the groupoid of pairs of points $(x,y)\in X^2$, and an isomorphism
between their formal neighborhoods.

For a choice of $x\in X$, and formal coordinate at $x$, the isomorphism
$
\oneaffgr|x\simeq\affgr
$
provides an inverse 
$$
\rho^{-1}:\catp_{G, \hat{\mathcal S}}(\oneaffgr)\to\catp_{G(\co)}(\affgr)
$$
defined by restriction to the fiber 
$$
\rho^{-1}(\psh)=\psh|x.
$$ 
By the definition of $\catp_{G, \hat{\mathcal S}}(\oneaffgr)$,
the inverse functors for different choices of $x\in X$, and formal coordinate at $x$,
are canonically isomorphic.




\end{subsection}


\begin{subsection}{Sheaves constant along $X_\R$}

Here we work with the $\Aut^0(\cor)$-torsor $\encurve_\R^0\to X_\R$ 
of orientation-preserving
smooth maps $\Spec(\cor)\to\encurve_\R$,
and the canonical isomorphism
$$
\encurve_\R^0 \times_{\Aut^0(\cor)}\raffgr\simeq\oneraffgr.
$$


Consider the diagram
$$
X_\R\times\raffgr\stackrel{p}{\gets} 
\encurve_\R^0 \times\raffgr\stackrel{q}{\to}
\encurve_\R^0 \times_{\Aut^0(\cor)}\raffgr\simeq\oneraffgr.
$$

By Lemma~\ref{lffeq}, we may define a functor
$$
\rho_\R:\redcat\to\catd_{G_\R, \mathcal S}(\oneraffgr)
$$
by the formula
$$
\rho_\R(\sh)=\C_{X_\R}\tilde\boxtimes\sh
$$ 
where  
$$
q^*(\C_{X_\R}\tilde\boxtimes\sh)= p^*(\C_{X_\R}\boxtimes\sh).
$$
Note that here we are able to define the functor 
on the entire derived
category, not only on the catgory of perverse sheaves. 

Define the category $\catd_{G_\R, \hat{\mathcal S}}(\oneraffgr)$ to be the 
strict full subcategory of $\catd_{G_\R, \mathcal S}(\oneraffgr)$ 
whose objects are isomorphic to objects of the form $\rho_\R(\sh)$,
where $\sh$ runs through all objects of $\redcat$.
An object in the category $\catd_{G_\R, \mathcal S}(\oneraffgr)$
is in the subcategory $\catd_{G_\R, \hat{\mathcal S}}(\oneraffgr)$
if and only if it is equivariant with respect to the groupoid of pairs
of points $(x,y)\in X_\R^2$, and an orientation-preserving 
isomorphism between their
formal neighborhoods.
(Since the group $\Aut^0(\cor)$ is contractible,
such equivariance is a property not a structure.)

For the choice of $x\in X_\R$, and formal coordinate at $x$, the
resulting isomorphism
$
\oneraffgr|x\simeq\raffgr
$
provides an inverse
$$
\rho_\R^{-1}:\catd_{G_\R, \hat{\mathcal S}}(\oneraffgr)\to\redcat
$$
defined by restriction to the fiber 
$$
\rho_\R^{-1}(\sh)=\sh|x.
$$ 
By the definition of $\catd_{G_\R, \hat{\mathcal S}}(\oneraffgr)$,
the inverse functors for different choices of $x\in X_\R$, and formal coorcinate at $x$,
are canonically isomorphic.

%



\end{subsection}

\begin{subsection}{Local specialization}\label{ssloc}

To define the local specialization
$$
\on R:\epcat\to\redcat,
$$
we use the following proposition.

\begin{Proposition}
The global specialization descends to a functor
$$
\on{R}_X:\catp_{G,\hat
{\mathcal S}}(\oneaffgr)\to\catd_{G_\R,\hat{\mathcal S}}(\oneraffgr).
$$

\end{Proposition}

\begin{proof}
A sheaf $\psh$ in the category 
$\catp_{G,\hat{\mathcal S}}(\oneaffgr)$ is equivariant with respect to the groupoid
of a pair of points $(x,y)\in X^2$, and an isomorphism between their formal neighborhoods.
For points $x,y\in X_\R$, 
we may find analytic disks $D(x),D(y)\subset X$, and a conjugation-equivariant
isomorphism between them preserving the orientations of their real forms 
$D_\R(x),D_\R(y)\subset X_\R$ such that it induces an isomorphism
$$
\affgr^{(\sigma)}_{\R}|D(x)\simeq\affgr^{(\sigma)}_{\R}|D(y).
$$
Therefore the global specialization $\on R_X(\psh)$ is equivariant with respect
to the groupoid of points $(x,y)\in X_\R^2$, analytic disks $D(x),D(y)\subset X$, 
and a conjugation-equivariant
isomorphism between them preserving the orientations of their real forms 
$D_\R(x),D_\R(y)\subset X_\R$.
It follows that the sheaf $\on R_X(\psh)$ is equivariant 
with respect to the groupoid
of a pair of points $(x,y)\in X_\R^2$, and an isomorphism between their formal neighborhoods,
and thus it is in the category $\catd_{G_\R,\hat{\mathcal S}}(\oneraffgr).$
\end{proof}

By the proposition, we may define the local specialization
by
$$
\on R(\psh)= \rho_\R^{-1} (\on{R}_X (\rho(\psh))).
$$
Here we use Lemma~\ref{lffeq} to
obtain a 
$G_\R(\cor)$-equivariant sheaf from a $G_\R$-equivariant, $\strat$-constructible
sheaf.

\end{subsection}


\begin{subsection}{Perverse specialization}\label{ssperv}

Recall from Section~\ref{scompref} 
that $\raffgr^0$ is the union of certain components
of $\raffgr$ defined as follows. 
The exact sequence of component groups
$$
\pi_0(\affgr)\stackrel{}{\to} \pi_0(\raffgr)\stackrel{}{\to} \pi_0(G_\R),
$$
may be identified with the exact sequence
$$
\cowt/ Q\stackrel{\sigma}{\to} \rcowt/ (\rcowt\cap Q) 
\stackrel{\partial}{\to} \pi_0(G_\R)
$$
where $Q$ is the coroot lattice of $G$, $\sigma(\lambda)=\theta(\lambda)+\lambda$ is
the projection, 
and $\partial(\lambda)=[\lambda(-1)]$ is the boundary map.
By definition, $\raffgr^0$ is the union of the components
of $\raffgr$ in the image of $\sigma$, or equivalently in
the kernel of $\partial$.

\begin{Lemma}\label{lcomp}
For $\psh\in\epcat$, the support of $\on R(\psh)\in\redcat$ lies in $\raffgr^0.$
\end{Lemma}

\begin{proof}
Consider the real form $\affgr^{(\sigma)}_{T_\R}\subset\raffgr$
of the Beilinson-Drinfeld Grassmannian $\taffgr_T\subset\twoaffgr$
of the torus $T\subset G$.
By Proposition~\ref{prstrat}, it suffices to understand the limits of points
in the family $\affgr^{(\sigma)}_{T_\R}\to X$.
Note that the projection $\pi:\affgr^{(\sigma)}_{T_\R}\to X$
has discrete fibers.
It is easy to check that the limit of a point $\lambda\in\cowt\simeq
\affgr^{(\sigma)}_{T_R}|x_+,$ for $x_+\in X_+,$
is the point $\sigma(\lambda)\in\rcowt\simeq
\affgr^{(\sigma)}_{T_\R}|x,$ for $x\in X_\R$.
\end{proof}

By Lemma~\ref{lcompinc}, the components $\raffgr^0$ are in fact the union of certain components of $\raffgr^+$.
Since there is a perverse $t$-structure on the
category $\catd_{G_\R(\cor)}(\raffgr^+)$,
we may ask whether the specialization $\on R$
takes perverse
sheaves to perverse sheaves.
It is a fundamental result
that the nearby cycles in a complex
algebraic family are perverse. 
(See~\cite[Section 4.4]{BBD82} or \cite[Section 6.5]{GM83b}.)
In general this is not true for real algebraic families.  
The family $\sraffgr\to X$
is a good example of this: in general the specialization $\on R$ does not take
perverse sheaves to perverse sheaves.

\begin{Theorem}\label{tqs}
The specialization $\on R:\epcat\to\catd_{G_\R(\cor)}(\raffgr^+)$ 
is perverse if and only
if $G_\R$ is quasi-split.
\end{Theorem}

To prove the theorem we need to know more about
the specialization. The result is not used in what follows,
and we postpone the proof. See Corollary~\ref{cpiffqs}.

We define $\repcat_\Z$ to be the category $\repcat\otimes\vect_\Z$,
and identify it with the subcategory $\sum_k \repcat[k]$
of the derived category $\redcat$ from which
it inherits a convolution product.

We define the perverse specialization
$$
\pR:\epcat\to\catp_{G_\R(\cor)}(\raffgr^+)_\Z
$$
to be the sum 
$$
\pR= \sum_{k} \pH^k \circ\on R
$$
of the perverse homology sheaves of the specialization
where the $\Z$-grading corresponds to the degree of perverse homology.

%

\end{subsection}






%

%
\end{section}


\begin{section}{Monoidal structure for specialization}\label{smon}

The aim of this section is to 
construct an isomorphism 
$$
{}^p \hspace{-0em}r:\pR(\cdot\odot\cdot)\risom\pR(\cdot)\odot\pR(\cdot)
$$
for the perverse specialization 
$$\pR:\epcat\to\repcat_\Z.$$
We shall refer to such an isomorphism as a {\em monoidal structure}
for the functor.
Observe that the composition of two functors with monoidal structures
inherits a monoidal structure.


\begin{subsection}{Monoidal structure for local specialization $\Rightarrow$ 
monoidal structure for perverse specialization}\label{ssloctoperv}

Recall that the convolution 
$$
\odot :\catd_{G(\cor)}(\raffgr^+)\times\catd_{G(\cor)}(\raffgr^+)\to\catd_{G(\cor)}(\raffgr^+)
$$ 
is $t$-exact since
the multiplication map is stratified semismall.
This provides isomorphisms
$$
{}^p h^k:\pH^k(\cdot\odot\cdot)\risom\sum_{m+n=k}\pH^m(\cdot)\odot\pH^n(\cdot)
$$
for the perverse homology $\sum_k\pH^k:\catd_{G_\R(\cor)}(\raffgr^+)\to\catp_{G_\R(\cor)}(\raffgr^+)_\Z$.
%


If we have a monoidal structure
$$
r:\on R(\cdot\odot\cdot)\risom\on R(\cdot)\odot\on R(\cdot)
$$
for the local specialization
$\on R:\epcat\to\catd_{G_\R(\cor)}(\raffgr^+)$,
then we obtain a monoidal structure
$$
{}^p r:\pR(\cdot\odot\cdot)
\simeq \pR(\cdot)\odot\pR(\cdot)
$$
for the perverse specialization $\pR= \sum_{k} \pH^k \circ\on R:\epcat\to\repcat_\Z$.

\end{subsection}


\begin{subsection}{Global convolution}
Following~\cite[Section 5]{MV04},
we recall the global version of the convolution product.
For a $\C$-algebra $\alg$, and $x\in X(\alg)$, 
let $\widehat{(X_\alg)}_{x}$ denote the formal neighborhood of 
the graph of $x$ in the product $X_\alg=X\times\spec(\alg)$.

Consider the global convolution diagram
\begin{eqnarray*}\label{gcd}
\oneaffgr\times\oneaffgr \stackrel{p}{\gets}
\widetilde{\oneaffgr\times\oneaffgr} \stackrel{q}{\to}
\oneaffgr\tilde{\times}\oneaffgr \stackrel{m}{\to}
\twoaffgr\stackrel{d}{\gets}\oneaffgr.
\end{eqnarray*}
The ind-scheme
$\widetilde{\oneaffgr\times\oneaffgr}$ represents the functor
$\alg\mapsto
\{\left(
x_1,x_2, 
 \tors_1,\tors_2, 
\nu_1,
\nu_2,
\mu_1 
\right)\}
$
where for $i=1,2$, $x_i\in X(\alg)$, 
$\tors_i$ is a $G$-torsor on $X_\alg$, and 
$\nu_i$ is a trivialization of $\tors_i$ 
over $X_\alg\setminus x_i$, 
and
$\mu_1$ is a trivialization of $\tors_1$ over $\widehat{(X_\alg)}_{x_2}$.
The ind-scheme ${\oneaffgr\tilde\times\oneaffgr}$ represents the functor 
$\alg\mapsto
\{\left(
x_1, x_2,
\tors_1, \tors,
\nu_1, 
\eta 
\right)\}
$
where $x_1,x_2\in X(\alg)$,
$\tors_1,\tors$ are $G$-torsors on $X_\alg$,
$\nu_1$ is a trivialization of $\tors_1$ over $X_\alg\setminus x_1$,
and
$\eta$ is an isomorphism from $\tors_1$ to $\tors$ over $X_\alg\setminus x_2$.
The map $p$ forgets the trivialization $\mu_1$.
The map $q$ is given by
$
(\tors_1,\tors_2,\nu_1,\nu_2,\mu_1)\mapsto 
(\tors_1,\tors,\nu_1,\eta)
$
where $\tors$ is obtained by gluing $\tors_1$ over $X_\alg\setminus x_2$
with $\tors_2$ over $\widehat{(X_\alg)}_{x_2}$ via the 
isomorphism $\nu_2\circ\mu_1^{-1}|(X_\alg\setminus x_2)\cap\widehat{(X_\alg)}_{x_2}$.
The map $m$ is given by
$
(\tors_1,\tors,\nu_1,\eta)\mapsto
(\tors,\nu)
$
where $\nu=\eta\circ\nu_1|X\setminus(x_1\cup x_2)$.
The map $d$ is the inclusion
$
\oneaffgr\stackrel{\sim}{\to}\twoaffgr|\Delta\subset\twoaffgr
$
where $\Delta$ is the diagonal in $X^2$.

The global convolution
$$
\odot_X:\catp_{G,\hat\strat}(\oneaffgr)\times\catp_{G,\hat\strat}(\oneaffgr)\to
\catp_{G,\hat\strat}(\oneaffgr)$$
is defined to be
$$
\sh_1\odot_X\sh_2=d^*m_!(\sh_1\tilde\boxtimes\sh_2)
$$
where $\sh_1\tilde\boxtimes\sh_2$ is the unique sheaf such that
$
q^*(\sh_1\tilde\boxtimes\sh_2)=p^*(\sh_1\boxtimes\sh_2).
$

We are not immediately guaranteed that 
there is a unique sheaf 
$\sh_1\tilde\boxtimes\sh_2$ 
with the required property since $p$ and $q$ are not the projections
of torsors for a group-scheme 
which is the limit of linear algebraic groups. 
Fortunately, we may extend the diagram
to obtain maps with this property. Consider the diagram
$$
\begin{array}{ccccc}
& & \widetilde{\oneaffgr\times \hat X}\times\affgr & & \\
& \swarrow & \downarrow & \searrow & \\
\oneaffgr\times\oneaffgr & \stackrel{p}{\gets} &
\widetilde{\oneaffgr\times\oneaffgr} & \stackrel{q}{\to} &
\oneaffgr\tilde{\times}\oneaffgr. 
\end{array}
$$
The ind-scheme
$\widetilde{\oneaffgr\times \hat X}$ represents the functor 
$\alg\mapsto
\{\left(
x_1,x_2,
 \tors_1, \nu_1,
\mu_1,
\phi
\right)\}
$
where $x_1,x_2\in X(\alg),
 \tors_1$ is a $G$-torsor on $X_\alg$,
$\nu_1$ is a trivialization of $\tors_1$ 
   over $X_\alg \setminus x_1$,
$\mu_1$ is a trivialization of $\tors_1$ over $\widehat{(X_\alg)}_{x_2}$,
and
$\phi$ is an isomorphism $\Spec(\alg[ [t] ])\to\widehat{(X_\alg)}_{x_2}$.
The vertical map 
is the projection of an $\Aut(\co)$-torsor, and
the diagonal maps are the projections of
$G(\co)\rtimes\Aut(\co)$-torsors.
%

\end{subsection}


\begin{subsection}{Real global convolution}

Consider the real global convolution diagram
\begin{eqnarray*}\label{rgcd}
\oneraffgr\times\oneraffgr \stackrel{p}{\gets}
\widetilde{\oneraffgr\times\oneraffgr} \stackrel{q}{\to}
\oneraffgr\tilde{\times}\oneraffgr \stackrel{m}{\to}
\tworaffgr\stackrel{d}{\gets} \oneraffgr.
\end{eqnarray*}

The real global convolution
$$
\odot_{X_\R}:\catd_{G_\R,\hat\strat}(\oneraffgr)\times\catd_{G_\R,\hat\strat}(\oneraffgr)\to
\catd_{G_\R,\hat\strat}(\oneraffgr)
$$
is defined to be
$$
\sh_1\odot_X\sh_2=d^*m_!(\sh_1\tilde\boxtimes\sh_2)
$$
where $\sh_1\tilde\boxtimes\sh_2$ is the unique sheaf such that
$
q^*(\sh_1\tilde\boxtimes\sh_2)=p^*(\sh_1\boxtimes\sh_2).
$
Here we may use the extended diagram 
$$
\begin{array}{ccccc}
& & \widetilde{\oneraffgr\times \hat X_\R^0}\times\raffgr & & \\
& \swarrow & \downarrow & \searrow & \\
\oneraffgr\times\oneraffgr & \stackrel{p}{\gets} &
\widetilde{\oneraffgr\times\oneraffgr} & \stackrel{q}{\to} &
\oneraffgr\tilde{\times}\oneraffgr
\end{array}
$$
to see that there exists a unique sheaf $\sh_1\tilde\boxtimes\sh_2$
with the required property.

\end{subsection}


\begin{subsection}{Monoidal structure for global specialization $\Rightarrow$ 
Monoidal structure for local specialization}\label{ssgtoloc}

\begin{Lemma} 
There is a canonical isomorphism
$$
\rho(\cdot)\odot_X\rho(\cdot) \simeq  \rho(\cdot\odot\cdot)
$$ 
for $\rho:\epcat\to\catp_{G,\hat\strat}(\oneaffgr)$. 
\end{Lemma}

\begin{proof}
Consider the diagram
$$
\begin{array}{ccccccc}
\oneaffgr\times\oneaffgr & \stackrel{p}{\gets} &
\widetilde{\oneaffgr\times\oneaffgr} & \stackrel{q}{\to} &
\oneaffgr\tilde{\times}\oneaffgr & \stackrel{m}{\to} &
\twoaffgr   \\

d\uparrow & & d\uparrow & & d\uparrow & & d\uparrow   \\

\oneaffgr\times\oneaffgr|\Delta & \stackrel{p}{\gets} &
\widetilde{\oneaffgr\times\oneaffgr}|\Delta & \stackrel{q}{\to} &
\oneaffgr\tilde{\times}\oneaffgr|\Delta & \stackrel{m}{\to} &
\twoaffgr|\Delta   \\

s_\Delta \uparrow & &s_\Delta \uparrow & &s_\Delta \uparrow & & s_\Delta\uparrow   \\

\encurve\times \affgr\times\affgr & \stackrel{p}{\gets} &
\encurve \times G(\ck)\times \affgr & \stackrel{q}{\to} &
\encurve\times G(\ck)\times_{G(\co)}\affgr & \stackrel{m}{\to} &
\encurve\times \affgr 

\end{array}
$$
The first row is the global convolution diagram,
and the second is its restriction to the diagonal.
The third row is the product of the local convolution diagram
with the enhanced curve $\encurve$.
The arrows $d$ are the inclusions,
and the arrows $s_\Delta$ are the projections of $\aut(\co)$-torsors.

Let $r$ be the composition $d\circ s_\Delta$.
The proof 
is a diagram chase using the standard identity
$
r^*p^*\simeq p^*r^*$,
$
r^*q^*\simeq q^*r^*
$
for the composition of maps, 
the base change isomorphism 
$r^*m_!\simeq m_!r^*,$
and the canonical isomorphism
$
r^*(\rho(\cdot)\boxtimes\rho(\cdot))\simeq \C_{\encurve}\boxtimes(\cdot)\boxtimes(\cdot)
$
which results from a chase in the diagram
$$
\begin{array}{ccc}
\oneaffgr\times\oneaffgr|\Delta & \stackrel{d}{\to} & \oneaffgr\times\oneaffgr\\
s_\Delta \uparrow & & \uparrow s\times s \\
\encurve \times \affgr\times\affgr & \stackrel{d}{\to} & \encurve\times\encurve \times \affgr\times\affgr \\
 \downarrow & & \downarrow  \\
\affgr\times\affgr & = & \affgr\times\affgr.\\
\end{array}
$$
\end{proof}

\begin{Lemma} 
There is a canonical isomorphism
$$
\rho^{-1}_\R(\cdot)\odot_{}\rho^{-1}_\R(\cdot)\simeq  \rho^{-1}_\R(\cdot\odot_{X_\R}\cdot)
$$
for $\rho^{-1}_\R: \catd_{G_\R,\hat\strat}(\oneraffgr)\to\redcat$.
\end{Lemma}

\begin{proof}
Fix $x\in X_\R$, and a formal coordinate at $x$.
Consider the diagram
$$
\begin{array}{ccccccc}
\oneraffgr\times\oneraffgr & \stackrel{p}{\gets} &
\widetilde{\oneraffgr\times\oneraffgr} & \stackrel{q}{\to} &
\oneraffgr\tilde{\times}\oneraffgr & \stackrel{m}{\to} &
\tworaffgr   \\

r\uparrow & & r\uparrow & & r\uparrow & & r\uparrow   \\

\oneraffgr\times\oneraffgr|(x,x) & \stackrel{p}{\gets} &
\widetilde{\oneraffgr\times\oneraffgr}|(x,x) & \stackrel{q}{\to} &
\oneraffgr\tilde{\times}\oneraffgr|(x,x) & \stackrel{m}{\to} &
\tworaffgr|(x,x)   \\

\wr\uparrow & & \wr\uparrow & & \wr\uparrow & & \wr\uparrow   \\

\raffgr\times\raffgr & \stackrel{p}{\gets} &
G_\R(\ckr)\times \raffgr & \stackrel{q}{\to} &
\raffgr\tilde\times \raffgr & \stackrel{m}{\to} &
 \raffgr 

\end{array}
$$
The first row is the real global convolution diagram,
and the second is its restriction to $(x,x)\in X_\R^2$.
The third row is the real convolution diagram.
The arrows $r$ are the inclusions,
and the isomorphisms are provided by the formal coordinate at $x$.

The proof 
is a diagram chase using the standard identity
$
r^*p^*\simeq p^*r^*$,
$
r^*q^*\simeq q^*r^*
$
for the composition of maps, and
the base change isomorphism 
$r^*m_!\simeq m_!r^*.$
\end{proof}

Now if we have a monoidal structure
$$
{r}_X:\on {R}_X(\cdot\odot\cdot)\risom\on{R}_X(\cdot)\odot \on{R}_X(\cdot)
$$
for the global specialization
$\on{R}_X:\catd_{G,\hat\strat}(\oneaffgr)\to\catd_{G_\R,\hat\strat}(\raffgr)$,
then we obtain a monoidal structure
$$
r:\on R(\cdot\odot\cdot)
\risom
\on R(\cdot)\odot\on R(\cdot)
$$
for the local specialization $\on R=\rho_\R^{-1}\circ \on{R}_X \circ \rho:\epcat\to\redcat$.

\end{subsection}


\begin{subsection}{Monoidal structure for global specialization}\label{ssglobmon}

\begin{Proposition}\label{tglobmon}
There is a canonical isomorphism
$$
{r}_X: \on{R}_X(\cdot\odot_X\cdot)\simeq  \on{R}_X(\cdot)\odot_{X_\R}  \on{R}_X(\cdot)
$$
for the global specialization
$ \on{R}_X:\catp_{G, {\hat\strat}}(\oneaffgr)\to\catd_{\GR, {\hat\strat}}(\oneraffgr)$. 
\end{Proposition}

\begin{proof}
Consider the intertwining diagram
$$
\twoaffgr\times\twoaffgr \stackrel{p}{\gets}
\widetilde{\twoaffgr\times\twoaffgr} \stackrel{q}{\to}
\twoaffgr\tilde{\times}\twoaffgr \stackrel{m}{\to}
\affgr^{(4)}\stackrel{d}{\gets}\twoaffgr.
$$
The ind-scheme $\widetilde{\twoaffgr\times\twoaffgr}$
represents the functor 
$\alg\mapsto
\{\left(
x_1,x_2,x_3,x_4,
\tors_1,\tors_2,
\nu_1,
\nu_2,
\mu_1
\right)\}$
where
$x_1,x_2,x_3,x_4\in X(\alg),$
$\tors_1,\tors_2$ are $G$-torsors on $X_\alg$,
$\nu_1$ is a trivialization of $\tors_1$ over 
$X_\alg\setminus(x_1\cup x_{2})$,
$\nu_2$ is a trivialization of $\tors_2$ over 
$X_\alg\setminus(x_3\cup x_{4})$,
and
$\mu_1$ is a trivialization of $\tors_1$ over  
$\widehat{(X_\alg)}_{x_3\cup x_4}.$
The ind-scheme $\twoaffgr\tilde\times\twoaffgr$
represents the functor 
$\alg\mapsto
\{\left(
x_1,x_2,x_3,x_4,
\tors_1,\tors,
\nu_1, 
\eta
\right)\}$
where
$x_1,x_2,x_3,x_4\in X(\alg),$
$\tors_1,\tors$ are $G$-torsors on $X_\alg$,
$\nu_1$ is a trivialization of $\tors_1$ over 
$X_\alg\setminus(x_1\cup x_2),$
$\eta$ is an isomorphism from  
$\tors_1$ to $\tors$ over $X_\alg\setminus(x_3\cup x_4)$.
The maps are straightforward generalizations of the maps
in the global convolution diagram.


To construct the isomorphism of the proposition, we work with the 
real intertwining diagram
\begin{eqnarray*}\label{id}
\sraffgr\times\sraffgr \stackrel{p}{\gets}
\widetilde{\sraffgr\times\sraffgr} \stackrel{q}{\to}
\sraffgr\tilde{\times}\sraffgr \stackrel{m}{\to}
\stsraffgr\stackrel{d}{\gets}\sraffgr
\end{eqnarray*}
which is the real form of the intertwining diagram
with respect to the conjugation associated
to the product element 
$\sigma\times\sigma\in\Sigma_2\times\Sigma_2\subset\Sigma_4$.

To simplify the notation in what follows, write 
$
Y_+=X_+\times X_+$ and $Y_\R=X_\R\times X_\R$
for the products, and
$
Y=Y_+\cup Y_\R
$
for their union.
Let 
$\Delta\subset Y,$
$\Delta^+\subset Y^+,$ and 
$\Delta_\R\subset Y_\R$
denote the diagonals, and
$
j:Y^+\to Y$, $i: Y_\R\to Y$, and $d:\Delta\to Y$ 
the inclusions.

Observe that the restriction of the real intertwining diagram to $Y_\R$
is canonically identified with the real global
convolution diagram.
Thus we may use it to
calculate the real global convolution.
The restriction of the real intertwining 
diagram to $Y^+$ is canonically identified with
the restriction of two copies of the
global convolution diagram to $X_+$.
The following lemma confirms that we may use
this to calculate the global
convolution. 

\begin{Lemma}
There is a canonical isomorphism
$
(\cdot \odot_X \cdot)|X_+\simeq (\cdot|X_+)\odot_X (\cdot|X_+).
$


\end{Lemma}

\begin{proof}
We break the global convolution into two steps. 

First, consider the diagram
$$
\begin{array}{ccccc}

\oneaffgr\times\oneaffgr & \stackrel{p}{\gets} &
\widetilde{\oneaffgr\times\oneaffgr} & \stackrel{q}{\to} &
\oneaffgr\tilde{\times}\oneaffgr \\

j\uparrow && j\uparrow && j\uparrow \\

\oneaffgr\times\oneaffgr|Y^+ & \stackrel{p}{\gets} &
\widetilde{\oneaffgr\times\oneaffgr} |Y^+ & \stackrel{q}{\to} &
\oneaffgr\tilde{\times}\oneaffgr |Y^+ 

\end{array}
$$
We have the standard identity
$
j^*p^*\simeq p^*j^*,
$
$
j^*q^*\simeq q^*j^*
$
for the composition of maps.
It is a diagram chase to see that this provides an isomorphism 
$
(\cdot \tilde\boxtimes \cdot)|Y^+\simeq 
(\cdot|X_+)\tilde\boxtimes (\cdot|X_+).
$

Next, consider the diagram
$$
\begin{array}{ccccc}

\oneaffgr\tilde{\times}\oneaffgr & \stackrel{m}{\to} &
\twoaffgr & \stackrel{d}{\gets} & \oneaffgr \\

j\uparrow && j\uparrow && j\uparrow \\

\oneaffgr\tilde{\times}\oneaffgr |Y^+ & \stackrel{m}{\to} &
\twoaffgr|Y^+ & \stackrel{d}{\gets} & \oneaffgr|Y^+

\end{array}
$$
We have the base change isomorphism
$
j^*m_!\simeq m_!j^*,
$
and the standard identity
$
j^*d^*\simeq d^*j^*
$ 
for the composition of maps.
The required isomomorphism is the composition
$$
(\cdot \odot_X \cdot)|X_+=j^*d^*m_!(\cdot\tilde\boxtimes\cdot)\simeq d^*m_!j^*(\cdot\tilde\boxtimes\cdot)\simeq d^*m_!((\cdot|X_+)\tilde\boxtimes(\cdot|X_+))=
(\cdot|X_+)\odot_X (\cdot|X_+).
$$
\end{proof}

Now to construct the isomorphism of the proposition, 
we construct a series of isomorphisms 
relating global convolution to global specialization.
We begin with the left-hand portion of the real intertwining diagram
$$
\begin{array}{ccccc}

\sraffgr\times\sraffgr|Y^+ & \stackrel{p}{\gets} &
\widetilde{\sraffgr\times\sraffgr}|Y^+ & \stackrel{q}{\to} & 
\sraffgr\tilde{\times}\sraffgr|Y^+ \\

j\downarrow &  & j\downarrow &  & j\downarrow \\

\sraffgr\times\sraffgr|Y & \stackrel{p}{\gets} &
\widetilde{\sraffgr\times\sraffgr}|Y & \stackrel{q}{\to} & 
\sraffgr\tilde{\times}\sraffgr|Y \\

i\uparrow &  & i\uparrow &  & i\uparrow \\

\sraffgr\times\sraffgr|Y_\R & \stackrel{p}{\gets} &
\widetilde{\sraffgr\times\sraffgr}|Y_\R & \stackrel{q}{\to} &
\sraffgr\tilde{\times}\sraffgr|Y_\R 

\end{array}
$$ 

\begin{Lemma}
There is a canonical isomorphism 
$
i^*j_*(\cdot\tilde\boxtimes\cdot)\simeq (i^*j_*\cdot)\tilde\boxtimes (i^*j_*\cdot).
$

\end{Lemma}

\begin{proof}
This is a diagram chase using the 
standard identity
$i^*p^*\simeq p^*i^*$,
$i^*q^*\simeq q^*i^*$
for the composition of maps,
and the smooth base change isomorphism
$j_*p^*\simeq p^*j_*$, $j_*q^*\simeq q^*j_*.$
The maps $p$ and $q$ are smooth since
they are the projections of 
torsors for the smooth relative group-scheme 
$\rgrsc\to X_R^{(\sigma)}$ which
is a real form of the smooth relative group-scheme
$\grsc\to X^2$ that
represents the functor
$\alg\mapsto
\{\left(
x_1,x_2,
\mu\right)\}
$
where
$x_1,x_2\in X(\alg)$
and
$\mu$ is a trivialization of the trivial $G$-torsor 
over $\widehat{(X_\alg)}_{(x_1\cup x_2)}$.
\end{proof}

Next, consider the middle portion of the real intertwining diagram
$$
\begin{array}{ccc}

\sraffgr\tilde{\times}\sraffgr|Y^+ & \stackrel{m}{\to} &
\stsraffgr|Y^+ \\

j\downarrow &  & j\downarrow \\

\sraffgr\tilde{\times}\sraffgr|Y & \stackrel{m}{\to} &
\stsraffgr|Y \\

i\uparrow &  & i\uparrow \\

\sraffgr\tilde{\times}\sraffgr|Y_\R & \stackrel{m}{\to} &
\stsraffgr|Y_\R
\end{array}
$$ 

\begin{Lemma}
There is a canonical isomorphism 
$
i^*j_*m_!\simeq m_! i^*j_*.
$

\end{Lemma}

\begin{proof}
This follows from the base change isomorphism
$m_!i^*\simeq i^*m_!,$
the standard identity $m_*j_*\simeq j_*m_*$ for the composition of maps,
and the fact that $m_*=m_!$ for the sheaves under consideration
since $m$ is proper on their supports.  
\end{proof}

Finally, consider the right-hand portion of the 
real intertwining diagram





$$
\begin{array}{ccc}
\stsraffgr|Y^+ & \stackrel{d}{\gets} & \stsraffgr|\Delta^+ \\

 j\downarrow & & j\downarrow  \\

\stsraffgr|Y^2 & \stackrel{d}{\gets} & \stsraffgr|\Delta \\

i\uparrow &  & i\uparrow  \\

\stsraffgr|Y_\R & \stackrel{d}{\gets} & \stsraffgr|\Delta_\R  \\

\end{array}
$$ 

\begin{Lemma}
There is a canonical isomorphism 
$
d^* i^*j_*\simeq i^*j_*d^*.
$

\end{Lemma}

\begin{proof}
To obtain the candidate for the isomorphism,
we apply $d^*j_*$ to the adjunction morphism
$\id\to d_*d^*$. We obtain the morphism
$
d^*j_*\to d^*j_*d_*d^*\stackrel{\sim}{\to}d^*d_*j_*d^*
\stackrel{\sim}{\to} j_*d^*,$
where the first isomorphism is the standard identity for the composition
of maps, and the second follows from the fact that 
the adjunction morphism $d^*d_*\to\id$ is an isomorphism
since $d$ is an embedding.
Now applying $i^*$ to the above composite morphism, we obtain
$$
d^*i^*j_*\stackrel{\sim}{\to} i^*d^*j_*\to i^*j_*d^*,
$$
where the initial isomorphism is the standard identity
for the composition of maps. 
We shall prove this morphism is an isomorphism.

Write $W^+=Y^+\setminus\Delta^+$ for the complement,
and 
$k:\stsraffgr|W^+\to \stsraffgr|Y^+$
for the inclusion.
The adjunction morphism $\id\to d_*d^*$ fits into a distinguished triangle
$
k_!k^{*}\to \id\to d_*d^*\stackrel{[1]}{\to}.
$
Applying $i^*d^*j_*$ and using the isomorphism $d^*j_*d_*d^*\stackrel{\sim}{\to} j_*d^*$,
we obtain the distinguished triangle
$
i^*d^*j_*k_!k^*\to i^*d^*j_*\to i^*j_*d^*\stackrel{[1]}{\to}. 
$
We shall prove 
$$i^*d^*j_*k_!k^*=0,$$ 
which immediately implies the composition
$d^*i^*j_*\stackrel{\sim}{\to} i^*d^*j_*\to i^*j_*d^*$ is an isomorphism.

To establish this, we show that for any sheaf 
$\sh \in\catd_{G_\R,\strat}(\stsraffgr|Y^+)$
the stalk cohomology vanishes
$
H(j_*k_!k^*\sh)_x =0,$ 
for $x\in\stsraffgr|\Delta_\R.$  
To be specific, we show that for 
any open
neighborhood $N\subset\stsraffgr|Y$ 
containing $x$
there is an open neighborhood $N'\subset N$ containing $x$ such that
the hypercohomology vanishes
$$
\hc(N',j_*k_!k^*\sh)=0.
$$
For the remainder of the proof,
we fix a point $x\in\stsraffgr|\Delta_\R$ and open 
neighborhood $N\subset\stsraffgr|Y$ containing $x$.

\begin{Lemma}
There is a
closed ball $B\subset N$ with $x\in B^o$, and an open ball 
$D\subset Y$ with $\pi(x)\in D$ such that 
$\pi:B|D\to D$
is a proper, stratified submersion of Whitney stratified sets.
\end{Lemma}

\begin{proof}
The argument is standard in the theory of nearby cycles.
The essential point is that the projection 
$\pi:\stsraffgr|Y\to Y$ is Thom stratified by Proposition~\ref{pthom}.
This implies that if $B$ is chosen so that its boundary $\partial B$
is transverse to all the strata in $\stsraffgr|X_\R$,
then $\partial B$ will also be 
transverse to all the strata in the nearby fibers of $\pi$.
Thus one can find $D$ such that 
$B|D$ is Whitney stratified by the restriction of the strata
of $\stsraffgr$.
\end{proof}

Fix a closed ball $B$ and an open ball $D$ as provided by the lemma.
We first assert that it is possible to choose
an open ball $D'\subset D$ containing $x$
so that there is a complete
vector field on $W^+\cap D'$ whose flow $f_t$ satisifes
$
\lim_{t\to0} f_t(w)\in\Delta^+
$
and 
$
\lim_{t\to\infty} f_t(w)\not \in\Delta^+,
$
for all $w\in W^+\cap D'$.
This is not difficult to see by examining the stratification
of $Y$.

We next assert that for the neighborhood $N'=B^o|D'$, we have the vanishing
$
\hc(N',j_*k_!k^*\sh)=0.
$
To see this, first note
$
\hc(N',j_*k_!k^*\sh)\simeq \hc(N'|Y^+,k_!k^*\sh).
$
Now by the lemma and the first Thom-Mather isotopy lemma (see \cite{Math70}
or \cite[Section 1.5]{GM88}),
we may lift the flow $f_t$ on $W^+\cap D'$ to a flow $F_t$ on $B|(W^+\cap D')$.
Then by restriction we obtain a lift to $N'|W^+$.
It follows that 
$
\hc(N'|Y^+,k_!k^*\sh)=0,
$
since any cycle in this group must be supported away from 
$N'|\Delta^+$, and hence may be realized as a boundary using 
the flow $F_t$ on $N'|W^+$.
\end{proof}

Finally, to complete the proof of the proposition, the monoidal structure 
$$
{r}_X: \on{R}_X(\cdot\odot_X\cdot)\simeq  \on{R}_X(\cdot)\odot_{X_\R}  \on{R}_X(\cdot)
$$
for the global specialization
$ \on{R}_X:\catp_{G, {\hat\strat}}(\oneaffgr)\to\catd_{\GR, {\hat\strat}}(\oneraffgr)$
is obtained from the isomorphisms of the lemmas via a diagram chase.
\end{proof}

\end{subsection}

\end{section}



\begin{section}{Image category I}\label{sim1}

%


First, we define $\catq_{\pR}$ to be the strict full subcategory of $\repcat_\Z$
whose objects are isomorphic to subquotients of objects
of the form $\pR(\psh)$, where $\psh$ runs through all objects of $\epcat$,
and $\pR$ is the perverse specialization
$$
\pR:\epcat\to\repcat_\Z.
$$

Next, we define $\catq$ to be the strict full subcategory of $\repcat$
whose objects are isomorphic to objects of the form $\ff(\CQ)$,
where $\CQ$ runs
through all objects of $\catq_{\pR}$,
and $\ff$ is the forgetful functor
$$
\ff:\repcat_\Z\to\repcat.
$$ 

Finally, we define
$\catq_\Z$ to be  the category $\catq\otimes\vect_\Z$,
and identify it with the strict full subcategory $\sum_k\catq[k]$ of the category $\repcat_\Z$.

Each of the three categories $\catq,\catq_{\pR},$ and $\catq_\Z$
inherits a convolution product from $\repcat_\Z$.
We may organize the categories into a commutative diagram of functors with monoidal structure
$$
\begin{array}{ccccc}
& & \catq_\pR & & \\
& & \ef\downarrow & \stackrel{\ff}{\searrow}  & \\
\catq & \stackrel{\ef}{\to} & \catq_\Z & \stackrel{\ff}{\to} & \catq. \\
\end{array}
$$
Here the arrows labelled $\ef$ are the obvious fully faithful functors,
and the arrows labelled $\ff$ are the obvious forgetful functors.

%


%

\end{section}


\begin{section}{Character functors}\label{schar}

In this section, we discuss the weight functors first introduced 
in~\cite{MV00} and studied in~\cite{MV04}. See also \cite{NP01}.

%

\begin{subsection}{Semi-infinite orbits}

Let $U$ be the unipotent radical of the Borel subgroup $B\subset G$.

The group $U(\ck)$ 
acts on the loop Grassmannian $\affgr$ by left-multiplication.
Recall that each coweight $\nu\in\cowt$ defines a point $\nu\in\Gr$.
Let $\mv$ be the $U(\ck)$-orbit $U(\ck)\cdot \nu \subset\affgr$
through $\nu$. 
By the Iwasawa decomposition~\cite[Proposition 2.33]{IM65}, 
each $U(\ck)$-orbit in $\affgr$
is of the form $\mv$, for some coweight $\nu\in\cowt$,
and the orbits are distinct for distinct coweights.
For the rest of the following, see~\cite[Proposition 3.1]{MV04}.

\begin{Proposition}\label{pmv}
The loop Grassmannian $\affgr$ is the disjoint union
of the $U(\ck)$-orbits $\mv$ through $\nu\in\cowt$.
The closure of the orbit $\mv$ 
is the union of the orbits $S^\mu$ with
$\mu\leq\nu$.
\end{Proposition}

Note that there is a similar statement for
the orbits $T^\nu=U^o(\CK)\cdot\nu\subset\affgr$,
where $U^o\subset G$ is the unipotent subgroup opposite to $U$.

\end{subsection}


\begin{subsection}{Real semi-infinite orbits}

Let $U_{P_\R}$ be the unipotent radical of the minimal
parabolic subgroup $P_\R\subset G_\R$.

The group $U_{P_\R}(\ckr)$ 
acts on the real loop Grassmannian $\raffgr$ by left-multiplication.
Recall that each real coweight $\nu\in\rcowt$ defines a point $\nu\in\raffgr$.
Let $\rmv$ be the $U_{P_\R}(\ckr)$-orbit $U_{P_\R}(\ckr)\cdot \nu \subset\raffgr$
through $\nu$. 
%


%


\begin{Lemma}\label{laltrmvdef}
For $\nu\in\rcowt$, the $U_{P_\R}(\ckr)$-orbit $\rmv$ equals the intersection $\mv\cap \raffgr$,
and for $\nu\in\cowt\setminus\rcowt$, the intersection $\mv\cap \raffgr$ is empty.
\end{Lemma}

\begin{proof}
We show that if a $U_P(\CK)$-orbit in $\affgr$ intersects $\raffgr$, then it must
contain an element of $\rcowt$.
By Proposition~\ref{pmv}, the group $P(\CK)$ acts transitively on $\affgr$,
and thus each $U_P(\CK)$-orbit intersects the loop Grassmannian $\affgr_M$
of the Levi factor $M\subset P$. Therefore it suffices to show that 
the real loop Grassmannian of the Levi factor ${M_\BR}\subset P_\BR$ is equal to $\rcowt$,
or equivalently, that the real loop Grassmannian of 
the derived group ${M'_\BR}\subset M_\BR$ is equal to a single point.
Since $M'_\BR$ is compact, this follows immediately from Proposition~\ref{prstrat}.
\end{proof}

We have the following Iwasawa decomposition. 

\begin{Proposition}\label{prmv}
The real loop Grassmannian $\raffgr$ is the disjoint union
of the $U_{P_\R}(\ckr)$-orbits $\rmv$ through $\nu\in\rcowt$.
The closure of the orbit $\rmv$ 
is the union of the orbits $S_\R^\mu$ with
$\mu\leq\nu$. 
\end{Proposition}

\begin{proof}
Thanks to the previous lemma and Proposition~\ref{pmv},
it only remains to prove that
if $\nu-\mu$ is a non-negative integral
linear combination of positive coroots of $G$,
then $S_\R^\mu\subset\overline S_\R^\nu$.
The action on $\raffgr$ of elements of the real coweight lattice $\rcowt$ translate the orbits,
so we may assume that $\nu$ is the trivial coweight $0$.
By Lemma~\ref{lrcowtord}, it suffices to prove the assertion when $-\mu$
is a positive coroot $\alpha \in R^\pos_{}$, or when $-\mu$ is of the form
$\theta(\alpha)+\alpha$, for a positive coroot $\alpha\in R^\pos_{}$,
but $-\mu$ is not a multiple of a positive coroot.
In the first case, we may find $\on{SL}_2(\R)\subset G_\R$ such that $\alpha$
is its positive coroot. Then the orbit through $0$ of the one parameter subgroup
$U_{\alpha}(rt^{-1})\subset U_{P_\R}(\ckr)$, for $r\in \R$,
is isomorphic to $\R$, and its closure is isomorphic to $\R\BP^1$
with $\mu$ the point at infinity. 
In the second case, we may find $\on{SL}_2(\C)\subset G_\R$ such that $\theta(\alpha)+\alpha$
is its positive coroot. Then the orbit through $0$ of the one parameter subgroup
$U_{\theta(\alpha)+\alpha}(ct^{-1})\subset U_{P_\R}(\ckr)$, for $c\in\C$,
is isomorphic to $\C$, and its closure is isomorphic to $\C\BP^1$
with $\mu$ the point at infinity.
\end{proof}


Note that there are similar results for
the orbits $T_\R^\nu=U_{P^o_\R}(\ckr)\cdot\nu\subset\raffgr$,
where $U_{P^o_\R}\subset G_\R$ is the unipotent radical of
the parabolic subgroup $P^o_\R\subset G_\R$ opposite to ${P_\R}$.

\end{subsection}


\begin{subsection}{Flows}

%

%

%

%

For $\xi\in\fg$, we have the flow 
$$
\varphi_\xi:\R\times\affgr\to\affgr 
$$
$$
\varphi_\xi^t(x)=\exp(t\xi)\cdot x 
$$
generated by $\xi$ via the action of $G$ on $\Gr$.

Choose a Cartan involution of $G$ which commutes with the conjugation $\theta$,
and let $G_c\subset G$ be the resulting maximal compact subgroup.

For the following, see~\cite[Sections 3 and 6]{MV04} or~\cite[Section 4]{NP01} for more details.

\begin{Lemma}
Let $\xi\in\ft\cap i\fg_c$ be in the interior of the dominant Weyl chamber.
Then the fixed points in $\affgr$ of the flow $\varphi_\xi$  are the coweights $\nu\in\cowt$,
and the orbits $\mv$ are the ascending sets
$$
\mv=\{x\in \affgr | \lim_{t\to+\infty}\varphi_\xi^t(x)=\nu\},
$$
and the orbits $T^\nu$ are the descending sets
$$
T^\nu=\{x\in \affgr | \lim_{t\to-\infty} \varphi_\xi^t(x)=\nu\}.
$$
\end{Lemma}

\begin{proof}
The torus $T$ centralizes $\xi$, and
the weights of $\xi$ in $\fu_{}$ are positive,
and the weights in $\fu_{}$ are negative.
\end{proof}

For $\xi\in\fg_\R$, the flow 
generated by $\xi$ preserves the real form $\raffgr$ since the action of $G_\R$
preserves $\raffgr$. 
Recall that the torus $S_\R\subset G_\R$ is split, so that $\fs_\R\subset i\fg_c$.

\begin{Lemma}\label{lrflow}
Let $\xi\in\fs_\R$ be
in the interior of the dominant Weyl chamber for $\CW_{G_\R}$.
Then the fixed points in $\raffgr$ of the flow $\varphi_\xi$ are the real coweights
$\nu\in\rcowt$, and
the orbits $\rmv$ are the ascending sets
$$
\rmv=\{x\in \raffgr | \lim_{t\to+\infty}\varphi^t_\xi(x)=\nu\},
$$
and the orbits $T^\nu_\R$ are the descending sets
$$
T^\nu_\R=\{x\in \raffgr | \lim_{t\to-\infty}\varphi^t_\xi(x)=\nu\}.
$$
\end{Lemma}

\begin{proof}
The torus $S_\R$ centralizes $\xi$, and
the weights of $\xi$ in $\fu_{P_\R}$ are positive,
and the weights in $\fu_{P^o_\R}$ are negative.
\end{proof}

Let $P'\subset P$ be the derived group of the parabolic subgroup $P\subset G$.
Let $Q_M\subset\cowt$ be the coroot lattice of the Levi factor $M\subset P$.

For $\bar\nu\in\cowt/Q_M$, choose $\nu\in\cowt$ projecting to $\bar\nu$.
Let $S^{\bar\nu}_P$ be the connected component of the 
orbit $P'(\ck)\cdot \nu \subset\affgr$
through $\nu$,
and let $\affgr_{M,\bar\nu}$ 
be the connected component 
of the loop Grassmannian $\affgr_M$ 
containing $\nu$. 
The notation is justified by the fact that
for $\bar\nu\in\cowt/Q_M$, and $\nu\in\cowt$ projecting to $\bar\nu$,
 the connected component $S^{\bar\nu}_P$ is the disjoint union
of the $U(\ck)$-orbits $S^\mu$, for $\mu=\nu\mod Q_M$,
and the connected component $\affgr_{M,\bar\nu}$
contains $\mu\in\cowt$ if and only if $\mu=\nu\mod Q_M$.

Similarly, for $\bar\nu\in\cowt/Q_M$, we have the connected components
$T^{\bar\nu}_{P}$ of the 
orbits ${P^o}'(\ck)\cdot \nu \subset\affgr$
where $P^o\subset G$ is the parabolic subgroup opposite to $P$,
and ${P^o}'$ is its derived group.

\begin{Lemma}\label{lpflow}
Let $\xi\in\fs_\R$ be
in the interior of the dominant Weyl chamber for $\CW_{G_\R}$.
Then the fixed points in $\affgr$ of the flow $\varphi_\xi$
are the components $\affgr_{M,\bar\nu}$,
and the components $S^{\bar\nu}_P$ are the ascending sets
$$
S^{\bar\nu}_P=\{x\in \affgr | \lim_{t\to+\infty}\varphi^t_\xi(x)\in\affgr_{M,{\bar\nu}}\},
$$
and the conponents $T^{\bar\nu}_P$ are the descending sets
$$
T^{\bar\nu}_P=\{x\in \raffgr | \lim_{t\to-\infty}\varphi^t_\xi(x)\in\affgr_{M,{\bar\nu}}\}.
$$
\end{Lemma}

\begin{proof}
The Levi factor $M$ centralizes $\xi$,
the weights of $\xi$ in $\fu_{P}$ are positive,
and the weights in $\fu_{P^o}$ are negative.
\end{proof}

\end{subsection}


\begin{subsection}{Weight functors}\label{sslocchar}


Although in general the $U(\CK)$-orbits $\mv\subset\affgr$ are neither finite-dimensional
nor finite-codimensional, their intersections with the strata
provide perverse cell decompositions of the strata. 
The results of this section are all from~\cite[Section 3]{MV04}.

\begin{Proposition}[\cite{MV04}, Theorem 3.2]\label{pdimest}
For $\lambda\in\dcowt$, and $\nu\in\cowt$, the orbit $\mv$
meets the stratum $\affgr^\lambda$ if and only if $\nu\in\overline\affgr^\lambda$.
In this case, 
we have
$$
\dim_\C(\affgr^\lambda\cap \mv)=\langle\ch\rho,\lambda+\nu\rangle.
$$
\end{Proposition}

This and the subsequent proposition are the primary ingredients 
in the proof of the following. 

\begin{Theorem}[\cite{MV04}, Theorem 3.5]\label{twt}
For $\nu\in\cowt$, and $\psh\in\epcat$,
we have the vanishing
$$
H^k_c(\mv,\psh)=0 \mbox{ if } k\neq 2\langle\ch\rho,\nu\rangle.
$$
\end{Theorem}


\begin{Proposition}[\cite{MV04}, Theorem 3.5]\label{plocdual}
For $\psh\in\epcat$, there is a canonical isomorphism
$$
H^k_c(\mv,\psh	)\simeq H^k_{T^\nu}(\affgr,\psh),
\mbox{ for all } k.
$$
\end{Proposition}

For $\nu\in\cowt$, we call the functor
\begin{eqnarray*}
&F^\nu:\epcat\to\vect & \\
& F^\nu(\cdot)=H^{2\langle \ch\rho,\nu\rangle}_c(\mv,\cdot) &
\end{eqnarray*}
a weight functor.

\begin{Corollary}[\cite{MV04}, Theorem 3.5]
For $\nu\in\cowt$, the weight functor $F^\nu:\epcat\to\vect$ is exact.
\end{Corollary}

Recall that $\vect_\cowt$ is the category of finite-dimensional
$\cowt$-graded vector spaces.
We collect the weight functors into a functor
\begin{eqnarray*}
& \tchar:\edcat\to\vect_\cowt & \\
& \tchar=\sum_{\nu\in\cowt} F^\nu &
\end{eqnarray*}
we call the character functor.

%

Recall that $\hc$ is the hypercohomology functor,
and $\ff$ the forgetful functor which
forgets the grading of a vector space. 

\begin{Corollary}[\cite{MV04}, Theorem 3.6]\label{cevenpar}
There is a canonical isomorphism
$$
\hc\simeq \ff\circ\tchar 
:\epcat\to\vect.
$$
More generally,
for any set $Y\subset\cowt$ of coweights,
there is a canonical isomorphism
$$
H_c^*(\bigsqcup_{\nu\in Y}\mv,\cdot)\simeq \sum_{\nu\in Y} H_c^*(\mv,\cdot):\epcat\to\vect.
$$

\end{Corollary}

\begin{Corollary}[\cite{MV04}, Corollary 3.7]
The hypercohomology functor $\hc:\epcat\to\vect$ is exact and faithful.
\end{Corollary}

\end{subsection}


\begin{subsection}{Real weight functors}


In the real case, the basic dimension estimate takes the following form.

\begin{Proposition}\label{prdimest}
For $\lambda\in\rdcowt$, and $\nu\in\rcowt$, the orbit $\rmv$
meets the stratum $\affgr^\lambda_\R$ if and only if 
$\nu\in\overline\affgr^\lambda_\R$.
In this case, 
we have
$$
\dim_\R(\affgr^\lambda_\R\cap \rmv)\leq \langle\ch\rho,\lambda+\nu\rangle.
$$
\end{Proposition}

\begin{proof}
By Lemma~\ref{laltrmvdef}, $\rmv$ equals the intersection $\mv\cap\raffgr$,
and so $\rmv\cap\raffgr^\lambda$ equals the intersection $(\mv\cap\affgr^\lambda)\cap\raffgr$.
Thus the result follows immediately
from Proposition~\ref{pdimest}.
\end{proof}

\begin{Theorem}\label{trwt}
For $\nu\in\rcowt$, and $\psh\in\repcat$,
we have the vanishing
$$
H^k_c(\rmv,\psh)=0 \mbox{ if } k\neq \langle\ch\rho,\nu\rangle.
$$
\end{Theorem}

\begin{proof}
The dimension estimate implies the vanishings 
$$
H^k_c(\rmv,\psh)=0 \mbox{ for } k> \langle\ch\rho,\nu\rangle,
$$
$$
H^k_{T^\nu_\R}(\raffgr,\psh)=0 \mbox{ for } k< \langle\ch\rho,\nu\rangle.
$$
Therefore the following proposition implies the theorem.
\end{proof}

\begin{Proposition}\label{prhyp}
For $\sh\in\redcat$, there is a canonical isomorphism
$$
H^k_c(\rmv,\sh)\simeq H^k_{T^\nu_\R}(\raffgr,\sh), \mbox{ for all } k.
$$
\end{Proposition}

\begin{proof}
For $\xi\in \fs_\R$ in the interior of the dominant Weyl chamber for $\CW_{G_\R}$,
the map $\varphi_\xi:\raffgr\to\raffgr$ defined by $\varphi_\xi(x)=\exp(\xi)\cdot x$
is weakly hyperbolic in the sense of \cite[Section 1.3]{GM93}. Therefore by Lemma~\ref{lrflow}, and 
\cite[Proposition 9.2]{GM93},
we have the asserted isomorphism.
\end{proof}

For $\nu\in\rcowt$, $z\in\Z$, we call the functor
\begin{eqnarray*}
& F^{\nu,z}_\R:\redcat\to\vect & \\
& F^{\nu,z}_\R(\cdot)=H^{\langle\ch\rho,\nu\rangle+z}_c({\rmv},\cdot) &
\end{eqnarray*}
a real weight functor. 
When $F^{\nu,z}_\R$ is applied to perverse sheaves,
the integer $z$ becomes superfluous by the theorem. We sometimes then 
forget it and simply write $F^\nu_\R:\repcat\to\vect_\rcowt$ for the functor
$F^\nu_\R(\cdot)=H^*_c(\rmv,\cdot)$.

\begin{Corollary}\label{crwtexact}
For $\nu\in\rcowt$, the weight functor $F^{\nu}_\R:\repcat\to\vect$
is exact.
\end{Corollary}

We collect the real weight functors into a functor
\begin{eqnarray*}
& \schar:\redcat\to\vect_{\rcowt\times\Z} & \\
& \schar=\sum_{\stackrel{\nu\in\rcowt}{z\in\Z}} F_\R^{\nu,z} &
\end{eqnarray*}
we call the real character functor.

\begin{Remark}
In the real setting, the critical degree $k$ for the cohomology
group $H^k_c(\rmv,\psh)$, for $\psh\in\repcat$, is 
half of what it is in the complex case.
In particular, it is not always of even parity.
\end{Remark}

\end{subsection}


\begin{subsection}{Two useful facts}

In contrast to the complex case,
the intersections of the $U_{P_\R}(\ckr)$-orbits $\rmv\subset\raffgr$
with the strata are in general not
of pure dimension.
In addition, in some cases their dimension is strictly less than 
the upper bound seen in Proposition~\ref{prdimest}.
We note here the further vanishing this implies.

\begin{Proposition}\label{pextravanish}
For $\lambda\in\rdcowt$, and $\nu\in\rcowt$, if we have
$$
\dim_\R(\affgr^\lambda_\R\cap \rmv)<\langle\ch\rho,\lambda+\nu\rangle,
$$ 
then for $\psh\in\repcat$, we have the vanishing
$$
H^k_c(\rmv,\psh)=0 \mbox{ for all } k.
$$
\end{Proposition}

\begin{proof}
The dimension assumption implies the vanishing 
$$
H^k_c(\rmv,\psh)=0 \mbox{ for } k\geq \langle\ch\rho,\nu\rangle,
$$
and the dimension estimate of Proposition~\ref{prdimest} implies the vanishing
$$
H^k_{T^\nu_\R}(\raffgr,\psh)=0 \mbox{ for } k< \langle\ch\rho,\nu\rangle,
$$
so by Proposition~\ref{prhyp}, we have the asserted vanishing.
\end{proof}

We have the following symmetry of the real weight functors 
analogous to \cite[Section 5]{NP01}.
Recall that the small Weyl group $\CW_{G_\R}$ acts on the real coweights $\rcowt$.

\begin{Proposition}\label{pweylsym}
For $\nu\in\rdcowt$, and $w\in \CW_{G_\R}$,
we have
$$
F^{w(\nu)}_\R(\ic^\nu)\simeq \C 
$$
where $\ic^\nu$ is the intersection cohomology sheaf of the closure
of the stratum $\raffgr^\nu$.
\end{Proposition}

\begin{proof}
The intersection $\raffgr^\nu\cap S_\R^{\omega(\nu)}$ is a real affine space,
and by Proposition~\ref{prmv}, the intersection $\raffgr^\mu\cap S_\R^{\omega(\nu)}$ 
is empty for any $\mu$ in the closure of $\raffgr^\nu$.
\end{proof}

\end{subsection}


\subsection{Global preliminaries}\label{ssglobchar}

For $x\in X$, and a choice of formal coordinate at $x$,
the $U(\ck_x)$-orbits in the local loop Grassmannian
$\xaffgr$ are taken to the
$U(\ck)$-orbits in $\affgr$
under the isomorphism
$
\xaffgr\stackrel{\sim}{\to}\affgr.
$
Since the $U(\ck)$-orbits in $\affgr$ are $\Aut(\co)$-invariant,
we may unambiguously index the $U(\ck_x)$-orbits ${}_x\mv\subset \xaffgr$ by coweights
$\nu\in\cowt$. 

Similarly, for $x\in X_\R$, 
we may unambiguously index
the $U_{P_\R}({\ck_\R}_x)$-orbits ${}_x\rmv\subset \xraffgr$ by real coweights
$\nu\in\rcowt$.

Recall the factorization of the real form $\wraffgr$ with respect
to the projection 
$$\pi:\wraffgr\to X^{(\omega)}.$$
Fixing an identification $X^{(\omega)}\simeq X_\R^r\times X^s$, for 
$(x_1,\ldots,x_r,z_1,\ldots,z_s)\in X^r_\R\times {X}^{s},$
we have a canonical isomorphism
$$
\wraffgr|(x_1,\ldots,x_r,z_1,\ldots,z_s)\stackrel{\sim}{\to}  
\prod_{i=1}^k {{}_{y_i}\hspace{-0.3em}\affgr_{\R}}\times
\prod_{j=1}^\ell {}_{w_j}\hspace{-0.3em}\affgr.
$$
where $y_1,\ldots,y_k\in X_\R, w_1,\ldots,w_\ell\in X$ are the distinct points with an
equality of sets
$$\{y_1,\ldots,y_k,w_1,\ldots, w_\ell\}=\{x_1,\ldots,x_r,z_1,\ldots,z_s\}.$$

For $\nu\in\rcowt$, let $S^{(\omega)}_\nu\subset\wraffgr$ be the union over all fibers 
of the subspaces
$$
\prod_{i=1}^k {}_{y_i}{S^{\nu_i}_\R}\times\prod_{j=1}^\ell {}_{w_j}{S^{\mu_j}}
\subset\prod_{i=1}^k {}_{y_i}\hspace{-0.15em}{\affgr_{\R}}\times
\prod_{j=1}^\ell {}_{w_j}\hspace{-0.25em}\affgr
\mbox{ with }
\nu=\sum_{i=1}^k\nu_i+\sum_{j=1}^\ell \sigma(\mu_j),
$$  
and let 
$s_\nu:S_{\nu}^{(\omega)}\to \wraffgr$
be the inclusion. 
Here as usual $\sigma:\cowt\to\rcowt$ 
is the projection
$\sigma(\lambda)=\theta(\lambda)+\lambda$.

Similarly, 
for $\nu\in\rcowt$,
let $T^{(\omega)}_\nu\subset\wraffgr$ be the union over all fibers 
of the subspaces
$$
\prod_{i=1}^k {}_{y_i}{T^{\nu_i}_\R}\times\prod_{j=1}^\ell {}_{w_j}{T^{\mu_j}}
\subset\prod_{i=1}^k {}_{y_i}\hspace{-0.15em}{\affgr_{\R}}\times
\prod_{j=1}^\ell {}_{w_j}\hspace{-0.25em}\affgr
\mbox{ with }
\nu=\sum_{i=1}^k\nu_i+\sum_{j=1}^\ell \sigma(\mu_j),
$$  
and let
$t_\nu:
T_{\nu}^{(\omega)}\to \wraffgr$ 
be the inclusion.

Since the conjugation $\theta$ restricts to a Cartan involution of the derived group $M'$
of the Levi factor $M\subset P$, 
each connected component of the real form $\affgr_{M,\R}$
is a single real coweight $\nu\in\rcowt$.
Note that the connected components of the real form $\affgr^{(\omega)}_{M,\R}$
may be identified with the connected components of $\affgr_{M,\R}$
by collecting together points in the base $X^{(\omega)}$. 
Therefore each connected component $\affgr^{(\omega)}_{M,\nu,\R}$
of the real form
$\affgr^{(\omega)}_{M,\R}$ may be indexed by $\nu\in\rcowt$.

\begin{Lemma}
Let $\xi\in\fs_\R$ be
in the interior of the dominant Weyl chamber for $\CW_{G_\R}$.
Then the fixed points in $\wraffgr$ of the flow $\varphi_\xi$ are the components
$\affgr^{(\omega)}_{M,\nu,\R}$,
and the sets $S^{(\omega)}_\nu$ are the ascending sets
$$
S^{(\omega)}_\nu=\{x\in \wraffgr | \lim_{t\to+\infty}\varphi^t_\xi(x)\in \affgr^{(\omega)}_{M,\nu,\R}\},
$$
and the sets $T^{(\omega)}_\nu$ are the descending sets
$$
T^{(\omega)}_\nu=\{x\in \wraffgr | \lim_{t\to-\infty}\varphi^t_\xi(x)\in \affgr^{(\omega)}_{M,\nu,\R}\}.
$$
\end{Lemma}

\begin{proof}
The flow $\varphi_\xi$ 
respects the factorization, so the statement reduces
to the local Lemmas ~\ref{lrflow} and \ref{lpflow}.
\end{proof}

We have the following global version of Proposition~\ref{prhyp}.

\begin{Proposition}\label{prglobdual}
For $\nu\in\rcowt$,
there is a canonical isomorphism
$$
\pi_!s_\nu^*\simeq \pi_*t_\nu^!:\catd_\strat(\wraffgr)\to\catd(\base).
$$
\end{Proposition}

\begin{proof}
Consider the contractions
$$p:S^{(\omega)}_\nu\to\affgr^{(\omega)}_{M,\nu,\R} 
\mbox{ and }
q:T^{(\omega)}_\nu\to\affgr^{(\omega)}_{M,\nu,\R}$$
provided by the previous lemma.
For any sheaf $\sh\in\catd_\strat(\wraffgr)$, the supports of $p_! s_\nu^*\sh$
and $q_*t_\nu^!\sh$ are compact. 
Thus by the standard identity for the composition of maps,
it suffices to construct a canonical
isomorphism $p_! s_\nu^*\simeq q_*t_\nu^!$.
For $\xi\in \fs_\R$ in the interior of the dominant Weyl chamber for $\CW_{G_\R}$,
the map $\varphi_\xi:\wraffgr\to\wraffgr$ defined by $\varphi_\xi(x)=\exp(\xi)\cdot x$
is weakly hyperbolic in the sense of \cite[Section 1.3]{GM93}. Therefore by the previous lemma, the
asserted isomorphism may be deduced from \cite[Proposition 9.2]{GM93}.
\end{proof}


\begin{subsection}{Relation to specialization I}\label{sscharspec1}
Recall that we have the tensor functor
\begin{eqnarray*}
& \Sigma:\vect_\cowt \to\vect_{\rcowt\times\Z} & \\
&\displaystyle \Sigma (V)^{(\nu,z)} = \sum_{\Sigma(\lambda)=(\nu,z)} V^\lambda &
\end{eqnarray*}
where
$\Sigma:\cowt\to\rcowt\times\Z$ is the projection 
$\Sigma(\lambda)=(\theta(\lambda)+\lambda, \langle 2\ch\rho_M,\lambda\rangle)$.


\begin{Theorem}\label{tmainisom}
There is a canonical isomorphism 
$$
\schar\circ\pR\simeq \Sigma \circ\tchar:\epcat\to \vect_{\rcowt\times\Z}.
$$
\end{Theorem}

\begin{proof}
By the following lemma, it suffices to prove a similar
statement for 
the local specialization.

\begin{Lemma}\label{lpervss}
There is a canonical isomorphism
$$
\schar\circ {\pH}^*\simeq \schar:\catd_{G_\R(\cor)}(\raffgr^+)\to\vect_{\rcowt\times\Z}.
$$
\end{Lemma}

\begin{proof}
By Corollary~\ref{crwtexact},
the real weight functors are $t$-exact, so 
the real character functor is $t$-exact.
\end{proof}

By the lemma, it suffices to construct an isomorphism
$$
 \schar \circ \on R\simeq \Sigma \circ\tchar:\epcat\to \vect_{\rcowt\times\Z},
$$
or in other words, 
to construct isomorphisms
$$
F^{\nu,z}_\R\circ \on R\simeq \sum_{\Sigma(\lambda)=(\nu,z)}F^\lambda:\epcat\to\vect, \mbox{ for }
\nu\in\rcowt, z\in\Z.
$$

Recall the specialization diagram
$$
\begin{array}{ccccccccc}
\oneaffgr|X_+  & \simeq & 
\sraffgr|X_+ & \stackrel{j}{\to} 
& \sraffgr|X & \stackrel{i}{\gets} & \sraffgr|X_\R 
& \simeq & \oneraffgr \\
 & \searrow & \downarrow & & \downarrow & & \downarrow & \swarrow &\\
  && X_+ & \to & X & \gets & X_\R && 
\end{array}
$$ 

Define a functor
\begin{eqnarray*}
& \ls^\nu:\epcat\to\catd(X_+\cup X_\R) & \\
& \ls^\nu(\psh)=\pi_! s_\nu^*j_*(\rho(\psh)|{X_+}) &
\end{eqnarray*}
where 
$s_\nu:S_{\nu}^{(\sigma)}\to \sraffgr$ is the inclusion, and
$\pi:\sraffgr\to X$ is the projection. 

\begin{Lemma}
The sheaf $\ls^\nu(\psh)$
is constant.
\end{Lemma}

\begin{proof} 
First, observe that the sheaf $\ls^\nu(\psh)|X_+$ 
is the local system associated to the $\Aut(\CO)$-torsor
$\encurve_+\to X_+$ and the discrete $\Aut(\CO)$ representation $H^k_c(\mv,\psh)$.
Since $\Aut(\co)$ is connected, the sheaf $\ls^\nu(\psh)|X_+$ is constant.

Next, note the following chain of isomorphisms
$$
\ls^\nu(\psh)=
\pi_! s_\nu^*j_*(\psh_{X_+})\simeq \pi_* t_\nu^! j_*(\psh_{X_+})\simeq
j_*\pi_* t_\nu^!(\psh_{X_+})\simeq 
j_* \pi_! s_\nu^*(\psh_{X_+})\simeq j_*(\ls^\nu(\psh)|X_+).
$$
Here we used Proposition~\ref{prglobdual},
the base change isomorphism,
and the standard identity for the composition of maps.

Since each $x\in X_\R$ has a neighborhood that intersects $X_+$ in a contractible set,
the pushforward $ j_*(\ls^\nu(\psh)|X_+)$ is constant.
\end{proof}

\begin{Lemma}
We have canonical identifications of cohomology stalks
$$
H^{\langle\ch\rho,\nu\rangle+z}(\ls^\nu(\psh))_x \simeq F^{\nu,z}_\R(\on R(\psh)), \mbox{ for } x\in X_\R,  
$$
$$
H^{\langle\ch\rho,\nu\rangle+z}(\ls^\nu(\psh))_x \simeq \sum_{\Sigma(\mu)=(\nu,z)}
 F^\mu(\psh), \mbox{ for } x\in X_+.
$$ 

\end{Lemma}

\begin{proof}
For $x\in X$, we have
$$
H^*(\ls^\nu(\psh))_x
\simeq H^*_c(S^{(\sigma)}_{\nu}|x, (j_*\rho(\psh)|{X_+})_x).
$$

By definition, we have identifications
$$
S^{(\sigma)}_{\nu}|x \simeq {}_x\rmv, \mbox{ for } x\in X_\R,
$$
$$
S^{(\sigma)}_{\nu}|x \simeq 
\bigsqcup_{\sigma(\mu)=\nu} {}_x S^{\mu}, \mbox{ for } x\in X_+.$$
By choosing formal coordinates, we obtain isomorphisms
$$
H^*_c({}_x\rmv, j_*(\psh_{X_+})_x) \simeq H^*_c(\rmv, \on R(\psh)),   \mbox{ for } x\in X_\R,
$$
$$
H^*_c({}_xS^\mu, (\psh_{X_+})_x) \simeq H^*_c(S^\mu, \psh),   \mbox{ for } x\in X_+.
$$
Since $\Aut^0(\cor)$, respectively $\Aut(\CO)$, is connected, the 
first, respectively second, isomorphism is independent
of the choice of formal coordinate.

This proves the first assertion of the lemma. For the second,
thanks to Corollary~\ref{cevenpar}, it only remains to check that
the $\Z$-gradings agree. 
We must check that 
if $\nu=\theta(\mu)+\mu$, then
we have
$$
\langle 2\ch\rho,\mu\rangle=\langle\ch\rho,\nu\rangle+\langle2\ch\rho_M,\mu\rangle.
$$
Using the identities $\ch\theta(\ch\rho-\ch\rho_M)=\ch\rho-\ch\rho_M$, 
and $\langle 2\ch\rho_M, \nu\rangle=0$, we have
$$
\langle \ch\rho,\nu\rangle=
\langle\ch\rho-\ch\rho_M,\theta(\mu)-\mu\rangle=
\langle 2(\ch\rho-\ch\rho_M),\mu\rangle.
$$
\end{proof}

The identification of the stalks of the constant cohomology
sheaf $H^z(\ls^\nu(\psh))$, for $\nu\in\rcowt$, $z\in\Z$, provides the sought-after isomorphism.
\end{proof}

We have the following corollaries.
Recall that $\bigimlat\subset\rcowt\times\Z$ is the image
of the projection $\Sigma:\cowt\to\rcowt\times\Z$.

\begin{Corollary}\label{ccharim}
The real character functor 
descends to a functor
$$\schar:\catq_\pR\to\vect_{\bigimlat}.$$
\end{Corollary}

\begin{proof}
By definition, 
an object $\CQ\in\catq_\pR$ is isomorphic to a subquotient of an object of the form $\pR(\psh)$,
for $\psh\in\epcat$. Since each real weight functor is exact,
the real character functor is exact, and so $\schar(\CQ)$ is a subquotient
of $\schar(\pR(\psh))$. By the theorem, we have 
$\schar(\pR(\psh))\simeq\Sigma(\tchar(\psh))$, and so $\schar(\CQ)$ is a subquotient
of $\Sigma(\tchar(\psh))$.
\end{proof}


Recall that $\hc$ denotes the hypercohomology functor,
and $\ff$ the forgetful functor
which forgets the grading of a vector space.

\begin{Corollary}\label{crevenpar}
There are canonical isomorphisms
\begin{eqnarray*}
& \hc\simeq\ff\circ\schar:\catq_\pR\to \vect & \\
& \hc\simeq \hc\circ \pR:\epcat\to \vect. &
\end{eqnarray*}
\end{Corollary}

\begin{proof}
By definition, 
an object $\CQ\in\catq_\pR$ is isomorphic to a subquotient of an object of the form $\pR(\psh)$,
for $\psh\in\epcat$. 
We may assume without loss of generality that $\psh$ is supported on a single
component of $\affgr$, and that $\CQ$ is a subquotient
of $\pH^k(\on R(\psh))$, for some $k$.

By Theorem~\ref{twt} and Proposition~\ref{plgstrat},
for $\psh\in\epcat$ supported on a single component of $\affgr$,
$\tchar(\psh)$ is non-zero in only even or odd
degrees. 
Thus by the theorem, $\schar(\pH^k(\on R(\psh)))$ is non-zero
in only even or odd degrees, and since $\schar$ is exact, $\schar(\CQ)$ 
is non-zero in only even or odd degrees.
Therefore the 
compactly
supported cohomology along the orbits $\rmv$ provides a filtration of 
the hypercohomology functor. Similarly, the local cohomology along the
orbits $T^\nu_\R$ also provides a filtration. The isomorphism of Proposition~\ref{prhyp}
implies that the two filtrations split each other.

This proves the first assertion, 
and now the second is obtained from the theorem by applying the forgetful functor
and using the first assertion and Corollary~\ref{cevenpar}.
\end{proof}

%

%


%

\begin{Corollary}
The hypercohomology $\hc:\catq_\pR\to\vect$ is exact and faithful.
\end{Corollary}

\begin{proof}
The real weight functors are exact, and so the real character functor is exact.
By the previous corollary, the hypercohomology functor is exact as well.
Since the hypercohomology functor is exact, 
to see that it is faithful,
it suffices to check that
it does not vanish on any $\CQ\in\catq_\pR$.
But for any $\CQ\in\catq_\pR$, some weight functor will not vanish on $\CQ$.
For example, we could take the weight functor $F^{\nu}_\R$ associated 
to the dominant real coweight $\nu$ in an open stratum in the support of $\CQ$. 
\end{proof}
\end{subsection}

\begin{Corollary}\label{cpiffqs}
The specialization $\on R:\epcat\to\catd_{G_\R(\cor)}(\raffgr^+)$ 
is perverse if and only
if $G_\R$ is quasi-split.
\end{Corollary}

\begin{proof}
By the theorem, the $\Z$-grading on the perverse specialization is trivial
if and only $G_\R$ is quasi-split.
\end{proof}


\begin{subsection}{Monoidal structure for real character functor}\label{sscharmon}
In this section,
we 
equip the real character functor $\schar:\catq_\pR\to\vect_\bigimlat$, and the hypercohomology
$\hc:\catq_\pR\to\vect$ with monoidal structures.

\begin{Theorem}\label{trhypmon}
There is a canonical isomorphism
$$
\hc(\cdot\odot\cdot)\simeq\hc(\cdot)\otimes\hc(\cdot)
$$
for the hypercohomology functor $\hc:\catq_\pR\to\vect_{}$.
\end{Theorem} 

\begin{proof}
Recall the real global convolution diagram
$$
\oneraffgr\times\oneraffgr \stackrel{p}{\gets}
\widetilde{\oneraffgr\times\oneraffgr} \stackrel{q}{\to}
\oneraffgr\tilde \times\oneraffgr \stackrel{m}{\to}
\tworaffgr\stackrel{d}{\gets}\oneraffgr.
$$

Define a functor
\begin{eqnarray*}
& \ls_\R:\catq_\pR\times \catq_\pR\to\catd(X_\R^2) & \\ 
& \ls_\R(\CQ_1,\CQ_2)=\pi_! (\rho_\R(\CQ_1)\tilde\boxtimes \rho_\R(\CQ_2)) &
\end{eqnarray*}
where $\pi:\oneraffgr\tilde\times\oneraffgr\to X^2_\R$ is the projection.

\begin{Proposition}\label{lhyperconst}
The sheaf $\ls_\R(\CQ_1,\CQ_2)$
is constant with cohomology stalk 
$$
H^*(\ls_\R(\CQ_1,\CQ_2))_x\simeq
\hc(\raffgr,\CQ_1)\boxtimes\hc(\raffgr,\CQ_2), \mbox{ for } x\in X^2_\R.
$$
\end{Proposition}

\begin{proof}

Recall that $\oneraffgr$ classifies data
$
(x,\tors,\nu),
$ 
where
$x\in X_\R,$ $\tors$ is a $G_\R$-torsor on $X_\R$, and
$\nu$ is a trivialization of $\tors$ over $X_\R\setminus x$.
By Lemma~\ref{lcomp}, 
a sheaf $\CQ_1\in\catq_\pR$ is supported on the components $\raffgr^0\subset\raffgr$
that are in the kernel of the map $\partial:\pi_0(\raffgr)\to\pi_0(G_\R)$.
Thus the sheaf $\rho_\R(\CQ_1)$ is supported on the components 
of $\oneraffgr$ for which the trivialization $\nu$
extends to a trivialization of the torsor $\tors/G_\R^0\to X_\R$
for the component group $G_\R/G_\R^0$, or
in other words, for which
the torsor $\tors$ is induced from a torsor for the connected group $G_\R^0$.

Recall that ${\oneraffgr\tilde\times\oneraffgr}$ 
classifies data
$
(x_1, x_2,
\tors_1, \tors,
\nu_1, 
\eta)
$
where $x_1,x_2\in X_\R$,
$\tors_1,\tors$ are $G_\R$-torsors on $X_\R$,
$\nu_1$ is a trivialization of $\tors_1$ over $X_\R\setminus x_1$,
and
$\eta$ is an isomorphism from $\tors_1$ to $\tors$ over $X_\R\setminus x_2$.
We have also seen that it is the twisted product 
constructed from 
the action of $G_\R(\cor)\rtimes\Aut^0(\cor)$ on $\raffgr$,
and the
$G_\R(\cor)\rtimes\Aut^0(\cor)$-torsor 
$$
\widetilde{\oneraffgr\times \hat X^0_\R} {\to} \oneraffgr
$$
that classifies data
$
(x_1, x_2, \tors_1, \nu_1, \mu_1, \phi),
$
where $x_1, x_2\in X_\R$,
$\tors_1$ is a $G$-torsor on $X_\R$,
$\nu_1$ is a trivialization of $\tors_1$ 
   over $X_\R\setminus x_1$,
$\mu_1$ is a trivialization of $\tors_1$ over the formal neighborhood of $x_2$,
and
$\phi$ is an orientation-preserving isomorphism from the abstract formal disk 
to the formal neighborhood of $x_2$.

Consider the projection
$$
r:\oneraffgr\tilde\times\oneraffgr\to\oneraffgr
$$
defined by $(x_1,x_2,\tors_1,\tors,\nu_1,\eta)\mapsto(x_1,\tors_1,\nu_1)$,
with fiber the $G_\R(\cor)\rtimes\Aut^0(\cor)$-twist of $\raffgr$.
We claim that there is a canonical isomorphism
$$
r_!(\rho_\R(\CQ_1)\tilde\boxtimes \rho_\R(\CQ_2))
\simeq
\rho_\R(\CQ_1)\boxtimes\hc(\raffgr,\CQ_2).
$$
Assuming this isomorphism for the moment, 
we shall then have 
$$
\ls_\R(\CQ_1,\CQ_2)\simeq
\pi_!r_!(\rho_\R(\CQ_1)\tilde\boxtimes \rho_\R(\CQ_2))\simeq
\pi_!(\rho_\R(\CQ_1)\boxtimes\hc(\raffgr,\CQ_2))\simeq
\hc(\raffgr,\CQ_1)\boxtimes\hc(\raffgr,\CQ_2).
$$
Here the last isomorphism $\pi_!(\rho_\R(\CQ_1))\simeq\hc(\raffgr,\CQ_1)$
follows from the fact that $\Aut^0(\cor)$ is connected.
Thus to prove the lemma, we are left to establish
the asserted isomorphism.

Observe that the restriction of the $G_\R(\cor)\rtimes\Aut^0(\cor)$-torsor
$$
\widetilde{\oneraffgr\times \hat X^0_\R} {\to} \oneraffgr
$$
to the components of $\oneraffgr$ 
that classify data $(x,\tors,\nu)$ where $\tors$ is induced
from a $G_\R^0$-torsor, is itself induced from a torsor for the 
connected group $G_\R(\cor)^0\rtimes\Aut^0(\cor)$.
The sheaf
$
r_!(\rho_\R(\CQ_1)\tilde\boxtimes \rho_\R(\CQ_2))
$
is the product of $\rho_\R(\CQ_1)$ with the 
$G_\R(\cor)^0\rtimes\Aut^0(\cor)$-twist of $\hc(\raffgr,\CQ_2)$,
and since $G_\R(\cor)^0\rtimes\Aut^0(\cor)$ is connected, the twist must be trivial.
\end{proof}

The identification of the stalks of the constant cohomology sheaf $H^*(\ls_\R(\CQ_1,\CQ_2))$
provides the sought-after isomorphism.
\end{proof}



\begin{Theorem}\label{tcharmon}
There is a canonical isomorphism
$$
\schar(\cdot\odot\cdot)\simeq\schar(\cdot)\otimes\schar(\cdot)
$$
for the real character functor $\schar:\catq_\pR\to\vect_{\bigimlat}$.
\end{Theorem} 

\begin{proof}

For $\nu\in\rcowt$,
define functors
\begin{eqnarray*}
& \ls^\nu_\R:\catq_\pR\times \catq_\pR\to\catd(X_\R^2) & \\
& \ls^\nu_\R(\CQ_1,\CQ_2)=\pi_! s_\nu^*m_!(\rho_\R(\CQ_1)\tilde\boxtimes \rho_\R(\CQ_2)) &
\end{eqnarray*}
where $m:\oneraffgr\tilde \times\oneraffgr{\to}\tworaffgr$ is the multiplication, 
$s_\nu:S^{(2)}_{\nu}\to\affgr^{(2)}_\R$ is the inclusion, and
$\pi:\tworaffgr\to X_\R^2$ is the projection.

\begin{Lemma}\label{lrstalk}
There are canonical identifications of cohomology stalks
$$
H^*(\ls_\R^\nu(\CQ_1,\CQ_2))_x \simeq F^\nu(\CQ_1\odot\CQ_2), \mbox{ for } x\in \Delta_\R,
$$
$$
H^*(\ls_\R^\nu(\CQ_1,\CQ_2))_x \simeq \sum_{\lambda+\mu=\nu} F^\lambda(\CQ_1)\otimes F^\mu(\CQ_2),
\mbox{ for }  x\in X_\R\setminus \Delta_\R.
$$
\end{Lemma}

\begin{proof}
For $x\in X_\R^2$, we have
$$
H^*(\ls_\R^\nu(\CQ_1,\CQ_2))_x= 
H^*_c(S^{(2)}_{\nu}|x,m_!(\rho_\R(\CQ_1)\tilde\boxtimes\rho_\R(\CQ_2))_x ).
$$

By definition, we have identifications
$$
S^{(2)}_{\nu}|x \simeq {}_x\rmv,\mbox{ for }
x\in\Delta_\R,
$$
$$
S^{(2)}_{\nu}|x \simeq \bigsqcup_{\lambda+\mu=\nu} 
{}_{x_1}S^{\lambda}_\R\times {}_{x_2}S^{\mu}_\R,
\mbox{ for } (x_1,x_2)\in X^2_\R\setminus \Delta_\R.
$$ 
By choosing formal coordinates, we obtain 
isomorphisms
$$
H^*_c({}_x\rmv,m_!(\rho_\R(\CQ_1)\tilde\boxtimes\rho_\R(\CQ_2))_x )\simeq
H^*_c(\rmv,\CQ_1\odot\CQ_2), \mbox{ for } x\in\Delta_\R,
$$
$$
H^*_c({}_{x_1}S_\R^\lambda\times{}_{x_2}S_\R^\mu,
m_!(\rho_\R(\CQ_1)\tilde\boxtimes\rho_\R(\CQ_2))_{(x_1,x_2)} )\simeq
H^*_c(S^\lambda_\R,\CQ_1)\otimes H^*_c(S^\mu_\R,\CQ_2), 
\mbox{ for } (x_1,x_2) \in X^2_\R\setminus\Delta_\R.
$$
Since $\Aut^0(\cor)$ is connected, the 
isomorphisms are indepedendent
of the choice of formal coordinates. 
\end{proof}

\begin{Lemma}\label{lwtconst}
The sheaf $\ls^\nu_\R(\CQ_1,\CQ_2)$
is constant.
\end{Lemma}

\begin{proof}
For $i=1,2$,
an object $\CQ_i\in\catq_\pR$ is isomorphic to a subquotient of an object of the form $\pR(\psh_i)$,
for $\psh_i\in\epcat$. 
We may assume without loss of generality that $\psh_i$ is supported on a single
component of $\affgr$, so that the hypercohomology of $\psh_i$ is possibly non-zero
only in degrees congruent to some $d_i\mod 2$. We may also assume 
that $\CQ_i$ is a subquotient
of $\pH^{k_i}(\on R(\psh_i))$, for some $k_i$, so that by the proof of Corollary~\ref{crevenpar}, 
$F_\R^{\nu}(\CQ_i)$ is possibly non-zero only in degrees congruent to $d_i+k_i \mod 2$.
The product $\CQ_1\odot\CQ_2$ is isomorphic to a subquotient of the sheaf 
$\pH^{k_1}(\on R(\psh_1))\odot\pH^{k_2}(\on R(\psh_2))$, so that
by the proof of Corollary~\ref{crevenpar}, 
$F^\nu_\R(\CQ_1\odot\CQ_2)$ is possibly non-zero
only in degrees congruent to $d_1+k_1+d_2+k_2\mod 2$.
Therefore by the previous lemma,
the sheaves $\ls_\R^\nu(\CQ_1,\CQ_2)$
provide a filtration of the sheaf $\ls_\R(\CQ_1,\CQ_2)$,
and similarly the sheaves constructed for the opposite minimal
parabolic subgroup also provide a filtration. The isomorphism
of Proposition~\ref{prglobdual} implies that the two filtrations
split each other.
%
The lemma now follows from the fact that a direct summand of a 
constant sheaf must itself be constant.
\end{proof}

The identification of the stalks of the 
constant cohomology sheaf $H^*(\ls_\R^\nu(\CQ_1,\CQ_2))$, for $\nu\in\rcowt$, 
provides the sought-after isomorphism.
\end{proof}

\begin{Remark}
Note that by construction, the monoidal structure for the hypercohomology functor
$\hc:\catq_\pR\to\vect$ 
is obtained from that of the real chracter functor $\schar:\catq_\pR\to\vect_\bigimlat$
by applying the forgetful functor.
\end{Remark}
%


%

\end{subsection}


\begin{subsection}{Relation to specialization II}\label{sscharspec2}
In Theorem~\ref{tmainisom}, we constructed a canonical isomorphism
$$
\schar \circ\pR\simeq \Sigma \circ\tchar:\epcat\to \vect_{\bigimlat}
$$
intertwining the perverse specialization and the character functors.
In Section~\ref{smon}, we constructed a monoidal structure for
the perverse specialization $\pR$, and in the previous section, 
we constructed a monoidal structure for the real character functor $\schar$.
Therefore the functor on the left hand side has a monoidal structure.
One may give the
character functor $\tchar$ a monoidal structure 
in the same way as we did 
the real character functor $\schar$
in the previous section. (See~\cite[Proposition 6.4]{MV04}.)
Therefore the functor on the right hand side has a monoidal structure as well.
The following proposition confirms that the monoidal structures agree.

\begin{Proposition}\label{tmainisommon}
We have a commutative square
$$
\begin{array}{ccc}
\schar(\pR(\cdot\odot\cdot)) & \risom &  
\Sigma(\tchar(\cdot\odot\cdot)) \\

\wr\downarrow & & \wr\downarrow \\

\schar(\pR(\cdot))\otimes\schar(\pR(\cdot)) & \risom &  
\Sigma(\tchar(\cdot))\otimes \Sigma(\tchar(\cdot))
\end{array}
$$
\end{Proposition}

The proof is an elaboration of the techniques of previous sections,
and we leave it to the reader. 
Applying the forgetful functor, and using Corollaries~\ref{cevenpar} and~\ref{crevenpar},
we obtain the following corollary confirming that under the isomorphism
$$
\hc\circ\pR\simeq\hc:\epcat\to\vect,
$$
the monoidal structures on each side agree.

\begin{Corollary}\label{chcpr}
We have a commutative square
$$
\begin{array}{ccc}
\hc(\pR(\cdot\odot\cdot))  & \risom &  
\hc(\cdot\odot\cdot) \\

\wr\downarrow & & \wr\downarrow \\

\hc(\pR(\cdot))\otimes \hc(\pR(\cdot)) & \risom &  
\hc(\cdot)\otimes\hc(\cdot) 
\end{array}
$$

\end{Corollary}

\end{subsection}

\end{section}



\begin{section}{The image category II}\label{sim2}

\begin{subsection}{Tannakian dictionary}

We review the dictionary between monoidal categories
and bi-algebras developed in \cite{Saav72} and \cite{DM82}. 
The results are valid over any field $\BK$ which we will take to be $\C$.

Let
 $\sC$ be a
$\BK$-linear abelian category equipped with a $\BK$-linear
exact faithful functor $\omega:\sC\to\Vect$.
For an object $X\in\Ob(\sC)$, let $\langle X\rangle$ be the strict full subcategory
of $\sC$ whose objects are isomorphic to subquotients of the objects $X^{\oplus n}$,
for some $n$. Define the $\BK$-algebra 
$$
A(X)=\on{End}(\omega|\langle X\rangle)
$$
to be the algebra of 
endomorphisms of the functor $\omega$ restricted to the subcategory $\langle X\rangle$.
Define the $\BK$-coalgebra 
$$
B(X)=\Hom(A(X),\BK)
$$
to be the coalgebra dual to the algebra $A(X)$.
Let $B(\sC)$ be the inverse limit of the coalgebras $B(X)$ over all objects $X\in\Ob(\sC)$.

\begin{Proposition}[\cite{DM82}, Proposition 2.14]
There is an equivalence of categories $\sC\to\on{Comod}_{B(\sC)}$
under which $\omega$ corresponds to the forgetful functor.
\end{Proposition}

A homomorphism $f:B\to B'$ defines the functor
$$
\phi^f:\Comod_{B}\to\Comod_{B'}
$$
that takes a $B$-comodule $X$ with coaction $a_X:X\to X\otimes B$
to the $B'$-comodule $X$ with coaction
$$
a'_X:X\stackrel{a_X}{\to} X\otimes B\stackrel{1\otimes f}{\to} X\otimes B'.
$$

\begin{Proposition}[\cite{DM82}, p 135]
The map $f\mapsto \phi^f$ defines a bijection from 
the set of homomorphisms $f:B\to B'$
to the set of functors 
$\phi:\Comod_{B}\to\Comod_{B'}$ which
cover the identity on the underlying vector spaces.
\end{Proposition}

For a $\BK$-coalgebra $B$, a homomorphism $u:B\otimes B\to B$
defines the functor 
$$
\phi^u:\Comod_B\times \Comod_B\to\Comod_B
$$
that takes a pair of $B$-comodules $X,Y$, with coactions $a_X:X\to X\otimes B$,
$a_Y:Y\to Y\otimes B$,
to the $B$-comodule $X\otimes Y$
with coaction
$$
a_{X\otimes Y}:X\otimes Y\stackrel{a_X\otimes a_Y}{\to} X\otimes B\otimes Y\otimes B
\stackrel{1\otimes u}{\to} X\otimes Y\otimes B.
$$

\begin{Proposition}[\cite{DM82}, Proposition 2.16]
The map 
$
u\mapsto \phi^u
$
defines a bijection from 
the set of homomorphisms $u:B \otimes B\to B$
to
the set of functors 
$
\phi:\Comod_{B}\times \Comod_{B}\to \Comod_{B}
$ which
cover the ordinary tensor product on the underlying vector spaces.

The natural associativity and commutativity constraints
on $\Vect$ induce similar constraints on $\Comod_{B}$ equipped with 
the product $\phi^u$
if and only if $u$ is associative and commutative. There is an identity
object for $\phi^u$ in $\Comod_{B}$ with underlying vector space $\BK$
if and only if there is an identity element for $u$ in $B$.
\end{Proposition}

When there is an associative and commutative homomorphism
$u:B\otimes B\to B$ and an identity element for $u$ in $B$, then 
$\Spec(B)$ is an affine monoid-scheme and there is an equivalence
$$
\Rep(\Spec(B))\simeq\Comod_B.
$$
To obtain $\Spec(B)$ directly from the category $\Comod_B$,
one checks that it is the monoid of tensor endomorphisms of the forgetful functor.

\begin{Proposition}[\cite{Ulb90}]
Suppose the homorphism $u:B\otimes B\to\ B$ is associative
and there is an identity element for $u$ in $B$.
Then the tensor category $\Comod_B$ is rigid if and only if there is an antipode $S:B\to B$
making $B$ a Hopf algebra.
\end{Proposition}

When the multiplication $u:B\otimes B\to B$ is commutative
and there is an antipode $S:B\to B$ making $B$ a commutative Hopf algebra,
then $\Spec(B)$ is an affine group-scheme. One checks that in this
case every tensor endomorphism of the forgetful functor
is in fact an automorphism.

\end{subsection}


\begin{subsection}{Coherence constraints}

The category $\epcat$ equipped with the convolution product is a rigid tensor
category~\cite[Theorem 7.3]{MV04}, 
and the hypercohomology $\hc:\epcat\to\vect$ is an exact faithful tensor
functor~\cite[Corollary 3.7, Proposition 6.3]{MV04}.

We have seen that the category $\catq_\pR$ is equipped with a convolution product,
and that the perverse specialization $\pR:\epcat\to\catq_\pR$ is equipped with
a monoidal structure such that there is a canonical isomorphism
of functors
$\hc\simeq\hc\circ\pR:\epcat\to\vect$ respecting the monoidal structures.

\begin{Proposition}
There exist unique associativity and commutativity constraints for 
the category $\catq_\pR$ equipped with the convolution product
such that 
the perverse specialization
$\pR:\epcat\to\catq_\pR$ and hypercohomology $\hc:\catq_\pR\to\vect$
equipped with their monoidal structures
respect the constraints.
\end{Proposition}

\begin{proof}
First, since the canonical isomorphism $\hc\simeq\hc\circ\pR$ 
respects monoidal structures,
the composite functor $\hc\circ\pR$ 
is an exact faithful
tensor functor. We may apply the Tannakian dictionary to the category $\epcat$
and the functor $\hc\circ\pR$ 
to obtain a commutative Hopf algebra $B({\epcat})$
with identity such that $\epcat$ is equivalent to the category of $B({\epcat})$-comodules.

Next, we may apply the Tannakian dictionary to the category $\catq_\pR$
and the functor $\hc$ to obtain a coalgebra $B({\catq_\pR})$ with a multiplication
such that $\catq_\pR$ is equivalent to the category of $B({\catq_\pR})$-comodules.
By the Tannakian dictionary,
the perverse specialization $\pR$ provides a coalgebra
homomorphism 
$$
B({\epcat})\to B({\catq_\pR})
$$ 
which is surjective
since every object of $\catq_\pR$ is isomorphic to a subquotient of
an object of the form $\pR(\psh)$, for some object $\psh$ in the category $\epcat$.

\begin{Lemma}
For coalgebras $B, B'$, 
given a homomorphism $f:B\to B'$ inducing the functor
$\phi^f:\Comod_B\to\Comod_{B'}$, and
given homomorphisms
$u:B\otimes B\to B,u':B'\otimes B'\to B'$
inducing the functors $\phi^u:\Comod_B \otimes\Comod_B\to\Comod_B,
\phi^{u'}:\Comod_{B'} \otimes\Comod_{B'}\to\Comod_{B'},$
there is an identity of homomorphisms
$$
f\circ u= u'\circ(f\otimes f) 
$$
if and only if the identity isomorphism in $\Vect$ induces
an isomorphism of functors
$$
\phi^f\circ \phi^u\simeq \phi^{u'}\circ (\phi^f\times\phi^f).
$$
\end{Lemma}

\begin{proof}
Follows directly from the definitions.
\end{proof}

Applying the lemma to the perverse specialization $\pR$,
we see that the coalgebra homomorphism $B({\epcat})\to B({\catq_\pR})$
respects multiplication. Since it is surjective, we conclude
that the multiplication of $B({\catq_\pR})$ is associative
and commutative. By the Tannakian dictionary, we obtain unique
associativity and commutativity constraints on the category
$\catq_\pR$ equipped with the convolution product
such that the perverse specialization $\pR$ and hypercohomology $\hc$
respect the constraints.
\end{proof}

\end{subsection}


\begin{subsection}{Rigidity}
The following is a result of Waterhouse. 

\begin{Proposition}[\cite{Nichols78}, Theorem 0]
Let $B$ be a finitely-generated Hopf algebra over a field $\BK$ with comultiplication
$\Delta:B\to B\otimes B$, counit $\varepsilon:B\to\C$, and antipode $S:B\to B$.
Let $I\subset B$ be an ideal such that $\Delta(I)\subset I\otimes B+ B\otimes I$.
Then $\varepsilon(I)=\{0\}$, and $S(I)=I$.
\end{Proposition}

When $B$ is commutative, the proposition
implies that a closed sub-semigroup $\Spec(B/I)$ of an affine
algebraic group $\Spec(B)$ automatically 
contains the identity element and is closed under taking inverses,
and therefore is itself a group.

Applying the proposition to the kernel of the surjective homomorphism $B(\epcat)\to B(\catq_\pR)$,
we conclude that the unit and antipode of $B(\epcat)$ descend to the quotient 
$B(\catq_\pR)$ making it a Hopf algebra. By the Tannakian dictionary, this confirms
that $\catq_\pR$ has an identity object and is rigid.

\end{subsection}


\begin{subsection}{Other categories}\label{sothercats}
%

By the Tannakian dictionary, 
from the commutative diagram of Section~\ref{sim1},
we obtain a commutative diagram of 
coalgebras with multiplication
$$
\begin{array}{ccccc}
& & B(\catq_\pR) & & \\
& & \downarrow & \stackrel{}{\searrow}  & \\
B(\catq) & \stackrel{}{\to} & B(\catq_\Z) & \stackrel{}{\to} & B(\catq). \\
\end{array}
$$

By definition,
every object of $\catq$ is isomorphic to a subquotient of
an object of the form $\ff(\CQ)$, for some object $\CQ$ in the category $\catq_\pR$,
where $\ff:\catq_\pR\to\catq$ is the forgetful functor.
Therefore the diagonal arrow 
is surjective, and
we conclude that $B(\catq)$ is a commutative Hopf algebra.

Since the category $\catq_\Z$ is the product category $\catq\otimes\vect_\Z$ with the product
convolution product, the coalgebra $B(\catq_\Z)$ is the product coalgebra
$B(\catq)\otimes B(\vect_\Z)$ with the product multiplication.
Noting that $B(\Z)$ is nothing but the group algebra $\C[\Z]$,
we conclude that $B(\catq_\Z)$ is a commutative Hopf algebra.

The vertical arrow and first horizontal arrow are the obvious inclusions,
and the second horizontal arrow is the obvious projection.

\end{subsection}
\end{section}


\begin{section}{The associated subgroup}\label{ssub}

In this section, we identify the Tannakian group associated to the
category $\catq$ with the fiber functor $\hc:\catq\to\vect$.


\begin{subsection}{Root data}

A {\em based root datum} $\Psi$ consists of a quadruple
$$
\Psi=(\ch\Lambda, \ch\Delta,\Lambda,\Delta)
$$
where $\ch\Lambda$ and $\Lambda$ are lattices in duality
with respect to a fixed perfect pairing
$$
\langle\cdot,\cdot\rangle:\ch\Lambda\times\Lambda\to\BZ,
$$
and $\ch\Delta\subset\ch\Lambda$ and $\Delta\subset\Lambda$ are 
subsets with a fixed bijection 
$$
\ch\Delta\simeq\Delta,
$$
such that the quadruple satisfies certain requirements~\cite[Sections 7.4.1, 16.2.1]{Spr98}.
The requirements are symmetric, and 
the dual based root datum $\ch\Psi$ is the quadruple
$$
\ch\Psi=(\Lambda, \Delta,\ch\Lambda,\ch\Delta),
$$
where $\Lambda$ and $\ch\Lambda$ are in duality
with respect to the same perfect pairing, and $\Delta$ and $\ch\Delta$
are related by the same bijection.

The based root datum $\Psi(G,B,T)$ of a connected reductive complex algebraic
group $G$ with respect to a Borel subgroup $B\subset G$, and maximal torus $T\subset B$,
is the quadruple
$$
\Psi(G, B, T)=(\wt, \ch\Delta_{B,T},\cowt,\Delta_{B,T}).
$$
For any other choice of Borel subgroup $B'\subset G$, and maximal torus $T'\subset B$,
there is a canonical isomorphism
$$
\Psi(G, B, T)\simeq\Psi(G,B',T').
$$
The based root datum
$
\Psi(G) 
$
is defined to be the projective limit of the based root data $\Psi(G,B,T)$
with respect to the canonical isomorphisms
over the set of choices of Borel subgroups $B\subset G$, and maximal tori $T\subset B$.

There is a short exact sequence
$$
1\to\Inn(G)\to\Aut(G)\stackrel{\psi}{\to}\Aut(\Psi(G))\to 1
$$
where $\Inn(G)$ is the group of inner automorphisms of $G$.
We may split the sequence as follows.
Choose a Borel subgroup $B\subset G$, a maximal torus $T\subset B$, and
for each simple root $\ch\alpha\in\ch\Delta_{B,T}$, 
a basis vector $X_{\ch\alpha}$
in the corresponding simple root space for $T$ acting on $\fb$.
Define the subgroup $\Aut(G,B,T,\{X_{\ch\alpha}\})\subset\Aut(G)$ 
to consist of the automorphisms that preserve $B$, $T$, and the set $\{X_{\ch\alpha}\}$.
Then the restriction of the homomorphism $\psi$ to 
$\Aut(G,B,T,\{X_{\ch\alpha}\})$
is an isomorphism. 
See~\cite[Sections 1 and 2]{Spr79} for more details.

\end{subsection}


\begin{subsection}{The dual group}\label{ssdualgroup}

A {\em dual group} for a connected reductive complex algebraic group $G$ 
is a connected reductive complex algebraic group $\ch G$
together with a fixed isomorphism 
$$
\Psi(\ch G)\simeq\ch\Psi(G).
$$
By the above short exact sequence,
any two dual groups are isomorphic, and the isomorphism
is canonical up to composition with inner automorphisms.
Note that there is a canonical isomorphism
$$
\Aut(\Psi(G))\simeq\Aut(\Psi(\ch G)),
$$
and for the torus $T$, the dual torus $\ch T$ is canonically a dual group.

Recall that 
the category $\epcat$ of $G(\CO)$-equivariant perverse sheaves on $\affgr$
is equivalent as a tensor category to the category $\Rep(\ch G)$ of finite-dimensional
representations of the dual group $\ch G$
so that the hypercohomology $\hc:\epcat\to\vect$ 
corresponds
to the forgetful functor~\cite[Theorem 7.3, Corollary 3.7, Proposition 6.3]{MV04}.
In other words, the group of tensor automorphisms $\taut(\hc)$ 
is a connected reductive complex algebraic group which we denote by $\ch G$,
and there is a canonical isomorphism of based root data
$$
\Psi(\ch G)\simeq\ch\Psi(G).
$$  

To see the root datum of $\ch G$, first recall that 
the character functor $\tchar:\epcat\to\vect_\cowt$ is a tensor functor~\cite[Proposition 6.4]{MV04},
and there is a canonical isomorphism 
$$
\hc\simeq\ff\circ\tchar:\epcat\to\vect
$$
where $\ff:\vect_\cowt\to\vect$ is the forgetful functor~\cite[Theorem 3.6]{MV04}.
Therefore we may identify the group $\ch G$ with the group of tensor automorphisms
of the composite functor $\ff\circ\tchar$,
and we have a canonical homomorphism
$$
\Aut^\otimes(\ff)\to \ch G. 
$$
This corresponds to 
an embedding of the dual torus 
$$
\ch T=\Aut^\otimes(\ff)
$$
as a maximal torus. (See~\cite[Theorem 7.3]{MV04} or Proposition~\ref{appmaximaltorus} of the appendix.) 

Next, recall that 
the simple objects in the category $\epcat$ are the intersection cohomology sheaves
$\ic^\lambda$ of the closures of the $G(\CO)$-orbits $\affgr^\lambda\subset\affgr$,
with coefficients in trivial one-dimensional local systems. (See~\cite[Lemma 7.1]{MV04}
or Proposition~\ref{appsemisimple} of the appendix.)
The weight functor $F^\lambda:\epcat\to\vect$ evaluated on $\ic^\lambda$ provides
a canonical line in the graded vector space $\tchar(\ic^\lambda)$.
By forgetting the grading, for each simple object $\ic^\lambda$,
we obtain a canonical line $\ff(F^\lambda(\ic^\lambda))$
in the vector space $\ff(\tchar(\ic^\lambda))$. 
The subgroup 
$$\ch B\subset\ch G$$ 
of those tensor automorphisms of $\ff\circ\tchar$
which preserve the collection of lines $\ff(F^\lambda(\ic^\lambda))$
is a Borel subgroup. By construction, it contains the maximal torus $\ch T$.

Finally, 
there is a canonical isomorphism of based root data
$$
\Psi(\ch G,\ch B,\ch T)\simeq\ch\Psi(G,B,T).
$$
See \cite[Theorem 7.3]{MV04} or the appendix for more details.

\end{subsection}


\begin{subsection}{Associated subgroups}

Via the Tannakian dictionary, we have seen that the category $\catq_\pR$ 
is equivalent as a tensor category to the category $\Rep(\ch H_\pR)$ of finite-dimensional
representations of a group-scheme $\ch H_\pR$
so that the hypercohomology $\hc:\catq_\pR\to\vect$ 
corresponds
to the forgetful functor.

We have also seen that the perverse specialization
$
\pR:\epcat\to\catq_\pR
$
is a tensor functor, and that there is a canonical isomorphism
of tensor functors
$$
\hc\simeq\hc\circ\pR:\epcat\to\vect.
$$
Therefore we may identify $\ch G$ with the group of tensor automorpshisms
of the composite functor $\hc\circ\pR$,
and we have a canonical homomorphism
$$
\ch H_\pR\to \ch G.
$$
Since every object of $\catq_\pR$ is isomorphic to a subquotient
of an object of the form $\pR(\psh)$, where $\psh$ runs through all objects
of $\epcat$, the homomorphism is an embedding.
%

Via the Tannakian dictionary, 
we have seen that the category $\catq$ is equivalent as a tensor category
to the category $\Rep(\ch H)$ of finite-dimensional representations
of a group-scheme $\ch H$ so that the hypercohomology $\hc:\catq\to\vect$
corresponds
to the forgetful functor.
From the commutative diagram of Section~\ref{sothercats}, we obtain
a commutative diagram of groups 
$$
\begin{array}{ccccc}
& & \ch H_\pR & & \\
& & \uparrow & \nwarrow  & \\
\ch H & {\leftarrow} & \ch H\times\C^\times & \stackrel{}{\leftarrow} & \ch H \\
\end{array}
$$
where the horizontal arrows are the obvious projection and inclusion,
the vertical arrow is a surjection, and the diagonal arrow is an inclusion.

%

\end{subsection}
%

\begin{subsection}{Embedding into Levi subgroup}\label{sslevi}


%

%


Recall that $2\ch\rho_M$ is the sum of the positive roots of the
Levi factor $M\subset P$ of
the complexification $P\subset G$ of the minimal parabolic
subgroup $P_\R\subset G_\R$.
We define the Levi subgroup $\ch L_0\subset\ch G$ to be the centralizer
of $2\ch\rho_M$ under the adjoint representation.

\begin{Proposition}
The embedding $\ch H_\pR\to \ch G$ factors through $\ch L_0$.
\end{Proposition}

\begin{proof}
The composite functor 
$$
\epcat\stackrel{\pR}{\to}\catq_\pR\stackrel{\ef}{\to} \catq_\Z\stackrel{\ff}{\to} \vect_\Z
$$
provides a composite homomorphism
$$
\langle 1\rangle\times\C^\times\to \ch H\times\C^\times\to \ch H_\pR\to\ch G
$$
in which 
the homomorphism
$
\ch H\times\C^\times\to \ch H_\pR
$
is surjective. 

The image of the composite homomorphism
is the subgroup
$\langle\exp(2c\ch\rho_M)|c\in\C\rangle\subset\ch G$ since the isomorphism
$$
\schar\circ\pR\simeq\Sigma\circ\tchar:\epcat\to\vect_{\rcowt\times\Z}
$$
implies that the $\Z$-grading on the perverse specialization 
$\pR$ corresponds to the $\Z$-grading on $\tchar$
which was defined by pairing with $2\ch\rho_M$.
\end{proof}

%

\end{subsection}


\begin{subsection}{Maximal torus}\label{ssmaxtorus}

We have seen that the real character functor 
$\schar:\catq_\pR\to\vect_\bigimlat$
is a tensor functor, and there is a canonical isomorphism of tensor functors
$$
\hc\simeq\ff\circ\schar:\catq_\pR\to\vect
$$
where $\ff:\vect_\bigimlat\to\vect$ is the forgetful functor.
Therefore we may identify the group $\ch H_\pR$ with the group
of tensor automorphisms of the composite functor $\ff\circ\schar$,
and we have a homomorphism
$$
\ch S_\bigimlat\to\ch H_\pR
$$
where $\ch S_\bigimlat$ is the torus 
$$
\ch S_\bigimlat=\Spec(\C[\bigimlat])
$$
where $\bigimlat$ is the image of the projection $\Sigma:\cowt\to\rcowt\times\Z$ 
defined by
$\Sigma(\lambda)=(\lambda+\theta(\lambda), \langle 2\ch\rho_M,\lambda\rangle)$.

We have also seen that there is a canonical isomorphism
of tensor functors
$$
\schar\circ\pR\simeq\Sigma\circ\tchar:\epcat\to\vect_\bigimlat.
$$
Therefore we have a commutative diagram
$$
\begin{array}{ccc}
\ch H_\pR&\to&\ch G\\
\uparrow && \uparrow \\
\ch S_\bigimlat &\to &\ch T
\end{array}
$$
where the embedding 
$
\ch S_\bigimlat\to\ch T
$
corresponds to the tensor functor $\Sigma:\vect_\cowt\to\vect_\bigimlat$. 
Since the homomorphism $\ch S_\bigimlat\to\ch T\to \ch G$ is an embedding,
we conclude that the homomorphism
$
\ch S_\bigimlat\to\ch H_\pR
$
is an embedding.

Our aim is to show this is the embedding of a maximal torus.

\begin{Lemma}\label{lsimple}
If $\CQ$ is a simple object in the category $\catq_\pR$, then it is 
isomorphic to the shifted intersection cohomology
sheaf $\ic^{\lambda}[z]$ of the closure of a stratum $\raffgr^\lambda$,
for a pair $(\lambda,z)\in\bigimlat$, 
with coefficients in a one-dimensional trivial local system.
\end{Lemma}

\begin{proof}
If $\CQ$ is a simple object in the category $\catq_\pR$,
then it must be the shifted intersection cohomology 
sheaf $\ic^\lambda(\ls)[z]$ of the closure of a stratum $\raffgr^\lambda$ with coefficients
in an irreducible $G_\R(\cor)$-equivarant local systems $\ls$.
Thanks to the isomorphism
$$
F^{\lambda}_\R(\ic^\lambda(\ls))\simeq\ls|\lambda,
$$
the real character functor
$\schar$ 
does not vanish on $\ic^\lambda(\ls)$. By 
Corollary~\ref{ccharim}, we conclude that $(\lambda,z)\in\bigimlat$.

To prove the
local system $\ls$ is trivial,
it suffices by Lemma~\ref{lffeq} 
to show that the stabilizer $P_\R^\lambda\subset G_\R$
acts trivially on the stalk $\ls|\lambda$.
As a representation of $P_\R^\lambda$, the stalk $\ls|\lambda$ 
is isomorphic
to the weight space $F^{\lambda}_\R(\ic^\lambda(\ls))$.
The weight space $F^{\lambda}_\R(\ic^\lambda(\ls))$
is a direct summand in the hypercohomology
$\hc(\raffgr,\ic^\lambda_\R(\ls))$.
The hypercohomology $\hc(\raffgr,\ic^\lambda_\R(\ls))$
is a subquotient of the hypercohomology
$\hc(\raffgr,\pR(\psh))$, for some object $\psh$ in the category $\epcat$.
The hypercohomology $\hc(\raffgr,\pR(\psh))$ is isomorphic to the hypercohomology
$\hc(\affgr,\psh)$.
Since $G$ is connected, it acts trivially on $\hc(\affgr,\psh)$,
and so the subgroup $P_\R^\lambda$ does as well.
%
\end{proof}

\begin{Proposition} \label{pcmaxt}
The group $H_\pR$ is connected, and the embedding $\ch S_\bigimlat\to\ch H_\pR$
is that of a maximal torus.
\end{Proposition}

\begin{proof}
For any dominant real coweight $\lambda\in\rdcowt$,
the support of the $n$-fold convolution of the intersection cohomology sheaf ${\ic^\lambda}$ of the closure of the stratum $\raffgr^\lambda$ contains the coweight $n\lambda$,
and so by the lemma and the criterion of \cite[Corollary 2.22]{DM82},
the group $\ch H_\pR$ is connected.

The lemma also implies that the rank of $\ch H_\pR$ is less than or equal to
the rank of the embedded torus $\ch S_\bigimlat$. Thus they are the same,
and $\ch S_\bigimlat$ is a maximal torus.
\end{proof}

\end{subsection}


\begin{subsection}{Root system}\label{ssrootsystem}

In this section, we identify the root system of
the subgroup $\ch H\subset\ch G$.
Note that due to the surjective homomorphism 
$
\ch H\times\C^\times\to\ch H_\pR,
$
the root system of $\ch H$ is equal to that of $\ch H_\pR$.

By forgetting the $\Z$-grading, 
from the embedding $\ch S_\bigimlat\to\ch H_\pR$,
we obtain a commutative diagram of embeddings
$$
\begin{array}{ccc}
\ch S_\bigimlat & \to & \ch H_\pR\\
\uparrow & & \uparrow \\
\ch S_\imlat & \to & \ch H 
\end{array}
$$
where the torus $\ch S_\imlat$ is defined to be
$$
\ch S_\imlat=\Spec(\C[\imlat])
$$
where $\sigma(\cowt)$ is the image of the projection $\sigma:\cowt\to\rcowt$ defined by
$\sigma(\lambda)=\lambda+\theta(\lambda)$.

\begin{Proposition}
The group $\ch H$ is connected, and the embedding $\ch S_\imlat\to\ch H_\pR$
is that of a maximal torus.
The rank of $\ch H$ equals the split-rank of $G_\R$.
\end{Proposition}

\begin{proof}
The first two statements are proved similarly to those of Proposition~\ref{pcmaxt}.
To see the last, 
note that the lattice $\imlat$ is a finite index sublattice in the real coweight lattice $\rcowt$,
and 
$S\subset G$ is the complexification of a maximal split torus $S_\R\subset G_\R$.
\end{proof}

By construction, the projection $\sigma:\cowt\to\imlat$ 
is the restriction of weights induced by the inclusion $\ch S_\imlat\to\ch T$.
Therefore the roots of $\ch H$ with respect to the 
maximal torus $\ch S_\imlat$ are a subset of the
projection $\sigma(R_0)$ of the roots
$R_0$ 
of the Levi subgroup $\ch L_0 $ with respect to the 
maximal torus $\ch T$.
Thanks to the identity $\ch\theta(2\ch\rho_M)=-2\ch\rho_M$
for the dual involution $\ch\theta$, 
the involution $\theta$ preserves the set of simple roots
$\Delta_0$
of the Levi subgroup $\ch L_0 $ with respect to the Borel subgroup
$\ch B\cap \ch L_0$. 
Therefore the projection $\sigma(R_0)$ 
forms a possibly non-reduced root system in $\imlat$ in which
the multiplicity of a projected root is one or two.
(A non-reduced root system is one in which there is a root $\alpha$
such that $2\alpha$ is also a root.
The multiplicity of a projected
root is the number
of roots which project to it.)
The classification in \cite{Ar62} immediately implies the following.


\begin{Lemma}\label{laraki}
The root system $\sigma(R_0)$ is non-reduced
if and only if the based root system $\Delta_0$ contains
a copy $\atwo$ of the based root system of $\liesl_3$ such
that the involution $\theta$ exchanges the simple roots of $\atwo$.
To be precise,
for a projected root $\alpha\in\sigma(R_0)$,
twice it is also a projected root $2\alpha\in\sigma(R_0)$ if and only if 
there are roots $\beta_1,\beta_2\in R_0$ such that $\alpha=\sigma(\beta_1)=\sigma(\beta_2)$ 
and the roots $\beta_1, \beta_2$ 
are the simple roots of a copy of the root system of $\liesl_3$
in the root system $R_0$.
\end{Lemma}

\begin{Remark}
By the classification, 
for $\lieg$ simple, the root system $\sigma(R_0)$ is non-reduced
if and only if $\lieg=\EuFrak{sl}_{n}(\C)$ and
$\lieg_\R=\EuFrak{su}(p,q)$ with $p+q=n$ odd.
\end{Remark}

Now the obvious candidate for the root system of $\ch H$ with respect to the maximal torus 
$\ch S_\imlat$ would be the projected root system
$\sigma(R_0)\subset\imlat$,
but we have seen that in general it is not reduced.
We define the root system
$$
\Xi\subset \sigma(\cowt)
$$
to be all of the projected roots $\alpha\in\sigma(R_0)$ satisfying $2\alpha\not\in\sigma(R_0).$
%
In other words, if $\sigma(R_0)$ is reduced,
then we take $\Xi$ to be all of the projected roots $\sigma(R_0)$, 
but if $\sigma(R_0)$ is non-reduced,
then we take $\Xi$ to be all of the projected roots $\sigma(R_0)$
except we discard the shorter root of any pair along a single ray.

\begin{Theorem} The group
$\ch H$ is reductive with 
root system $\Xi\subset \imlat$.
\end{Theorem}

\begin{proof}
We first identify the 
root system of the maximal reductive quotient $\ch H_\red$
with respect to the image 
of the maximal torus
$\ch S_\imlat$ under the natural projection $\ch H\to\ch H_\red$.

\begin{Proposition}\label{pweyl}
The root system of $\ch H_\red$ contains a root in every direction
in which the root system $\Xi$ contains a root.
\end{Proposition}

\begin{proof}
By definition, the root directions of the root systems $\Xi$ and $\sigma(R_0)$ coincide.
Since the roots of $\ch H_\red$ are a subset of $\sigma(R_0)$,
the Weyl group of $\ch H_\red$ is a subgroup of the Weyl group of
$\sigma(R_0)$.
By Proposition~\ref{pweylsym}, the Weyl group of $\ch H_\red$ is isomorphic to the small Weyl group 
$\CW_{G_\R}$ so that their actions on $\imlat$ agree. 
It is a standard fact that $\CW_{G_\R}$ is the Weyl group of the root system $\sigma(R)$,
and so
the Weyl group of $\sigma(R_0)$ is a subgroup of $\CW_{G_\R}$.
We conclude that the Weyl group of $\ch H_\red$, the Weyl group of $\sigma(R_0)$,
and the small Weyl group $\CW_{G_\R}$ must all coincide.
Since the root directions of a root system 
are the negative eigenspaces of the simple
reflections of its Weyl group, the assertion follows.
\end{proof}

Since the roots of $\ch H_\red$ are a subset of $\sigma(R_0)$,
if $\sigma(R_0)$ is reduced, so that $\Xi$ coincides with it, then
the proposition immediately implies that $\Xi$ is the root system of $\ch H_\red$.


If $\sigma(R_0)$ is non-reduced, then it contains roots $\alpha$ and $2\alpha$, and
by the proposition, the root system of $H_{\on{red}}$ contains either $\alpha$ or $2\alpha$.
By definition, the root system $\Xi$ contains $2\alpha$ not $\alpha$.
To arrive at a contradiction, suppose the root system of $\ch H_{\on{red}}$ 
contains $\alpha$ not $2\alpha$.
Then the zero weight space of
any representation of $\ch H_{\on{red}}$ with non-empty $\alpha$ weight space must also be
non-empty.
In particular, the irreducible representation of highest weight $\alpha_0$
a Weyl group translate of $\alpha$ must have non-empty zero weight space.
It must correspond to the intersection cohomology sheaf $\ic^{\alpha_0}$
of the stratum $\raffgr^{\alpha_0}$

\begin{Lemma}
Assume the root system $\sigma(R_0)$ 
is non-reduced and contains
both $\alpha$ and $2\alpha$. 
Then the intersection of the orbit $S^0_{\R}$ through the zero coweight 
with the stratum $\raffgr^{\alpha_0}$ satisfies
$$
\dim_\R(S^0_{\R}\cap\raffgr^{\alpha_0})<\langle\ch\rho,\alpha_0\rangle.
$$
\end{Lemma}

\begin{proof} 
It suffices to 
show that none of the components of the intersection $S^0\cap \affgr^{\alpha_0}$ 
is preserved by the conjugation $\theta$ of $\affgr$. 
By the explicit description provided by Lemma~\ref{laraki}, 
the coweight $\alpha_0$ is in fact a coroot of $G$,
and the intersection cohomology sheaf $\ic^{\alpha_0}$ of the stratum $\affgr^{\alpha_0}$
corresponds to a summand in the adjoint representation of $\LG$.
It is easy to check that each irreducible component of the intersection $S^0\cap \affgr^{\alpha_0}$
contains a unique simple coroot of $G$ in its closure.
But by the explicit description
provided by Lemma~\ref{laraki},
none of the simple coroots in the stratum $\affgr^{\alpha_0}$
are fixed by the conjugation $\theta$.
\end{proof}
   
By Proposition~\ref{pextravanish},
the real weight functor $F^0_\R$ for the zero weight space 
vanishes on $\ic^{\alpha_0}$, and we arrive at a contradiction. 
Therefore by Proposition~\ref{pweyl},
the root system of $\ch H_{\red}$ contains the root $2\alpha$ not $\alpha$,
and thus coincides with $\Xi$.

It remains to show that $\ch H$ is reductive,
or in other words, that its unipotent radical $\ch H_u$ is trivial.
Since $\ch H_u$ is connected, the following suffices.

\begin{Lemma}
The Lie algebra $\ch\fh_u$ is trival.
\end{Lemma}

\begin{proof}
Consider the Lie algebra $\ch l_0$ of the Levi subgroup $\ch L_0$
as a representation of $\ch H$.
It contains $\ch \fh$ as a subrepresentation which in turn 
contains $\ch \fh_u$ as a subrepresentation.
In what follows, references to the weights of these representations are with respect to the
torus $\ch S_\imlat$.

The zero weight space of $\ch \fl_0$ 
is the Cartan subalgebra $\ch \ft$. Since the intersection 
$\ch \fh_u\cap\ch \ft$ is trivial, 
the zero weight space of $\ch \fh_u$ is empty.
Therefore the weights of $\ch \fh_u$ 
lie in the projection $\sigma(R_0)$ of the roots of $\ch L_0$.
Furthermore, if $\alpha$ 
is a weight of $\ch \fh_u$, 
then $\alpha$ is not a root of $\ch H_{\red}$ since 
otherwise the zero weight space of $\ch \fh_u$ would be non-empty.
In summary, we see that the weights of $\ch \fh_u$
lie in the complement of the roots of $\ch H_{\red}$ 
within the projection $\sigma(R_0)$ of the roots of $\ch L_0$. 
Therefore if the root system $\sigma(R_0)$ is reduced, 
so that $\Xi$ coincides with it, then 
we are done.
 
It remains to consider the case when 
the root system $\sigma(R_0)$ is non-reduced. 
We may assume that $G$ is simple and adjoint,
so that by the classification we have $G=\on{PGL}_n\C$ and $G_\R=\on{PU}(p,q)$
with $p+q=n$ odd.
The standard representation of $\LG=\on{SL}_n\C$ corresponds to
the intersection cohomology sheaf of $\Bbb P(\C^n)$.
The specialization of this stratum is topologically
equivalent to the map that collapses the disjoint union
of the linear
subvarieties $\Bbb P(\C^p)$ and $\Bbb P(\C^q)$ to a point.
Thus the specialization in this case is clearly semisimple.
But if $\ch \fh_u$ were not trivial, then the restriction
of this representation to $\ch H$ would not be a
direct sum of irreducibles.
\end{proof}

This completes the proof of the theorem.
\end{proof}

\begin{Corollary}
The Weyl group of $\ch H$ is isomorphic to the small Weyl group of $G_\R$. 
\end{Corollary}

\begin{proof} 
Follows from Proposition~\ref{pweylsym} as was noted in 
the proof of Proposition~\ref{pweyl}.
\end{proof}

\begin{Corollary}
The category $\catq$ is semisimple.
Its object are sheaves isomorphic to a direct sum of 
intersection cohomology sheaves $\ic^\lambda$ of the closures of strata 
$\raffgr^\lambda$, for $\lambda\in\imlat$,
with coefficients in trivial local systems.
In particular, convolution and perverse specialization are semisimple.
\end{Corollary}

\end{subsection}


\begin{subsection}{Fixed points of automorphisms}\label{sfpofauto}
In this section,
we explain how the subgroup $\ch H\subset\ch G$ 
may be realized as the identity component of the fixed points of an involution 
of a certain Levi subgroup $\ch L_1\subset\ch G$.
To make the discussion self-contained, we begin
by collecting notation and results from previous sections.

We have fixed a minimal parabolic subgroup $P_\R\subset G_\R$,
with Levi factor $M_\R\subset P_\R$, 
maximal torus $T_\R\subset M_\R$,
and maximal split torus $S_\R\subset T_\R$.
We write $P$ for the complexification of $P_\R$,
$M$ for that of $M_\R$, $T$ for that of $T_\R$,
and $S$ for that of $S_\R$.

In Section~\ref{ssnot},
we have seen that
the conjugation of $G$ with respect to $G_\R$ induces an involution $\theta$
of the coweight lattice $\cowt$.
The real coweight lattice $\rcowt$
is precisely the fixed points of the involution $\theta$,
and we have a projection $\sigma:\cowt\to\rcowt$ defined by 
$\sigma(\lambda)=\theta(\lambda)+\lambda$.

In Section~\ref{ssdualgroup}, we have seen that the dual group $\ch G$ comes equipped 
with a distinguished Borel
subgroup $\ch B\subset\ch G$, and maximal torus $\ch T\subset \ch G$
which is identified with the torus dual to $T$.
We write $\Delta_{\ch B,\ch T}\subset\cowt$ for the simple roots of $\ch G$ with respect to $\ch B$
and $\ch T$.

In Section~\ref{sslevi},
we have seen that the subgroup $\ch H\subset\ch G$ is in fact a subgroup of
the Levi subgroup $\ch L_0\subset\ch G$ which is the centralizer
of $2\ch\rho_M$ under the adjoint representation.
Here as usual $2\ch\rho_M$ is the sum of the positive roots of
the Levi subgroup
$M\subset G$.

The simple roots $\Delta_0\subset\Delta_{\ch B,\ch T}$ 
of the Levi subgroup $\ch L_0\subset\ch G$
with respect to the Borel subgroup $\ch B\cap \ch L_0$,
and maximal torus $\ch T$, satisfy 
$$
\alpha\in\Delta_0\mbox{ if and only if }\langle2\ch\rho_M,\alpha\rangle=0.
$$
Note that we have  
$\theta(\Delta_0)=\Delta_0$
thanks to the identity $\ch\theta(2\ch\rho_M)=-2\ch\rho_M$,
where $\ch\theta$ is the dual involution. 

To define the Levi subgroup $\ch L_1\subset\ch G$,
for each simple root 
$
\alpha\in\Delta_{0},
$
choose 
a basis vector $X_\alpha$ in the simple root space 
$
(\fb\cap\fl_0)_\alpha.
$
Define elements $\tau(\alpha)\in\cowt$ to be 
$$
\tau(\alpha)=\alpha-\theta(\alpha), 
\mbox{ if } 
[X_\alpha,X_{\theta(\alpha)}]\neq 0
$$
$$
\tau(\alpha)=0,
\mbox{ if } 
[X_\alpha,X_{\theta(\alpha)}] = 0.
$$
Clearly the elements $\tau(\alpha)$ are independent of the choice
of basis vectors $X_\alpha$.

Define the Levi subgroup $\ch L_1\subset \ch G$ to be
the centralizer of the span of $2\rho_M$ and the elements $\tau(\alpha)$,
for all simple roots $\alpha\in\Delta_0$, under the coadjoint representation. 
Here as usual $2\rho_M$ is the sum of the positive coroots of
the Levi subgroup
$M\subset G$.

The simple roots $\Delta_1\subset\Delta_{\ch B,\ch T}$
of the Levi subgroup $\ch L_1\subset\ch G$ with respect to the Borel subgroup $\ch B\cap \ch L_1$,
and maximal torus $\ch T$, are of two types
$$
\alpha\in\Delta_1, \mbox{ for }
\alpha\in\Delta_0 \mbox{ with }
[X_\alpha,X_{\theta(\alpha)}] = 0,
\mbox{ and }
$$
$$
\alpha+\theta(\alpha)\in\Delta_1, \mbox{ for }
\alpha\in\Delta_0 \mbox{ with }
[X_\alpha,X_{\theta(\alpha)}] \neq 0.
$$
Note that we have a similar description for the simple coroots $\ch\Delta_1$.

The involution $\theta$ 
preserves the set $\Delta_1$,
and similarly the dual involution $\ch\theta$ preserves the set $\ch\Delta_1$.
Therefore they provide an involution of the based root datum
$$
\Psi(\ch L_1,\ch B\cap\ch L_1,\ch T)=(\cowt,\Delta_1,\wt,\ch\Delta_1)
$$
of the Levi subgroup $\ch L_1\subset \ch G$.
To lift the involution 
to an involution of $\ch L_1$,
we must choose for each simple root 
$
\beta\in\Delta_{1},
$ 
a basis vector $X_{\beta}$
in the corresponding simple root space 
$
(\ch \fb\cap \ch \fl_1)_\beta.
$
Then there is a unique lift of the involution 
to an involution of $\ch L_1$ which we denote by $\ch\theta_1$ such that
$$
\ch\theta_1(\ch B\cap\ch L_1)=\ch B\cap\ch L_1\qquad \ch\theta_1(\ch T)=\ch T\qquad
\ch\theta_1(X_\beta)=X_{\theta(\beta)}.
$$
For any other choice of basis vectors, the resulting lift will be conjugate
to $\ch\theta_1$ via an inner automorphism. Moreover, this inner automorphism
may be realized as conjugation by a unique element of $\ch T/ Z(\ch L_1)$,
where $Z(\ch L_1)$ is the center of $\ch L_1$.

\begin{Remark}\label{crit}
We mention here a criterion for $\ch L_1$ to equal $\ch L_0$.
For all simple roots $\alpha\in\Delta_{0}$,
we have
$$
[X_\alpha,X_{\theta(\alpha)}]= 0
$$ 
if and only if for each connected component $C$ in the Dynkin diagram
of $\ch L_0$ with $\ch \theta_0(C)=C$, there is a node $c\in C$
such that $\theta(c)=c$.
\end{Remark}

Alternatively, one can realize the Levi subgroup $\ch L_1$ as the fixed
points of an inner
automorphism $\ch\eta_0$ of the Levi subgroup $\ch L_0$.
For each simple root $\alpha\in\Delta_0$, consider the new basis vectors $X'_\alpha$ defined by
$$
X'_\alpha=X_\alpha, \mbox{ if } 
[X_\alpha,X_{\theta(\alpha)}]=0,
$$
$$
X'_\alpha=k_\alpha X_\alpha, \mbox{ if } 
[X_\alpha,X_{\theta(\alpha)}]\neq 0,
$$
for any non-zero numbers $k_\alpha$. 
Then there is a unique inner automorphism $\ch\eta_0$ of the group $\ch L_0$ 
such that
$$
\ch\eta_0(\ch B\cap\ch L_0)=\ch B\cap\ch L_0\qquad \ch\eta_0(\ch T)=\ch T\qquad
\ch\eta_0(X_\alpha)=X'_\alpha,
$$
and it does not depend on the choice of the basis vectors $X_\alpha$.
If one takes the numbers $k_\alpha$ to be different from $1$, then 
the Levi subgroup $\ch L_1$ is the fixed points in the Levi subgroup
$\ch L_0$ of the inner automorphism
$\ch \eta_0$. 

\begin{Theorem}
There is a choice of basis for the simple root spaces of the Levi
subgroup $\ch L_1\subset \ch G$ such that
the subgroup $\ch H\subset\ch G$ associated
to the real form $G_\R$
is the identity component of the fixed points 
of the involution $\ch\theta_1$.
\end{Theorem}
%

\begin{proof}
In Section~\ref{ssmaxtorus}, we have seen that the torus $\ch S_\imlat=\Spec(\C[\imlat])$ is a maximal
torus of $\ch H$.
Here as usual $\imlat$ denotes
the image of the projection $\sigma:\cowt\to\rcowt$.

\begin{Lemma}
The identity component of the fixed points of the
restriction of the involution $\ch\theta_1$ to the torus $\ch T$
is equal to the torus $\ch S_\imlat$.
\end{Lemma}

\begin{proof}
By definition, the restriction of the involution $\ch\theta_1$ to the torus $\ch T$
is induced by the involution $\theta$ of the lattice $\cowt$. 
\end{proof}

We write $R_0\subset\cowt$ for the roots of $\ch L_0$ with respect
to $\ch T$, and 
$R_1\subset \cowt$ for the roots of $\ch L_1$ with respect to $\ch T$.
In Section~\ref{ssrootsystem},
we have seen that the roots $\Xi\subset\imlat$ of the group $\ch H$
with respect to $\ch S_\imlat$
consist of those projected roots $\sigma(\alpha)\in\sigma(R_0)$ such that
$2\sigma(\alpha)\not\in\sigma(R_0)$.

\begin{Lemma}
For a root $\alpha\in R_0$, 
we have $\sigma(\alpha)\in\Xi$ if and only if $\alpha\in R_1$.
Therefore we have the inclusion $\ch H\subset \ch L_1$.
\end{Lemma}

\begin{proof}
Follows from the definition of $\ch L_1$, and Lemma~\ref{laraki}.
\end{proof}

It remains to show that there is a choice of basis vectors $X_\beta$ for
the simple root spaces of $\ch L_1$
such that the simple root spaces of $\ch H$ are the fixed subspaces
of the involution $\ch\theta_1$.
For a simple root in $\Xi$, 
there are two possibilities: one simple root $\beta\in \Delta_1$ projects to it,
or two simple roots $\beta_1,\beta_2\in\Delta_1$ project to it.
If only one simple
root $\beta$ projects to a root in $\Xi$, we take $X_\beta$ to be any non-zero
vector in the corresponding root space.

\begin{Lemma}
If two simple roots $\beta_1,\beta_2$ project to a root in $\Xi$, then 
the intersection of $\ch\fh$
with the span of the corresponding simple root spaces is exactly one-dimensional, and the
projection of the intersection
to either of the two root spaces is an isomorphism. 
\end{Lemma}

\begin{proof}
In Section~\ref{ssrootsystem}, we have seen that the intersection is at least
one-dimensional. Since $\ch S_\imlat$ is a maximal torus of $\ch H$,
the intersection is no greater than one-dimensional,
and can not lie entirely in either simple root space.
\end{proof}

If two simple roots $\beta_1,\beta_2$ project to a root in $\Xi$,
we choose any non-zero
vector in the intersection, and take $X_{\beta_1},
X_{\beta_2}$ to be its images under the projections.

It is straightforward to check that the subalgebra
$\ch\fh$ is exactly
the fixed points in $\ch\fl_1$ of the involution $\ch\theta_1$
constructed with respect to this basis.
\end{proof}

Alternatively,
one can define an automorphism of the Levi subgroup $\ch L_0\subset\ch G$ 
such that the subgroup $\ch H\subset\ch G$ 
is the identity component of its fixed points. 
In the definition of the inner automorphism $\ch\eta_0$, 
for each pair of simple roots $\alpha, \theta(\alpha)\in\Delta_0$
such that
$
[X_\alpha,X_{\theta(\alpha)}]\neq 0,
$
if one takes the numbers
$k_\alpha, k_{\theta(\alpha)}$ to be any pair such that $k_\alpha k_{\theta(\alpha)}=-1$,
%
then the subgroup $\ch H\subset\ch G$
is the identity component of the fixed points of the automorphism
$\ch\theta_0\circ\ch\eta_0$ of the Levi subgroup $\ch L_0$.
For example, for the choice $k_\alpha=1$ and $k_{\theta(\alpha)}=-1$, the automorphism
$\ch\theta_0\circ\ch\eta_0$ 
is of order $4$.




\end{subsection}

\end{section}

%



\begin{section}{Appendix: Identification of the dual group}
Let $G$ be a connected reductive complex algebraic group
with fixed Borel subgroup $B\subset G$, and maximal torus $T\subset B$.
We sketch here how to show the Tannakian group $\tangp$ of the category
$\epcat$ with fiber functor $\hc:\epcat\to\vect$ is a dual group for $G$.
We must show that $\tangp$ 
is a connected reductive complex algebraic group, and find a canonical
isomorphism of based root data 
$$
\Psi(\tangp)\simeq\ch\Psi(G).
$$  
In what follows, for $\lambda\in\dcowt$, we write $\ic^\lambda$ for the 
intersection cohomology  sheaf of the closure
of the $G(\co)$-orbit $\affgr^\lambda\subset\affgr$, with coefficients
in the trivial one-dimensional local system. 

\begin{Proposition}\label{appsemisimple}
The category $\epcat$ is semisimple, with simple objects 
isomorphic to $\ic^\lambda$, for $\lambda\in\dcowt$.
\end{Proposition}

\begin{proof}
The stabilizer in $G(\co)$ of a coweight $\lambda\in\affgr$
is the parahoric subgroup $\EuScript P^\lambda$ which is connected.
Therefore the simple objects of the category $\epcat$ are as asserted,
and there are no self-extensions of simple objects.
By Proposition~\ref{plgstrat}, the strata $\affgr^\lambda$ in
a given component of $\affgr$ 
are either all even-dimensional or all odd-dimensional.
By~\cite[Section 11, c)]{Lu83}, the stalks of  $\ic^\lambda$
have the parity vanishing property: they are non-zero only in
the parity of the dimension of $\affgr^\lambda$.
Therefore there are no other extensions, and we conclude that
the category $\epcat$ is semisimple.
\end{proof}

\begin{Corollary}\label{appconnectedreductivealgebraic}
The Tannakian group $\tangp$ is a connected reductive complex algebraic group.
\end{Corollary}

\begin{proof}
Choose generators $\lambda_1,\ldots,\lambda_r\in\dcowt$,
and 
for $\lambda\in\dcowt,$ write $\lambda=\sum_i n_i\lambda_i $.
Then the open part of the support of the convolution 
$({\ic^\lambda_1})^{\odot n_1}\odot\cdots\odot({\ic^\lambda_r})^{\odot n_r}$
contains $\lambda$,
and so $\ic^\lambda$ appears as a summand in the convolution.  
Therefore by the proposition, the direct sum $\oplus_i \ic^{\lambda_i}$
is a tensor generator for the category $\epcat$,
and so by~\cite[Proposition 2.20(b)]{DM82}, the group $\tangp$ is algebraic.

For $\lambda\in\dcowt$,
the open part of the support of the convolution $({\ic^\lambda})^{\otimes n}$  
contains $n\lambda$,
and so $\ic^{n\lambda}$ appears as a summand in the convolution.  
Therefore by the proposition and \cite[Corollary 2.22]{DM82},
the group $\tangp$ is connected.

Finally, by the proposition and~\cite[Proposition 2.23]{DM82}, the group $\tangp$ is reductive.
\end{proof}

The character functor $\tchar:\epcat\to\vect_\cowt$,
which is by definition the sum
$\tchar= \sum_{\lambda\in\cowt}F^\lambda$ of the weight functors
$F^\lambda:\epcat\to\vect$,
and the canonical isomorphism 
$
\hc\simeq\ff\circ\tchar:\epcat\to\vect_\cowt,
$
where $\ff:\vect_\cowt\to\vect$ is the forgetful functor,
provide a canonical homomorphism
from the dual torus 
$
\ch T=\Aut^\otimes(\ff)
$
to $\tangp.$

\begin{Proposition}\label{appmaximaltorus}
The homomorphism $\ch T\to\tangp$ is the embedding of a maximal torus.
\end{Proposition}

\begin{proof}
For $\lambda\in\dcowt$, we clearly have
$F^\lambda(\ic^\lambda)\simeq\C$,
so by~\cite[Proposition 2.21 (b)]{DM82}, the homomorphism is injective.
The rank of $\tangp$ equals the rank of $\epcat$,
and by the previous proposition, the rank of $\epcat$ equals the rank
of $G$. This in turn equals the rank of $\ch T$,
and so the image
of the homomorphism
is a maximal torus.
\end{proof}

A Borel subgroup $\ch B\subset\tangp$ containing $\ch T$ is equivalent to the choice
of a $\ch T$-invariant line $L^\lambda$
in each irreducible representation $V^\lambda$, for $\lambda\in\dcowt$
such that $L^\lambda\otimes L^\mu=L^{\lambda+\mu}$ under the projection
$V^\lambda\otimes V^\mu\to V^{\lambda+\mu}$. 
The Borel subgroup $\ch B\subset\tangp$
is the stabilizer of the lines, or conversely the lines are the
highest weight lines of the Borel subgroup.
Using the canonical isomorphism $\hc\simeq\ff\circ\tchar$,
we take $L^\lambda=F^\lambda(\ic^\lambda)$, for $\lambda\in\dcowt$.
 
It remains to show there is a canonical isomorphism of based root data
$$
\Psi(\tangp,\ch B,\ch T)\simeq\ch\Psi(G,B,T).
$$  

\begin{Proposition}\label{approotdatum}
With respect to the given Borel subgroups
and maximal tori, (i) the set of dominant weights
of $\tangp$ coincides with the set of dominant coweights of $G$
as subsets of $\cowt$, and (ii) the set of simple roots of $\tangp$
coincides with the set of simple coroots of $G$ as subsets of $\cowt$.
\end{Proposition}

\begin{proof}
Assertion (i) is immediate: 
the weights of the lines $L^\lambda\subset\hc(\affgr,\ic^\lambda)$
coincide with the dominant coweights of $G$.

For assertion (ii), it suffices to show that for $\lambda\in\dcowt$, and $\mu\in\cowt$,
the weights $F^\mu(\ic^\lambda)$ vanish if $\lambda-\mu$ is not 
a non-negative integral
linear combination of positive coroots $\alpha\in R^\pos$, 
and for a simple coroot $\alpha\in\Delta_{B,T}$,
the weights $F^{\lambda-\alpha}(\ic^\lambda)$ do not vanish.
Proposition~\ref{pmv} implies the vanishing, and Propositions~\ref{pmv} and~\ref{pdimest},
and
the explicit description of the weights given in \cite[Section 5]{MV00} implies
the non-vanishing.
\end{proof}

Now we have the asserted isomorphism of based root data.
The simple coroot directions of $\tangp$
are the extremal rays of the positive semigroup dual of 
the dominant weights of $\tangp$. Similarly, the simple root directions of $G$ are 
the extremal rays of the positive semigroup dual of 
the dominant coweights of $G$. By part (i) of the proposition,
these rays coincide. For a root $\ch\alpha$ and the corresponding coroot $\alpha$,
one has $\langle\ch\alpha,\alpha\rangle=2$. 
Therefore the simple coroots
of $\tangp$ coincide
with the simple roots of $G$
since they are in
the same rays and since, by part (ii) of the proposition, the simple roots of $\tangp$ coincide
with the simple coroots of $G$. 

%


%
\begin{Remark}
For Levi subgroups $L\subset G$ other than the torus $T$, 
one may generalize the character functor $\tchar:\epcat\to\vect_\cowt$
to tensor functors $\on{Ch}_{L}:\epcat\to\catp_{L(\co)}(\affgr_L)$. 
In particular, for each simple coroot of $G$,
this provides an explicit embedding of an $\on{\mathfrak{sl}_2}$-triple into
the Lie algebra of $\tangp$.
For more details, see for example~\cite[Section 4.3.1]{BG02}.
This is used in~\cite[Theorem 2.2]{BG01} to construct crystals for $\ch\fg$.
\end{Remark}

\end{section}


\bibliographystyle{alpha}
\bibliography{ref}


\end{document}